\documentclass[twoside, 11pt]{article}

\usepackage{jmlr2e} 
\usepackage[utf8]{inputenc} 
\usepackage[T1]{fontenc}
\usepackage{url}            
\usepackage{booktabs}       
\usepackage{amsfonts}       
\usepackage{nicefrac}       
\usepackage{microtype}      
\usepackage{thmtools,thm-restate}

\usepackage{enumitem} 

\usepackage{color}

\usepackage{subcaption}
\usepackage{float}

\usepackage{algorithm}
\usepackage{algorithmic}

\usepackage{mdframed}

\usepackage{amsmath}
\usepackage{bbm, bm}

\usepackage{hyperref}
\hypersetup{ %
	pdftitle={},
	pdfauthor={R\'emi Leblond, Fabian Pedregosa, Simon Lacoste-Julien},
	pdfkeywords={},
	pdfborder=0 0 0,
	pdfpagemode=UseNone,
	colorlinks=true,
	linkcolor=blue, 
	citecolor=blue, 
	filecolor=blue, 
	urlcolor=blue, 
	pdfview=FitH,
	hypertexnames=false,
}

\newcommand{\stepsize}{\gamma}
\newcommand{\strongconvex}{\mu}
\newcommand{\overlap}{\tau}
\newcommand{\contraction}{\rho}
\newcommand{\sparsity}{\Delta}
\newcommand{\sparsityr}{\Delta_r}
\newcommand{\lipschitz}{L}
\newcommand{\lyapunov}{\mathcal{L}}

\newcommand{\E}{\mathbb{E}}
\newcommand{\Econd}{\mathbf{E}}
\newcommand{\ind}{\mathbbm{1}}

\newcommand{\ASAGA}{\textsc{Asaga}}
\newcommand{\SAGA}{\textsc{Saga}}
\newcommand{\SAG}{\textsc{Sag}}

\newcommand{\SVRG}{\textsc{Svrg}}
\newcommand{\Hogwild}{\textsc{Hogwild}}
\newcommand{\SDCA}{\textsc{Sdca}}
\newcommand{\SGD}{\textsc{Sgd}}
\newcommand{\HSAG}{\textsc{Hsag}}
\newcommand{\KROMAGNON}{\textsc{Kromagnon}}
\newcommand{\AHSVRG}{\textsc{Ahsvrg}}

\newtheorem{asm}[theorem]{Assumption}
\newtheorem{prop}[theorem]{Property}

\renewcommand{\algorithmicloop}{\textbf{keep doing in parallel}} 
\renewcommand{\algorithmicendloop}{\algorithmicend\ \textbf{parallel loop}}

\usepackage{lastpage}
\jmlrheading{19}{2018}{1-\pageref{LastPage}}{11/17; Revised 7/18}{12/18}{17-650}{R\'emi Leblond, Fabian Pedregosa and Simon Lacoste-Julien}

\ShortHeadings{Improved Parallel Optimization Analysis for Stochastic Incremental Methods}{Leblond, Pedregosa and Lacoste-Julien}
\firstpageno{1}

\begin{document}
\title{Improved Asynchronous Parallel Optimization Analysis for Stochastic Incremental Methods}

\author{\name R\'emi Leblond \email remi.leblond@inria.fr \\
       \addr INRIA - Sierra Project-Team\\
       \'Ecole Normale Sup\'erieure, Paris\\       
       \AND
       \name Fabian Pedregosa \email f@bianp.net \\
       \addr INRIA - Sierra Project-Team\\
       \'Ecole Normale Sup\'erieure, Paris\\       
       \AND
       \name Simon Lacoste-Julien \email slacoste@iro.umontreal.ca \\
       \addr Department of CS \& OR (DIRO) \\
       Universit\'e de Montr\'eal, Montr\'eal\\ }

\editor{Tong Zhang}

\maketitle

\begin{abstract}%
As data sets continue to increase in size and multi-core computer architectures are developed, asynchronous parallel optimization algorithms become more and more essential to the field of Machine Learning.
Unfortunately, conducting the theoretical analysis asynchronous methods is difficult, notably due to the introduction of delay and inconsistency in inherently sequential algorithms. 
Handling these issues often requires resorting to simplifying but unrealistic assumptions. 
Through a novel perspective, we revisit and clarify a subtle but important technical issue present in a large fraction of the recent convergence rate proofs for asynchronous parallel optimization algorithms, and propose a simplification of the recently introduced ``perturbed iterate'' framework that resolves it.
We demonstrate the usefulness of our new framework by analyzing three distinct asynchronous parallel incremental optimization algorithms: \Hogwild\ (asynchronous \SGD), \KROMAGNON\ (asynchronous \SVRG) and \ASAGA, a novel asynchronous parallel version of the incremental gradient algorithm \SAGA\ that enjoys fast linear convergence rates.
We are able to both remove problematic assumptions and obtain better theoretical results.
Notably, we prove that \ASAGA\ and \KROMAGNON\ can obtain a theoretical linear speedup on multi-core systems even without sparsity assumptions.
We present results of an implementation on a 40-core architecture illustrating the practical speedups as well as the hardware overhead.
Finally, we investigate the overlap constant, an ill-understood but central quantity for the theoretical analysis of asynchronous parallel algorithms.
We find that it encompasses much more complexity than suggested in previous work, and often is order-of-magnitude bigger than traditionally thought.
\end{abstract}


\section{Introduction}
We consider the unconstrained optimization problem of minimizing a \emph{finite sum} of smooth convex functions:
\begin{equation} \label{eq:finiteSum}
\min_{x \in \mathbb{R}^d} f(x), \quad f(x) := \frac{1}{n} \sum_{i=1}^{n} f_i(x),
\end{equation}
where each $f_i$ is assumed to be convex with $\lipschitz$-Lipschitz continuous gradient, $f$ is $\mu$-strongly convex and $n$ is large (for example, the number of data points in a regularized empirical risk minimization setting).
We define a condition number for this problem as $\kappa := \nicefrac{\lipschitz}{\mu}$, as is standard in the finite sum literature.\footnote{Since we have assumed that each individual $f_i$ is $L$-smooth, $f$ itself is $L$-smooth – but its smoothness constant $L_f$ could be much smaller. While the more classical condition number is $\kappa_b := \nicefrac{L_f}{\strongconvex}$, our rates are in terms of this bigger $\nicefrac{L}{\strongconvex}$ in this paper.}
A flurry of randomized incremental algorithms (which at each iteration select $i$ at random and process only one gradient $f'_i$) have recently been proposed to solve~\eqref{eq:finiteSum} with a fast\footnote{Their complexity in terms of gradient evaluations to reach an accuracy of $\epsilon$ is $O((n+\kappa)\log(\nicefrac{1}{\epsilon}))$, in contrast to $O(n\kappa_b\log(\nicefrac{1}{\epsilon}))$ for batch gradient descent in the worst case.} linear convergence rate, such as \SAG~\citep{SAG}, \SDCA~\citep{SDCA}, \SVRG~\citep{svrg} and \SAGA~\citep{SAGA}.
These algorithms can be interpreted as variance reduced versions of the stochastic gradient descent (\SGD) algorithm, and they have demonstrated both theoretical and practical improvements over \SGD\ (for the \emph{finite sum} optimization problem~\ref{eq:finiteSum}).

In order to take advantage of the multi-core architecture of modern computers, the aforementioned optimization algorithms need to be adapted to the asynchronous parallel setting, where multiple threads work concurrently.
Much work has been devoted recently in proposing and analyzing asynchronous parallel variants of algorithms such as \SGD~\citep{hogwild}, \SDCA~\citep{asyncSDCA2015} and \SVRG~\citep{smola,mania,asySVRG}.
Among the incremental gradient algorithms with fast linear convergence rates that can optimize~\eqref{eq:finiteSum} in its general form, only \SVRG\ had had an asynchronous parallel version proposed.\footnote{We note that \SDCA\ requires the knowledge of an explicit $\mu$-strongly convex regularizer in~\eqref{eq:finiteSum}, whereas \SAG~/ \SAGA\ are adaptive to any local strong convexity of $f$~\citep{laggedsaga,SAGA}. The variant of \SVRG\ from~\citet{qsaga} is also adaptive (we review this variant in Section~\ref{ssec:svrgalgos}).}
No such adaptation had been attempted for \SAGA\ until~\cite{leblond2016Asaga}, even though one could argue that it is a more natural candidate as, contrarily to \SVRG, it is not epoch-based and thus has no synchronization barriers at all.
The present paper is an extended journal version of the conference paper from~\cite{leblond2016Asaga}.

The usual frameworks for asynchronous analysis are quite intricate (see Section~\ref{ssec:labelingIssue}) and thus require strong simplifying assumptions.
They are not well suited to the study of complex algorithms such as \SAGA.
We therefore base our investigation on the newly proposed ``perturbed iterate'' framework introduced in~\citet{mania}, which we also improve upon in order to properly analyze \SAGA.
Deriving a framework in which the analysis of \SAGA\ is possible enables us to highlight the deficiencies of previous frameworks and to define a better alternative.
Our new approach is not limited to \SAGA\ but can be used to investigate other algorithms and improve their existing bounds.

\paragraph{\textit{Contributions.}}
In Section~\ref{pif}, we propose a simplification of the ``perturbed iterate'' framework from~\citet{mania} as a basis for our asynchronous convergence analysis.
At the same time, through a novel perspective, we revisit and clarify a technical problem present in a large fraction of the literature on randomized asynchronous parallel algorithms (with the exception of~\citealt{mania}, which also highlights this issue): namely, they all assume unbiased gradient estimates, an assumption that is \textit{inconsistent} with their proof technique without further unpractical synchronization assumptions.

In Section~\ref{scs:sparse_saga}, we present a novel sparse variant of \SAGA\ that is more adapted to the parallel setting than the original \SAGA\ algorithm.
In Section~\ref{sec:ASAGA}, we present \ASAGA, a lock-free asynchronous parallel version of Sparse \SAGA\ that does not require consistent read or write operations.
We give a tailored convergence analysis for \ASAGA. Our main result states that \ASAGA\ obtains the same geometric convergence rate per update as \SAGA\ when the overlap bound $\overlap$ (which scales with the number of cores) satisfies
$\overlap \leq \mathcal{O}(n)$ and $\overlap \leq \mathcal{O}({\scriptstyle \frac{1}{\sqrt{\sparsity}}} \max\{1,\frac{n}{\kappa} \})$, where $\sparsity \leq 1$ is a measure of the sparsity of the problem.
This notably implies that a linear speedup is theoretically possible even without sparsity in the well-conditioned regime where $n \gg \kappa$.
This result is in contrast to previous analysis which always required some sparsity assumptions.

In Section~\ref{svrg}, we revisit the asynchronous variant of \SVRG\ from~\citet{mania}, \KROMAGNON, while removing their gradient bound assumption (which was inconsistent with the strongly convex setting).\footnote{Although the authors mention that this gradient bound assumption can be enforced through the use of a thresholding operator, they do not explain how to handle the interplay between this non-linear operator and the asynchrony of the algorithm. Their theoretical analysis relies on the linearity of the operations (e.g. to derive~\citep[Eq. (2.6)]{mania}), and thus this claim is not currently supported by theory (note that a strongly convex function over an unbounded domain always has unbounded gradients).}
We prove that the algorithm enjoys the same fast rates of convergence as \SVRG\ under similar conditions as \ASAGA\ -- whereas the original paper only provided analysis for slower rates (in both the sequential and the asynchronous case), and thus less meaningful speedup results.

In Section~\ref{sec:SGD}, in order to show that our improved ``after read'' perturbed iterate framework can be used to revisit the analysis of other optimization routines with correct proofs that do not assume homogeneous computation, we provide the analysis of the \Hogwild\ algorithm first introduced in~\citet{hogwild}.
Our framework allows us to remove the classic gradient bound assumption and to prove speedups in more realistic settings.

In Section~\ref{results}, we provide a practical implementation of \ASAGA\ and illustrate its performance on a 40-core architecture, showing improvements compared to asynchronous variants of \SVRG\ and \SGD.
We also present experiments on the overlap bound $\overlap$, showing that it encompasses much more complexity than suggested in previous work.

\paragraph{\textit{Related Work.}}
The seminal textbook of~\citet{bertsekasParalle1989} provides most of the foundational work for parallel and distributed optimization algorithms.
An asynchronous variant of \SGD\ with constant step size called \Hogwild\ was presented by~\citet{hogwild}; part of their framework of analysis was re-used and inspired most of the recent literature on asynchronous parallel optimization algorithms with convergence rates, including asynchronous variants of coordinate descent~\citep{asyncCD2015}, \SDCA~\citep{asyncSDCA2015}, \SGD\ for non-convex problems~\citep{taming,asyncSGDNonConvex2015}, \SGD\ for stochastic optimization~\citep{duchi} and \SVRG~\citep{smola,asySVRG}.
These papers make use of an unbiased gradient assumption that is not consistent with the proof technique, and thus suffers from technical problems\footnote{Except~\citep{duchi} that can be easily fixed by incrementing their global counter \emph{before} sampling.} that we highlight in Section~\ref{ssec:labelingIssue}.

The ``perturbed iterate'' framework presented in~\citet{mania} is to the best of our knowledge the only one that does not suffer from this problem, and our convergence analysis builds heavily from their approach, while simplifying it.
In particular, the authors assumed that $f$ was both strongly convex and had a bound on the gradient, two \emph{inconsistent} assumptions in the unconstrained setting that they analyzed.
We overcome these difficulties by using tighter inequalities that remove the requirement of a bound on the gradient. We also propose a more convenient way to label the iterates (see Section~\ref{ssec:labelingIssue}).
The sparse version of \SAGA\ that we propose is also inspired from the sparse version of \SVRG\ proposed by~\citet{mania}.

\citet{smola} presents a hybrid algorithm called \HSAG\ that includes \SAGA\ and \SVRG\ as special cases.
Their asynchronous analysis is epoch-based though, and thus does not handle a fully asynchronous version of \SAGA\ as we do.
Moreover, they require consistent reads and do not propose an efficient sparse implementation for \SAGA, in contrast to \ASAGA.

\citet{cyclades} proposes a black box mini-batch algorithm to parallelize \SGD-like methods while maintaining serial equivalence through smart update partitioning.
When the data set is sparse enough, they obtain speedups over ``\Hogwild'' implementations of \SVRG\ and \SAGA.\footnote{By ``\Hogwild'', the authors mean asynchronous parallel variants where cores independently run the sequential update rule.}
However, these ``\Hogwild'' implementations appear to be quite suboptimal, as they do not leverage data set sparsity efficiently: they try to adapt the ``lazy updates'' trick from~\citet{laggedsaga} to the asynchronous parallel setting -- which as discussed in Appendix~\ref{apx:DifficultyLagged} is extremely difficult -- and end up making several approximations which severely penalize the performance of the algorithms. 
In particular, they have to use much smaller step sizes than in the sequential version, which makes for worse results.

\citet{pedregosa2017maga} extend the \ASAGA\ algorithm presented in Section~\ref{sec:ASAGA} to the proximal setting.

\paragraph{\textit{Notation.}}
We denote by $\E$ a full expectation with respect to all the randomness in the system, and by $\Econd$ the \emph{conditional} expectation of a random~$i$ (the index of the factor $f_i$ chosen in \SGD\ and other algorithms), conditioned on all the past, where ``past'' will be clear from the context.
$[x]_v$ represents the coordinate~$v$ of the vector~$x \in \mathbb{R}^d$.
For \textit{sequential} algorithms, $x^+$ is the updated parameter vector after one algorithm iteration.

\section{Revisiting the Perturbed Iterate Framework for Asynchronous Analysis}\label{pif}
As most recent parallel optimization contributions, we use a similar hardware model to~\citet{hogwild}.
We consider multiple cores which all have read and write access to a shared memory.
The cores update a central parameter vector in an asynchronous and lock-free fashion.
Unlike~\citet{hogwild}, we \emph{do not} assume that the vector reads are consistent: multiple cores can read and write different coordinates of the shared vector at the same time. This also implies that a full vector read for a core might not correspond to any consistent state in the shared memory at any specific point in time.

\subsection{Perturbed Iterate Framework}
We first review the ``perturbed iterate'' framework recently introduced by~\citet{mania} which will form the basis of our analysis.
In the sequential setting, stochastic gradient descent and its variants can be characterized by the following update rule:
\begin{equation}
x_{t+1} = x_t -\stepsize g(x_t, i_t)\,,
\end{equation}
where $i_t$ is a random variable independent from $x_t$ and we have the unbiasedness condition $\Econd g(x_t, i_t) = f'(x_t)$ (recall that $\Econd$ is the relevant-past conditional expectation with respect to $i_t$).

Unfortunately, in the parallel setting, we manipulate stale, inconsistent reads of shared parameters and thus we do not have such a straightforward relationship.
Instead, \citet{mania} proposed to distinguish $\hat x_t$, the actual value read by a core to compute an update, from $x_t$, a ``virtual iterate'' that we can analyze and is \emph{defined} by the update equation:
\begin{equation}\label{eq:sgdupdate}
x_{t+1} := x_t -\stepsize g(\hat x_t, i_t)\,.
\end{equation}
We can thus interpret $\hat x_t$ as a noisy (perturbed) version of $x_t$ due to the effect of asynchrony.

We formalize the precise meaning of $x_t$ and $\hat x_t$ in the next section.
We first note that all references mentioned in the related work section that analyzed asynchronous parallel randomized algorithms assumed that the following unbiasedness condition holds:
\begin{equation} \label{eq:unbiasedness}
\left[{\text{unbiasedness condition}}\right] \quad \Econd [g(\hat x_{t}, i_t) | \hat x_t] = f'(\hat x_t). \,\, \footnote{We note that to be completely formal and define this conditional expectation more precisely, one would need to define another random vector that describes the entire system randomness, including all the reads, writes, delays, etc. Conditioning on $\hat{x}_t$ in~\eqref{eq:unbiasedness} is actually a shorthand to indicate that we are conditioning on all the relevant ``past'' that defines both the value of $\hat{x}_t$ as well as the fact that it was the $t^{\text{th}}$ labeled element. For clarity of exposition, we will not go into this level of technical detail, but one could define the appropriate sigma fields to condition on in order to make this equation fully rigorous.}
\end{equation}
This condition is at the heart of most convergence proofs for randomized optimization methods.\footnote{A notable exception is \SAG~\citep{SAG} which has biased updates and thus requires a significantly more complex convergence proof. Making \SAG\ unbiased leads to \SAGA~\citep{SAGA}, with a much simpler convergence proof.}
\citet{mania} correctly pointed out that most of the literature thus made the often implicit assumption that $i_t$ is independent of $\hat{x}_t$.
But as we explain below, this assumption is incompatible with a non-uniform asynchronous model in the analysis approach used in most of the recent literature.

\subsection{On the Difficulty of Labeling the Iterates} \label{ssec:labelingIssue}
Formalizing the meaning of $x_t$ and $\hat x_t$ highlights a subtle but important difficulty arising when analyzing \emph{randomized} parallel algorithms: what is the meaning of $t$?
This is the problem of \emph{labeling} the iterates for the purpose of the analysis, and this labeling can have randomness itself that needs to be taken in consideration when interpreting the meaning of an expression like $\E[x_t]$.
In this section, we contrast three different approaches in a unified framework.
We notably clarify the dependency issues that the labeling from~\citet{mania} resolves and propose a new, simpler labeling which allows for much simpler proof techniques.

We consider algorithms that execute in parallel the following four steps, where $t$ is a global labeling that needs to be defined:\footnote{Observe that contrary to most asynchronous algorithms, we choose to read the shared parameter vector \emph{before} sampling the next data point. We made this design choice to emphasize that in order for $\hat x_t$ and $i_t$ to be independent -- which will prove crucial for the analysis -- the reading of the shared parameter has to be independent of the sampled data point. Although in practice one would prefer to only read the necessary parameters \emph{after} sampling the relevant data point, for the sake of the analysis we cannot allow this source of dependence. We note that our analysis could also handle reading the parameter first and then sampling as long as independence is ensured, but for clarity of presentation, we decided to make this independence explicit.
	
\citet{mania} make the opposite presentation choice. In their main analysis, they explicitly assume that $\hat x_t$ and $i_t$ are independent, although they explain that it is not the case in practical implementations. The authors then propose a scheme to handle the dependency directly in their appendix. However, this ``fix'' can only be applied in a restricted setup: only for the \Hogwild\ algorithm, with the assumption that the norm of the gradient is uniformly bounded. Furthermore, even in this restricted setup, the scheme leads to worsened theoretical results (the bound on $\overlap$ is $\kappa^2$ worse). Applying it to a more complex algorithm such as \KROMAGNON\ or \ASAGA\ would mean overcoming several significant hurdles and is thus still an open problem.

In the absence of a better option, we choose to enforce the independence of $\hat x_t$ and $i_t$ with our modified steps ordering.\label{footnote:ordering}} 
\begin{equation} \label{eq:updates}
\begin{minipage}{0.9\textwidth}
\begin{enumerate}[leftmargin=3em, topsep=0mm, itemsep=0mm]
\item Read the information in shared memory ($\hat{x}_t$).
\item Sample $i_t$.
\item Perform some computations using ($\hat{x}_t, i_t$).
\item Write an update to shared memory.
\end{enumerate}
\end{minipage}
\end{equation}

\paragraph{\textit{The ``After Write'' Approach.}} We call the ``after write'' approach the standard global labeling scheme used in~\citet{hogwild} and re-used in all the later papers that we mentioned in the related work section, with the notable exceptions of~\citet{mania} and~\citet{duchi}. In this approach, $t$ is a (virtual) global counter recording the number of \emph{successful writes} to the shared memory $x$ (incremented after step~4 in~\ref{eq:updates}); $x_t$ thus represents the (true) content of the shared memory after $t$ updates.
The interpretation of the crucial equation~\eqref{eq:sgdupdate} then means that $\hat{x}_t$ represents the (delayed) local copy value of the core that made the $(t+1)^{\mathrm{th}}$ successful update; $i_t$ represents the factor sampled by this core for this update.
Notice that in this framework, the value of $\hat x_t$ and $i_t$ is unknown at ``time~$t$''; we have to wait to the later time when the next core writes to memory to finally determine that its local variables are the ones labeled by $t$.
We thus see that here~$\hat x_t$ and~$i_t$ are not necessarily independent -- they share dependence through the assignment of the $t$ label.
In particular, if some values of~$i_t$ yield faster updates than others, it will influence the label assignment defining $\hat x_t$.
We provide a concrete example of this possible dependency in Figure~\ref{fig:afterwritebias}.

\begin{figure}
	\centering
	\begin{tabular}{llll}
		& \multicolumn{1}{c}{$f_1$} & $f_2$   \\
		\midrule
		core 1 & \multicolumn{1}{c}{$\times$} &          \\
		core 2 & \multicolumn{1}{c}{$\times$} &          \\
		\midrule[0.9pt]
		\vspace{4mm}
		{\color{red}$f'_{i_0}(\hat{x}_0)$} & $f'_1(\hat{x}_0)$ &  \\
	\end{tabular}
	\begin{tabular}{llll}
		& \multicolumn{1}{c}{$f_1$} & $f_2$ \\
		\midrule
		core 1 & \multicolumn{1}{c}{$\times$} &          \\
		core 2 &  &  \multicolumn{1}{c}{$\times$}     \\
		\midrule[0.9pt]
		\vspace{4mm}
		 & $f'_1(\hat{x}_0)$ &  \\
	\end{tabular}
	\begin{tabular}{llll}
		& \multicolumn{1}{c}{$f_1$} & $f_2$ \\
		\midrule
		core 1 &  &    \multicolumn{1}{c}{$\times$}      \\
		core 2 & \multicolumn{1}{c}{$\times$} &         \\
		\midrule[0.9pt]
		\vspace{4mm}
		 & $f'_1(\hat{x}_0)$ &  \\
	\end{tabular}
	\begin{tabular}{llll}
		& \multicolumn{1}{c}{$f_1$} & $f_2$ \\
		\midrule
		core 1 & &    \multicolumn{1}{c}{$\times$}      \\
		core 2 & &     \multicolumn{1}{c}{$\times$}     \\
		\midrule[0.9pt]
		\vspace{4mm}
		 &  & $f'_2(\hat{x}_0)$ \\
	\end{tabular}
	\caption{\label{fig:afterwritebias}
		\footnotesize
		Suppose that we have two cores and that $f$ has two factors: $f_1$ which has support on only one variable, and $f_2$ which has support on $10^6$ variables and thus yields a gradient step that is significantly more expensive to compute.
		$x_0$ is the initial content of the memory, and we do not know yet whether $\hat{x}_0$ is the local copy read by the first core or the second core, but we are sure that $\hat{x}_0 = x_0$ as no update can occur in shared memory without incrementing the counter.
		There are four possibilities for the next step defining $x_1$ depending on which index $i$ was sampled on each core.
		If any core samples $i=1$, we know that $x_1 = x_0 - \stepsize f'_1(x_0)$ as it will be the first (much faster update) to complete.
		This happens in 3 out of 4 possibilities; we thus have that $\E x_1 = x_0 - \stepsize (\frac{3}{4} f'_1(x_0) + \frac{1}{4} f '_2(x_0))$. We see that this analysis scheme \emph{does not} satisfy the crucial unbiasedness condition~\eqref{eq:unbiasedness}.
		To understand this subtle point better, note that in this very simple example, $i_0$ and $i_1$ are not independent. We can show that $P(i_1=2 \mid i_0=2) =1$. They share dependency \emph{through the labeling assignment}.
	}
\end{figure}

The only way we can think to resolve this issue and ensure unbiasedness is to assume that the computation time for the algorithm running on a core is independent of the sample $i$ chosen.
This assumption seems overly strong in the context of potentially heterogeneous factors $f_i$'s, and is thus a fundamental flaw for analyzing non-uniform asynchronous computation that has mostly been ignored in the recent asynchronous optimization literature.\footnote{We note that \citet{bertsekasParalle1989} briefly discussed this issue (see Section~7.8.3), stressing that their analysis for SGD required that the scheduling of computation was independent from the randomness from SGD, but they did not offer any solution if this assumption was not satisfied. Both the ``before read'' labeling from~\citet{mania} and our proposed ``after read'' labeling resolve this issue.}

\paragraph{\textit{The ``Before Read'' Approach.}}
\citet{mania} address this issue by proposing instead to increment the global~$t$~counter just \emph{before} a new core starts to \emph{read} the shared memory (before step~1 in~\ref{eq:updates}).
In their framework, $\hat{x}_t$ represents the (inconsistent) read that was made by this core in this computational block, and $i_t$ represents the chosen sample.
The update rule~\eqref{eq:sgdupdate} represents a \emph{definition} of the meaning of $x_t$, which is now a ``virtual iterate'' as it does not necessarily correspond to the content of the shared memory at any point.
The real quantities manipulated by the algorithm in this approach are the $\hat{x}_t$'s, whereas $x_t$ is used only for the analysis -- consequently, the critical quantity we want to see vanish is $\E \|\hat x_t - x^*\|^2$.
The independence of~$i_t$ with~$\hat{x}_t$ can be simply enforced in this approach by making sure that the way the shared memory $x$ is read does not depend on $i_t$ (e.g. by reading all its coordinates in a fixed order). Note that this implies that we have to read all of $x$'s coordinates, regardless of the size of $f_{i_t}$'s support.
This is a much weaker condition than the assumption that all the computation in a block does not depend on $i_t$ as required by the ``after write'' approach, and is thus more reasonable.

\paragraph{\textit{A New Global Ordering: the ``After Read'' Approach.}}
The ``before read'' approach gives rise to the following complication in the analysis: $\hat{x}_t$ can depend on $i_r$ for $r > t$.
This is because $t$ is a global time ordering only on the assignment of computation to a core, not on when  $\hat{x}_t$ was finished being read.
This means that we need to consider both the ``future'' and the ``past'' when analyzing $x_t$.
To simplify the analysis, we thus propose a third way to label the iterates that we call ``after read'': $\hat{x}_t$ represents the $(t+1)^{\mathrm{th}}$ \emph{fully completed read} ($t$ incremented after step~1 in~\ref{eq:updates}).
As in the ``before read'' approach, we can ensure that $i_t$ is independent of $\hat{x}_t$ by ensuring that how we read does not depend on $i_t$.
But unlike in the ``before read'' approach, $t$ here now does represent a global ordering on the $\hat{x}_t$ iterates -- and thus we have that $i_r$ is independent of $\hat{x}_t$ for $r > t$.
Again using~\eqref{eq:sgdupdate} as the definition of the virtual iterate $x_t$ as in the perturbed iterate framework, we then have a very simple form for the value of $x_t$ and $\hat{x}_t$ (assuming atomic writes, see Property~\ref{eventconst} below):
\begin{equation} \label{eq:xhatUpdates}
\begin{aligned}
x_t &= x_0 - \stepsize \sum_{u = 0}^{t-1} g(\hat x_u, \hat \alpha^u, i_u) \, ;
\\
[\hat{x}_t]_v &= [x_0]_v - \stepsize
\mkern-36mu \sum_{\substack{u = 0 \\
		\text{u s.t. coordinate $v$ was written}\\
		\text{for $u$ before $t$} }
}^{t-1} \mkern-36mu  [g(\hat x_u, \hat \alpha^u, i_u)]_v \, .
\end{aligned}
\end{equation}
This proved crucial for our \ASAGA\ proof and allowed us to obtain better bounds for \Hogwild\ and the \KROMAGNON\ algorithm presented in~\citet{mania}.

The main idea of the perturbed iterate framework is to use this handle on $\hat x_t - x_t$ to analyze the convergence for $x_t$.
As $x_t$ is a virtual quantity, \citet{mania} supposed that there exists an index $T$ such that $x_T$ lives in shared memory ($T$ is a pre-set final iteration number after which all computation is completed, which means $x_T = \hat x_T$) and gave their convergence result for this $x_T$.

In this paper, we instead state the convergence results directly in terms of $\hat x_t$, thus avoiding the need for an unwieldy pre-set final iteration counter, and also enabling guarantees during the entire course of the algorithm.

\begin{remark}
	As mentioned in footnote~\ref{footnote:ordering}, \citet{mania} choose to sample a data point first and only then read the shared parameter vector in~\eqref{eq:updates}.
	One advantage of this option is that it allows for reading only the relevant dimensions of the parameter vector, although it means losing the crucial independence property between $\hat x_t$ and $i_t$.

	We can thus consider that their labeling approach is ``after sampling'' rather than ``before read'' (both are equivalent given their ordering).
	If we take this view, then by switching the order of the sampling and the reading steps in their setup, the ``after sampling'' approach becomes equivalent to our proposed ``after read'' labeling.
	
	However, the framework in which they place their analysis is the ``before read'' approach as described above, which results in having to take into account troublesome ``future'' terms in~\eqref{eq:xhatUpdates}.
	These additional terms make the analysis considerably harder and ultimately lead to worse theoretical results.
\end{remark}

\section{Asynchronous Parallel Sparse \SAGA}
We start by presenting Sparse \SAGA, a sparse variant of the \SAGA\ algorithm that is more adapted to the asynchronous parallel setting.
We then introduce \ASAGA, the asynchronous parallel version of Sparse \SAGA.
Finally, we state both convergence and speedup results for \ASAGA\ and give an outline of their proofs.
\subsection{Sparse \SAGA}\label{scs:sparse_saga}
Borrowing our notation from~\citet{qsaga}, we first present the original \SAGA\ algorithm and then describe our novel sparse variant.

\paragraph{\textit{Original \SAGA\ Algorithm.}} The standard \SAGA\ algorithm~\citep{SAGA} maintains two moving quantities to optimize~\eqref{eq:finiteSum}: the current iterate~$x$ and a table (memory) of historical gradients $(\alpha_i)_{i=1}^n$.\footnote{For linear predictor models, the memory $\alpha_i^0$ can be stored as a scalar. Following~\cite{qsaga}, $\alpha_i^0$ can be initialized to any convenient value (typically $0$), unlike the prescribed $f'_i(x_0)$ analyzed in~\citep{SAGA}.}
At every iteration, the \SAGA\ algorithm samples uniformly at random an index $i \in \{1,\ldots, n\}$, and then executes the following update on~$x$ and~$\alpha$ (for the unconstrained optimization version):
\begin{equation}\label{eq:SAGAupdate}
x^{+} = x - \stepsize \big(f'_i(x) - \alpha_i + \bar \alpha\big); \qquad  \alpha_i^+ = f'_i(x),
\end{equation}
where $\stepsize$ is the step size and $\bar \alpha := \nicefrac{1}{n} \sum_{i=1}^n \alpha_i$ can be updated efficiently in an online fashion. Crucially, $\Econd \alpha_i = \bar \alpha$ and thus the update direction is unbiased ($\Econd x^{+} = x - \stepsize f'(x)$).
Furthermore, it can be proven~\citep[see][]{SAGA} that under a reasonable condition on $\stepsize$, the update has vanishing variance, which enables the algorithm to converge linearly with a constant step size.

\paragraph{\textit{Motivation for a Variant.}}
In its current form, every \SAGA\ update is dense even if the individual gradients are sparse due to the historical gradient ($\bar \alpha$) term.
\citet{laggedsaga} introduced an implementation technique denoted lagged updates in which each iteration has a cost proportional to the size of the support of $f_i'(x)$.
However, this technique involves keeping track of past updates and is not easily adaptable  to the parallel setting (see Appendix~\ref{apx:DifficultyLagged}).
We therefore introduce Sparse \SAGA, a novel variant which explicitly takes sparsity into account and is easily parallelizable.

\paragraph{\textit{Sparse \SAGA\ Algorithm.}}
As in the Sparse \SVRG\ algorithm proposed in~\citet{mania}, we obtain Sparse \SAGA\ by a simple modification of the parameter update rule in~\eqref{eq:SAGAupdate} where $\bar{\alpha}$ is replaced by a sparse version equivalent in expectation:
\begin{equation} \label{eq:SparseSAGA}
x^{+} = x - \stepsize (f'_i(x) - \alpha_i + D_i \bar \alpha),
\end{equation}
where $D_i$ is a diagonal matrix that makes a weighted projection on the support of $f'_i$.
More precisely, let $S_i$ be the support of the gradient $f_i'$ function (i.e., the set of coordinates where $f_i'$ can be nonzero). Let $D$ be a $d\times d$ diagonal reweighting matrix, with coefficients~$\nicefrac{1}{p_v}$ on the diagonal, where $p_v$ is the probability that dimension~$v$ belongs to $S_i$ when $i$ is sampled uniformly at random in $\{1,...,n\}$.
We then define $D_i := P_{S_i} D$, where $P_{S_i}$ is the projection onto $S_i$. 
The reweighting by $D$ ensures that $\Econd D_i \bar \alpha = \bar{\alpha}$, and thus that the update is still unbiased despite the sparsifying projection.

\paragraph{\textit{Convergence Result for (Serial) Sparse \SAGA.}}
For clarity of exposition, we model our convergence result after the simple form of~\citet[Corollary~3]{qsaga}. Note that the rate we obtain for Sparse \SAGA\ is the same as the one obtained in the aforementioned reference for \SAGA.

\begin{theorem}\label{th1}
Let $\stepsize = \frac{a}{5\lipschitz}$ for any $a\leq1$.
Then \textnormal{Sparse \SAGA} converges geometrically in expectation with a rate factor of at least $\rho(a) = \frac{1}{5} \min\big\{\frac{1}{n}, a\frac{1}{\kappa}\big\}$, i.e., for $x_t$ obtained after $t$ updates, we have $\E \|x_t - x^*\|^2 \leq {(1-\rho)}^t \,  C_0$, where $C_0 := \|x_0 - x^*\|^2 + \frac{1}{5L^2} \sum_{i=1}^n\|\alpha_i^0 - f'_i(x^*)\|^2$.
\end{theorem}

\paragraph{\textit{Proof outline.}}\label{sparseoutline}
We reuse the proof technique from~\citet{qsaga}, in which a combination of classical strong convexity and Lipschitz inequalities is used to derive the inequality~\citep[Lemma~1]{qsaga}:
\begin{align}\label{eq:sparse}
\Econd \|x^{+} \! - \!x^*\|^2 \leq
&(1 \! - \! \stepsize\strongconvex) \|x \! -\! x^*\|^2
+ 2\stepsize^2 \Econd \|\alpha_i - f'_i(x^*)\|^2
+ (4 \stepsize^2 \lipschitz-2\stepsize)\big(f(x) - f(x^*)\big).
\end{align}
This gives a contraction term.
A Lyapunov function is then defined to control the two other terms.
To ensure our variant converges at the same rate as regular \SAGA, we only need to prove that the above inequality \citep[Lemma~1]{qsaga} is still verified.
To prove this, we derive close variants of equations $(6)$ and $(9)$ in their paper.
The rest of the proof can be reused without modification.
The full details can be found in Appendix~\ref{apxA}.

\paragraph{\textit{Comparison with Lagged Updates.}}
The lagged updates technique in \SAGA\ is based on the observation that the updates for component $[x]_v$ need not be applied until this coefficient needs to be accessed, that is, until the next iteration $t$ such that $ v \in S_{i_t}$.
We refer the reader to~\citet{laggedsaga} for more details.

Interestingly, the expected number of iterations between two steps where a given dimension $v$ is in the support of the partial gradient is $p_v^{-1}$, where $p_v$ is the probability that $v$ is in the support of the partial gradient at a given step.
$p_v^{-1}$ is precisely the term which we use to multiply the update to $[x]_v$ in Sparse \SAGA.
Therefore one may see the updates in Sparse \SAGA\ as \textit{anticipated} updates, whereas those in the~\citet{laggedsaga} implementation are \textit{lagged}.

The two algorithms appear to be very close, even though Sparse \SAGA\ uses an expectation to multiply a given update whereas the lazy implementation uses a random variable (with the same expectation). Sparse \SAGA\ therefore uses a slightly more aggressive strategy, which may explain the result of our experiments (see Section~\ref{sec:ssagacomp}): both Sparse \SAGA\ and \SAGA\ with lagged updates had similar convergence in terms of number of iterations, with the Sparse \SAGA\ scheme being slightly faster in terms of runtime.

Although Sparse \SAGA\ requires the computation of the $p_v$ probabilities, this can be done during a first pass throughout the data (during which constant step size \SGD\ may be used) at a negligible cost.

\subsection{Asynchronous Parallel Sparse \SAGA} \label{sec:ASAGA} \label{ssec:convergence}

\begin{figure*}[ttt!]
 \begin{minipage}[t]{0.49\textwidth}
   \begin{algorithm}[H]
     \caption{\ASAGA\ (analyzed algorithm)}
     \label{alg:theoretical}
     \label{theoreticalgo}
     \begin{algorithmic}[1]
	   \STATE Initialize shared variables $x$ and $(\alpha_i)_{i=1}^n$
	   \LOOP
	      \STATE $\hat x = $ inconsistent read of $x$
		  \STATE $\forall j$, $\hat \alpha_j = $ inconsistent read of $\alpha_j$
	   	  \STATE \textcolor{blue}{Sample $i$}  uniformly in $\{1,...,n\}$
	      \STATE Let $S_i$ be $f_i$'s support
	      \STATE $[\bar \alpha]_{S_i} = \nicefrac{1}{n} \sum_{k=1}^n [\hat \alpha_k]_{S_i}$
	      \STATE
    	  \STATE $[\delta x]_{S_i} = -\stepsize (f'_i(\hat x) - \hat \alpha_i + D_{i} [\bar \alpha]_{S_i})$
        	  \FOR{$v$ {\bfseries in} $S_i$}
        		 \STATE $[x]_v \leftarrow [x]_v + [\delta x]_v$      \hfill // atomic
	         \STATE $[\alpha_i]_v \leftarrow [f'_i( \hat x )]_v$
	         \STATE  // {\small `$\gets$' denotes a shared memory update.}
    		  \ENDFOR
	   \ENDLOOP
	  \end{algorithmic}
    \end{algorithm}
 \end{minipage}
 \hfill
 \begin{minipage}[t]{0.5\textwidth}
    \begin{algorithm}[H]
      \caption{\ASAGA\ (implementation)}
      \label{alg:sagasync}
      \begin{algorithmic}[1]
	    \STATE Initialize shared $x$, $(\alpha_i)_{i=1}^n$ and $\bar \alpha$
	    \LOOP
	      \STATE \textcolor{blue}{Sample $i$} uniformly in $\{1,...,n\}$
  	      \STATE Let $S_i$ be $f_i$'s support
 	      \STATE $[\hat x]_{S_i} = $ inconsistent read of $x$ on $S_i$
	      \STATE $\hat \alpha_i = $ inconsistent read of $\alpha_i$
	      \STATE $[\bar \alpha]_{S_i} = $ inconsistent read of $\bar \alpha$ on $S_i$
	      \STATE $[\delta \alpha]_{S_i} = f'_i([\hat x]_{S_i}) - \hat \alpha_i$
	      \STATE $[\delta x]_{S_i} = - \gamma ([\delta\alpha]_{S_i} + D_i [\bar \alpha]_{S_i})$
	      \FOR{$v$ {\bfseries in} $S_i$}
		      \STATE $[x]_v \gets [x]_v + [\delta x]_v$  \hfill // atomic
		      \STATE $[\alpha_i]_v \gets [\alpha_i]_v + [\delta \alpha]_v$ \hfill // atomic
	   	      \STATE $[\bar \alpha]_v \gets [\bar \alpha]_v + \nicefrac{1}{n}[\delta \alpha]_v$ \hfill // atomic
   	      \ENDFOR
	   \ENDLOOP
      \end{algorithmic}
    \end{algorithm}
 \end{minipage}
\end{figure*}

We describe \ASAGA, a sparse asynchronous parallel implementation of Sparse \SAGA, in Algorithm~\ref{alg:theoretical} in the theoretical form that we analyze, and in Algorithm~\ref{alg:sagasync} as its practical implementation.
We state our convergence result and analyze our algorithm using the improved perturbed iterate framework.

In the specific case of (Sparse) \SAGA, we have to add the additional read memory argument $\hat{\alpha}^t$ to our perturbed update~\eqref{eq:sgdupdate}:
\begin{equation}  \label{eq:PIupdate}
\begin{aligned}
x_{t+1} &:= x_t -\stepsize g(\hat x_t, \hat \alpha^t, i_t);
\\
g(\hat x_t, \hat \alpha^t, i_t) &:= f'_{i_t}(\hat x_t) - \hat \alpha_{i_t}^t + D_{i_t} \left({\textstyle \nicefrac{1}{n} \sum_{i=1}^n \hat{\alpha}_i^t }\right).
\end{aligned}
\end{equation}
Before stating our convergence result, we highlight some properties of Algorithm~\ref{alg:theoretical} and make one central assumption.

\begin{prop} [independence]
Given the ``after read'' global ordering, $i_r$ is independent of $\hat{x}_t$ $\forall r \geq t$.
\label{independence}
\end{prop}
The independence property for $r = t$ is assumed in most of the parallel optimization literature, even though it is not verified in case the ``after write'' labeling is used.
We emulate~\citet{mania} and \emph{enforce} this independence in Algorithm~\ref{alg:theoretical} by having the core read all the shared data parameters and historical gradients before starting their iterations.
Although this is too expensive to be practical if the data is sparse, this is required by the theoretical Algorithm~\ref{theoreticalgo} that we can analyze.
The independence for $r > t$ is a consequence of using the ``after read'' global ordering instead of the ``before read'' one.

\begin{prop}[unbiased estimator] \label{prop:unbiased}
The update, $g_t := g(\hat x_t, \hat \alpha^t, i_t)$, is an unbiased estimator of the true gradient at $\hat x_t$, i.e.~\eqref{eq:PIupdate} yields~\eqref{eq:unbiasedness} in conditional expectation.
\end{prop}
This property is crucial for the analysis, as in most related literature.
It follows by the independence of $i_t$ with $\hat{x}_t$ and from the computation of $\bar \alpha$ on line~7 of Algorithm~\ref{theoreticalgo}, which ensures that $\E \hat \alpha_i = 1/n \sum_{k=1}^n [\hat \alpha_k]_{S_i} = [\bar \alpha]_{S_i}$, making the update unbiased. In practice, recomputing $\bar \alpha$ is not optimal, but storing it instead introduces potential bias issues in the proof (as detailed in Appendix~\ref{apx:Bias}).

\begin{prop} [atomicity]
The shared parameter coordinate update of $[x]_v$ on line~11 is atomic.
\label{eventconst}
\end{prop}
Since our updates are additions, there are no overwrites, even when several cores compete for the same resources.
In practice, this is enforced by using \textit{compare-and-swap} semantics, which are heavily optimized at the processor level and have minimal overhead.
Our experiments with non-thread safe algorithms (i.e. where this property is not verified, see Figure~\ref{fig:cas_comparison} of Appendix~\ref{apxE}) show that compare-and-swap is necessary to optimize to high accuracy.

Finally, as is standard in the literature, we make an assumption on the maximum delay that asynchrony can cause -- this is the \emph{partially asynchronous} setting as defined in~\citet{bertsekasParalle1989}:
\begin{asm}[bounded overlaps]  \label{boundedoverlap}
We assume that there exists a uniform bound, called~$\overlap$, on the maximum number of iterations that can overlap together. We say that iterations $r$ and $t$ overlap if at some point they are processed concurrently.
One iteration is being processed from the start of the reading of the shared parameters to the end of the writing of its update.
The bound~$\overlap$ means that iterations $r$ cannot overlap with iteration $t$ for $r \geq t + \tau+1$,  and thus that every coordinate update from iteration $t$ is successfully written to memory before the iteration $t + \overlap+1$ starts.
\end{asm}
Our result will give us conditions on $\overlap$ subject to which we have linear speedups.
$\overlap$ is usually seen as a proxy for $p$, the number of cores (which lowerbounds it).
However, though $\overlap$ appears to depend linearly on $p$, it actually depends on several other factors (notably the data sparsity distribution) and can be orders of magnitude bigger than $p$ in real-life experiments.
We can upper bound $\overlap$ by $(p-1)R$, where $R$ is the ratio of the maximum over the minimum iteration time (which encompasses theoretical aspects as well as hardware overhead).
More details can be found in Section~\ref{sec:overlap}.

\paragraph{\textit{Explicit effect of asynchrony.}}
By using the overlap Assumption~\ref{boundedoverlap} in the expression~\eqref{eq:xhatUpdates} for the iterates, we obtain the following explicit effect of asynchrony that is crucially used in our proof:
\begin{align}\label{eq:async}
\hat x_t - x_t = \stepsize \sum_{u=(t - \overlap)_+}^{t-1}G_{u}^t g(\hat x_{u}, \hat \alpha^u, i_{u}),
\end{align}
where $G_{u}^t$ are $d\times d$ diagonal matrices with terms in $\{0, +1\}$.
From our definition of~$t$ and~$x_t$, it is clear that every update in $\hat x_t$ is already in $x_t$ -- this is the $0$ case.
Conversely, some updates might be late: this is the $+1$ case.
$\hat x_t$ may be lacking some updates from the ``past" in some sense, whereas given our global ordering definition, it cannot contain updates from the ``future".

\subsection{Convergence and Speedup Results}
We now state our main theoretical results. We give a detailed outline of the proof in Section~\ref{proofoutline} and its full details in Appendix~\ref{apxB}.

We first define a notion of problem sparsity, as it will appear in our results.

\begin{definition}[Sparsity]
As in~\citet{hogwild}, we introduce $\sparsityr := \max_{v=1..d} |\{i : v \in S_i\}|$. $\sparsityr$ is
the maximum right-degree in the bipartite graph of the factors and the dimensions, i.e.,
the maximum number of data points with a specific feature. For succinctness, we also define $\sparsity := \sparsityr / n$. We have $1 \leq \sparsityr \leq n$, and hence $1/n \leq \sparsity \leq 1$.
\end{definition}

\subsubsection{Convergence and Speedup Statements}

\begin{theorem}[Convergence guarantee and rate of \ASAGA]\label{thm:convergence}
Suppose $\overlap < \nicefrac{n}{10}$.\footnote{\ASAGA\ can actually converge for any $\overlap$, but the maximum step size then has a term of $\exp(\overlap/n)$ in the denominator with much worse constants. 
See Appendix~\ref{apxB:lma3}.} Let
\begin{equation}\label{eq:condition}
a^*(\overlap) := \frac{1}{32 \left(1+ \overlap  \sqrt \sparsity \right) \xi(\kappa, \sparsity, \overlap)} \quad
\begin{aligned} &\text{where } \xi(\kappa, \sparsity, \overlap) := \sqrt{1 + \frac{1}{8 \kappa}  \min\{\frac{1}{\sqrt{\sparsity}}, \overlap\} } \\
& \text{\small{(note that $\xi(\kappa, \sparsity, \overlap) \approx 1$ unless $\kappa < \nicefrac{1}{\sqrt{\sparsity}}  \,\, (\leq \sqrt{n})$)}}.
\end{aligned}
\end{equation}
For any step size $\stepsize = \frac{a}{L}$ with $a \leq a^*(\overlap)$, the inconsistent read iterates of Algorithm~\ref{alg:theoretical} converge in expectation at a geometric rate of at least: $\contraction(a) = \frac{1}{5} \min \big\{\frac{1}{n},  a \frac{1}{\kappa}\big\},$
i.e., $\E f(\hat x_t)-f(x^*) \leq (1-\rho)^t \,  \tilde C_0$, where $\tilde C_0$ is a constant independent of $t$  ($\approx \frac{n}{\stepsize}C_0$ with $C_0$ as defined in Theorem~\ref{th1}).
\end{theorem}
This result is very close to \SAGA's original convergence theorem, but with the maximum step size divided by an extra $1+ \overlap \sqrt{\sparsity}$ factor. Referring to~\citet{qsaga} and our own Theorem~\ref{th1}, the rate factor for \SAGA\ is $\min\{1/n, a/\kappa\}$ up to a constant factor. Comparing this rate with Theorem~\ref{thm:convergence} and inferring the conditions on the maximum step size $a^*(\overlap)$, we get the following conditions on the overlap $\overlap$ for \ASAGA\ to have the same rate as \SAGA\ (comparing upper bounds).
\begin{corollary}[Speedup condition]\label{thm:bigdata}\label{thm:illcondition}
Suppose $\overlap \leq \mathcal{O}(n)$ and $\overlap \leq \mathcal{O}({\scriptstyle \frac{1}{\sqrt{\sparsity}}} \max\{1,\frac{n}{\kappa} \})$. Then using the step size $\stepsize = \nicefrac{a^*(\overlap)}{L}\,$ from~\eqref{eq:condition}, \ASAGA\ converges geometrically with rate factor $\Omega( \min\{\frac{1}{n}, \frac{1}{\kappa}\})$ (similar to \SAGA), and is thus linearly faster than its sequential counterpart up to a constant factor.
Moreover, if $\overlap \leq \mathcal{O}(\frac{1}{\sqrt{\sparsity}})$, then a universal step size of $\Theta(\frac{1}{L})$ can be used for \ASAGA\ to be adaptive to local strong convexity with a similar rate to \SAGA\ (i.e., knowledge of $\kappa$ is not required).
\end{corollary}
Interestingly, in the well-conditioned regime ($n > \kappa$, where \SAGA\ enjoys a range of step sizes which all give the same contraction ratio), \ASAGA\ enjoys the same rate as \SAGA\ even in the non-sparse regime ($\sparsity = 1$) for $\overlap < \mathcal{O}(n/\kappa)$. This is in contrast to the previous work on asynchronous incremental gradient methods which required some kind of sparsity to get a theoretical linear speedup over their sequential counterpart~\citep{hogwild,mania}.
In the ill-conditioned regime ($\kappa > n$), sparsity is required for a linear speedup, with a bound on $\overlap$ of $\mathcal{O}(\sqrt{n})$ in the best-case (though degenerate) scenario where $\sparsity = 1/n$.

The proof for Corollary~\ref{thm:bigdata} can be found in Appendix~\ref{apx:proofASAGAspeedup}.

\paragraph{\textit{Comparison to related work.}}
\begin{itemize}[topsep=1mm, itemsep=-1mm]
\item We give the first convergence analysis for an asynchronous parallel version of \SAGA\ (note that \citet{smola} only covers an epoch based version of \SAGA\ with random stopping times, a fairly different algorithm).
\item Theorem~\ref{thm:convergence} can be directly extended to the a parallel extension of the \SVRG\ version from~\citet{qsaga} which is adaptive to the local strong convexity with similar rates (see Section~\ref{apx:SVRGext}).
\item In contrast to the parallel \SVRG\ analysis from~\citet[Thm. 2]{smola}, our proof technique handles inconsistent reads and a non-uniform processing speed across $f_i$'s. 
Our bounds are similar (noting that $\sparsity$ is equivalent to theirs), except for the adaptivity to local strong convexity: \ASAGA\ does not need to know $\kappa$ for optimal performance, contrary to parallel \SVRG\ (see Section~\ref{svrg} for more details).
\item In contrast to the \SVRG\ analysis from~\citet[Thm. 14]{mania}, we obtain a better dependence on the condition number in our rate ($1/\kappa$ vs. $1/\kappa^2$ on their work) and on the sparsity (they obtain $\overlap \leq \mathcal{O}(\sparsity^{\nicefrac{-1}{3}})$), while we furthermore remove their gradient bound assumption.
We also give our convergence guarantee on $\hat{x}_t$ \emph{during} the algorithm, whereas they only bound the error for the ``last'' iterate $x_T$.
\end{itemize}

\subsubsection{Proof Outline of Theorem~\ref{thm:convergence}}\label{proofoutline}
We give here an extended outline of the proof.
We detail key lemmas in Section~\ref{sec:keylemmas}.

\paragraph{\textit{Initial recursive inequality.}}
Let $g_t := g(\hat x_t, \hat \alpha^t, i_t)$. By expanding the update equation~\eqref{eq:PIupdate} defining the virtual iterate $x_{t+1}$ and introducing $\hat x_t$ in the inner product term, we obtain:
\begin{align}\label{eq:initrec}
\|x_{t+1} - x^*\|^2
&= \|x_t -\stepsize g_t -x^*\|^2
\nonumber\\
&= \|x_t -x^*\|^2 + \stepsize^2 \|g_t\|^2  -2\stepsize\langle x_t -x^*,  g_t\rangle
\nonumber\\
&= \|x_t -x^*\|^2 + \stepsize^2 \|g_t\|^2
-2\stepsize\langle \hat x_t -x^*,  g_t\rangle +2\stepsize\langle \hat x_t -x_t,  g_t\rangle  \, .
\end{align}
Note that we introduce $\hat x_t$ in the inner product because $g_t$ is a function of $\hat x_t$, not $x_t$.

In the sequential setting, we require $i_t$ to be independent of $x_t$ to obtain unbiasedness.
In the perturbed iterate framework, we instead require that $i_t$ is independent of $\hat x_t$ (see Property~\ref{independence}).
This crucial property enables us to use the unbiasedness condition~\eqref{eq:unbiasedness} to write:
$\E \langle \hat x_t -x^*,  g_t\rangle
= \E \langle \hat x_t -x^*,  f'(\hat x_t)\rangle$. Taking the expectation of~\eqref{eq:initrec} and using this unbiasedness condition we obtain an expression that allows us to use the $\mu$-strong convexity of $f$:\footnote{Note that here is our departure point with \citet{mania} who replaced the $f(\hat{x}_t)-f(x^*)$ term with the lower bound $\frac{\strongconvex}{2}\|\hat x_t - x^*\|^2$ in this relationship (see their Equation (2.4)), thus yielding an inequality too loose afterwards to get the fast rates for \SVRG.}
\begin{align} \label{eq:strongconvexity}
\langle \hat x_t -x^*,  f'(\hat x_t)\rangle &\geq f(\hat x_t) -f(x^*) +\frac{\strongconvex}{2}\|\hat x_t - x^*\|^2 .
\end{align}
With further manipulations on the expectation of~\eqref{eq:initrec}, including the use of the standard inequality $\|a + b\|^2 \leq 2\|a\|^2 + 2\|b\|^2$ (see Appendix~\ref{app:RecurDerivation}), we obtain our basic recursive contraction inequality:
\begin{align} \label{eq:RecursiveIneq1}
a_{t+1} &\leq
(1 -\frac{\stepsize \strongconvex}{2}) a_t
+ \stepsize^2 \E \|g_t\|^2
\underbrace{
	+ \stepsize\strongconvex \E\|\hat x_t - x_t\|^2
	+ 2\stepsize \E \langle \hat x_t -x_t,  g_t\rangle
}_{\text{\small additional asynchrony terms}}
-2\stepsize e_t  \, ,
\end{align}
where $a_t := \E \|x_t - x^*\|^2$ and $e_t := \E f(\hat x_t) - f(x^*)$.

Inequality~\eqref{eq:RecursiveIneq1} is a midway point between the one derived in the proof of Lemma~1 in~\citet{qsaga} and Equation~(2.5) in~\citet{mania}, because we use the tighter strong convexity bound~\eqref{eq:strongconvexity} than in the latter (giving us the important extra term $-2\stepsize e_t$).

In the sequential setting, one crucially uses the negative suboptimality term $-2\stepsize e_t$ to cancel the variance term $\stepsize^2 \E \|g_t\|^2$ (thus deriving a condition on $\stepsize$).
In our setting, we need to bound the additional asynchrony terms using the same negative suboptimality in order to prove convergence and speedup for our parallel algorithm -- this will give stronger constraints on the maximum step size.

The rest of the proof then proceeds as follows:
\begin{enumerate}
	\item By using the expansion~\eqref{eq:async} for $\hat{x}_t-x_t$, we can bound the additional asynchrony terms in~\eqref{eq:RecursiveIneq1} in terms of the past updates $(\E \|g_u\|^2)_{u\leq t}$. This gives Lemma~\ref{lma:1} below.
	\item We then bound the updates $\E \|g_t\|^2$ in terms of past suboptimalities $(e_v)_{v\leq u}$ by using standard \SAGA\ inequalities and carefully analyzing the update rule for $\alpha_i^+$~\eqref{eq:SAGAupdate} in expectation. This gives Lemma~\ref{lma:suboptgt} below.
	\item By applying Lemma~\ref{lma:suboptgt} to the result of Lemma~\ref{lma:1}, we obtain a master contraction inequality~\eqref{master} in terms of $a_{t+1}$, $a_t$ and $(e_u)_{u\leq t}$.
	\item We define a novel Lyapunov function $\lyapunov_t = \sum_{u=0}^t (1-\contraction)^{t-u}a_u$ and manipulate the master inequality to show that $\lyapunov_t$ is bounded by a contraction, subject to a maximum step size condition on~$\stepsize$ (given in Lemma~\ref{lma:3} below).
	\item Finally, we unroll the Lyapunov inequality to get the convergence Theorem~\ref{thm:convergence}.
\end{enumerate}

\subsubsection{Details}\label{sec:keylemmas}
We list the key lemmas below with their proof sketch, and pointers to the relevant parts of Appendix~\ref{apxB} for detailed proofs.

\begin{lemma}[Inequality in terms of $g_t := g(\hat x_{t}, \hat \alpha^t, i_{t})$]\label{lma:1}
	For all $t \geq 0$:
	\begin{equation} \label{eq:recursivegt}
	a_{t+1} \leq
	(1 - \frac{\stepsize\strongconvex}{2}) a_t + \stepsize^2 C_1\E\|g_t\|^2
	+ \stepsize^2 C_2\sum_{u=(t-\overlap)_+}^{t-1}\E\|g_{u}\|^2 - 2\stepsize e_t  \, ,
	\end{equation}
	\begin{equation}\label{eq:C12defs}
	\text{where} \quad
	C_1 := 1 + \sqrt{\sparsity}\overlap \quad
	\text{and} \quad
	C_2 :=  \sqrt{\sparsity} + \stepsize\strongconvex C_1 \,.\quad \footnote{Note that $C_2$ depends on $\stepsize$. In the rest of the paper, we write $C_2(\stepsize)$ instead of $C_2$ when we want to draw attention to that dependency.}
	\end{equation}
\end{lemma}
To prove this lemma we need to bound both $\E\|\hat x_t - x^*\|^2$ and $\E\langle \hat x_t -x_t,  g_t\rangle$ with respect to $(\E \|g_u\|^2)_{u\leq t}$.
We achieve this by crucially using Equation~\eqref{eq:async}, together with the following proposition, which we derive by a combination of Cauchy-Schwartz and our sparsity definition (see Section~\ref{apxB:lma1}).

\begin{proposition}\label{prop:1}
	For any $u \neq t$,
	\begin{align}\label{sparseproduct}
	\E |\langle g_{u}, g_t \rangle | &\leq \frac{\sqrt{\sparsity}}{2}(\E\|g_{u}\|^2 + \E\|g_{t}\|^2)  \, .
	\end{align}
\end{proposition}
To derive this essential inequality for both the right-hand-side terms of Eq.~\eqref{sparseproduct}, we start by proving a relevant property of $\sparsity$.
We reuse the sparsity constant introduced in~\citet{smola} and relate it to the one we have defined earlier, $\sparsityr$:
\begin{remark} \label{rmk:1}
	Let $D$ be the smallest constant such that:
	\begin{align} \label{sparsitycondition}
	\Econd \|x\|_i^2 = \frac{1}{n} \sum_{i=1}^n \|x\|_i^2 \leq D \|x\|^2 \quad  \forall x \in \mathbb{R}^d,
	\end{align}
	where $\|.\|_i$ is defined to be the $\ell_2$-norm restricted to the support $S_i$ of $f_i$.
	We have:
	\begin{equation}
	D = \frac{\sparsityr}{n} = \sparsity  \, .
	\end{equation}
\end{remark}

\begin{proof}
	We have:
	\begin{align}
	\Econd \|x\|_i^2 = \frac{1}{n} \sum_{i=1}^n \|x\|_i^2
	= \frac{1}{n} \sum_{i=1}^n \sum_{v \in S_i} [x]_v^2
	= \frac{1}{n} \sum_{v=1}^d \sum_{i \mid v \in S_i} [x]_v^2
	= \frac{1}{n} \sum_{v=1}^d \delta_v [x]_v^2 \, ,
	\end{align}
	where $\delta_v := |(i \mid v \in S_i)|$.
	This implies:
	\begin{align}
	D \geq \frac{1}{n} \sum_{v=1}^d \delta_v \frac{[x]_v^2}{\|x\|^2}  \, .
	\end{align}
	Since $D$ is the minimum constant satisfying this inequality, we have:
	\begin{align}
	D = \max_{x \in \mathbb{R}^d} \frac{1}{n} \sum_{v=1}^d \delta_v \frac{[x]_v^2}{\|x\|^2}  \, .
	\end{align}
	We need to find $x$ such that it maximizes the right-hand side term.
	Note that the vector $([x]_v^2 / \|x\|^2)_{v=1..d}$ is in the unit probability simplex, which means that an equivalent problem is the maximization over all convex combinations of $(\delta_v)_{v=1..d}$.
	This maximum is found by putting all the weight on the maximum $\delta_v$, which is $\sparsityr$ by definition.
	
	This implies that $\sparsity = \sparsityr / n$ is indeed the smallest constant satisfying~\eqref{sparsitycondition}.
\end{proof}

\paragraph{\textit{Proof of Proposition~\ref{prop:1}}}

Let $u \neq t$. Without loss of generality, $u < t$.\footnote{One only has to switch $u$ and $t$ if $u>t$.}
Then:

\begin{align}
\E |\langle g_{u}, g_t \rangle |
&\leq \E\|g_{u}\|_{i_t}\|g_t\|
\tag*{(Sparse inner product; support of $g_t$ is $S_{i_t}$)} \nonumber \\
&\leq \sqrt{\E\|g_{u}\|_{i_t}^2}\sqrt{\E\|g_{t}\|^2}
\tag*{(Cauchy-Schwarz for expectations)} \nonumber \\
&\leq \sqrt{\sparsity \E\|g_{u}\|^2}\sqrt{\E\|g_{t}\|^2}
\tag*{(Remark~\ref{rmk:1} and $i_t \perp\!\!\!\perp g_u, \forall u < t$)} \nonumber \\
&\leq \frac{\sqrt{\sparsity}}{2}(\E\|g_{u}\|^2 + \E\|g_{t}\|^2)  \, .
\tag*{(AM-GM inequality)}
\end{align}
All told, we have:
\begin{align}
\E |\langle g_{u}, g_t \rangle | &\leq \frac{\sqrt{\sparsity}}{2}(\E\|g_{u}\|^2 + \E\|g_{t}\|^2)  \, .
\end{align}

\begin{lemma} [Suboptimality bound on $\E \|g_t\|^2$]\label{lma:suboptgt}
	For all $t \geq 0$,
	\begin{equation}\label{gtbound}
	\E\|g_t\|^2
	\leq 4\lipschitz e_t
	+ \frac{4\lipschitz}{n} \sum_{u=1}^{t-1} (1 - \frac{1}{n})^{(t-2\overlap-u -1)_+} e_u
	+ 4\lipschitz (1 - \frac{1}{n})^{(t-\overlap)_+} \tilde e_0  \, ,
	\end{equation}
	where $\tilde e_0 := \frac{1}{2\lipschitz} \E\|\alpha_i^0 - f'_i(x^*)\|^2$.\footnote{We introduce this quantity instead of $e_0$ so as to be able to handle the arbitrary initialization of the $\alpha_i^0$.}
\end{lemma}
From our proof of convergence for Sparse \SAGA\ we know that (see Appendix~\ref{apxA}):
\begin{align}
\E\|g_t\|^2
&\leq 2 \E \|f'_{i_t}(\hat x_t)-f'_{i_t}(x^*)\|^2
+ 2 \E \|\hat \alpha_{i_t}^t - f'_{i_t}(x^*)\|^2 .
\end{align}
We can handle the first term by taking the expectation over a Lipschitz inequality (\citet[Equations 7 and 8]{qsaga}.
All that remains to prove the lemma is to express the $\E \|\hat \alpha_{i_t}^t - f'_{i_t}(x^*)\|^2$ term in terms of past suboptimalities.
We note that it can be seen as an expectation of past first terms with an adequate probability distribution which we derive and bound.

From our algorithm, we know that each dimension of the memory vector $[\hat \alpha_i]_v$ contains a partial gradient computed at some point in the past $[f'_i(\hat x_{u_{i, v}^{t}})]_v$\footnote{More precisely: $\forall t, i, v \hspace{0.5em}\exists u_{i, v}^t<t$ s.t. $[\hat \alpha_{i}^t]_v = [f'_{i}(\hat x_{u_{i, v}^{t}})]_v$.} (unless $u=0$, in which case we replace the partial gradient with $\alpha_i^0$).
We then derive bounds on $P(u_{i,v}^t = u)$ and sum on all possible $u$.
Together with clever conditioning, we obtain Lemma~\ref{lma:suboptgt} (see Section~\ref{apxB:lma2}).

\paragraph{\textit{Master inequality.}}
Let $H_t$ be defined as ${H_t: = \sum_{u=1}^{t-1} (1 - \frac{1}{n})^{(t-2\overlap-u -1)_+} e_u}$.
Then, by setting~\eqref{gtbound} into Lemma~\ref{lma:1}, we get (see Section~\ref{apxB:master}):
\begin{equation}\label{master}
\begin{aligned}
a_{t+1}
\leq &(1 - \frac{\stepsize\strongconvex}{2}) a_t
- 2\stepsize e_t
+ 4\lipschitz\stepsize^2 C_1 \big(e_t  + (1 - \frac{1}{n})^{(t-\overlap)_+} \tilde e_0 \big)
+ \frac{4\lipschitz\stepsize^2 C_1}{n} H_t
\\
&+4\lipschitz\stepsize^2 C_2\sum_{u=(t-\overlap)_+}^{t-1} (e_u +  (1 - \frac{1}{n})^{(u - \overlap)_+} \tilde e_0 \big)
+\frac{4\lipschitz\stepsize^2 C_2}{n} \sum_{u=(t-\overlap)_+}^{t-1} H_u  \, .
\end{aligned}
\end{equation}

\paragraph{\textit{Lyapunov function and associated recursive inequality.}}
We now have the beginning of a contraction with additional positive terms which all converge to $0$ as we near the optimum, as well as our classical negative suboptimality term.
This is not unusual in the variance reduction literature. One successful approach in the sequential case is then to define a Lyapunov function which encompasses all terms and is a true contraction~\citep[see][]{SAGA, qsaga}.
We emulate this solution here.
However, while all terms in the sequential case only depend on the current iterate, $t$, in the parallel case we have terms ``from the past'' in our inequality.
To resolve this issue, we define a more involved Lyapunov function which also encompasses past iterates:
\begin{align} \label{eq:LyapunovDefinition}
\lyapunov_t = \sum_{u=0}^t (1-\contraction)^{t-u}a_u, \quad 0<\contraction < 1,
\end{align}
where $\contraction$ is a target contraction rate that we define later.

Using the master inequality~\eqref{master}, we get (see Appendix~\ref{apxB:lyapunov}):
\begin{align}\label{Lyapunov}
\lyapunov_{t+1}
&= (1 - \contraction)^{t+1}a_0
+ \sum_{u=0}^t(1 - \contraction)^{t-u}a_{u+1}
\nonumber \\
&\leq (1 - \contraction)^{t+1}a_0 + (1-\frac{\stepsize\strongconvex}{2})\lyapunov_t + \sum_{u=1}^t r_u^t e_u + r_0^t \tilde e_0  \, .
\end{align}
The aim is to prove that $\lyapunov_t$ is bounded by a contraction.
We have two promising terms at the beginning of the inequality, and then we need to handle the last term.
Basically, we can rearrange the sums in~\eqref{master} to expose a simple sum of $e_u$ multiplied by factors $r_u^t$.

Under specific conditions on $\contraction$ and $\stepsize$, we can prove that $r_u^t$ is negative for all $u \geq 1$, which coupled with the fact that each $e_u$ is positive means that we can safely drop the sum term from the inequality.
The $r_0^t$ term is a bit trickier and is handled separately.

In order to obtain a bound on $e_t$ directly rather than on $\E \|\hat x_t - x^*\|^2$, we then introduce an additional $\stepsize e_t$ term on both sides of~\eqref{Lyapunov}.
The bound on $\stepsize$ under which the modified $r_t^t + \stepsize$ is negative is then twice as small (we could have used any multiplier between $0$ and $2\stepsize$, but chose $\stepsize$ for simplicity's sake).
This condition is given in the following Lemma.

\begin{lemma} [Sufficient condition for convergence]\label{lma:3}
	Suppose $\overlap < \nicefrac{n}{10}$ and $\contraction \leq \nicefrac{1}{4n}$. If
	\begin{align}\label{eq:lmacondition}
	\stepsize \leq \stepsize^* = \frac{1}{32\lipschitz (1 + \sqrt{\sparsity} \overlap) \sqrt{1 + \frac{1}{8\kappa} \min(\overlap, \frac{1}{\sqrt{\sparsity}})}}
	\end{align}
	then for all $u \geq 1$, the coefficients $r_u^t$ from~\eqref{Lyapunov} are negative.
	Furthermore, we have $r_t^t + \stepsize \leq 0$ and thus:
	\begin{align}
	\stepsize  e_t + \lyapunov_{t+1} \leq (1 - \contraction)^{t+1}a_0 + (1-\frac{\stepsize\strongconvex}{2})\lyapunov_t + r_0^t \tilde e_0 \,.
	\end{align}

\end{lemma}
We obtain this result after carefully deriving the $r_u^t$ terms.
We find a second-order polynomial inequality in $\stepsize$, which we simplify down to~\eqref{eq:lmacondition} (see Appendix~\ref{apxB:lma3}).

We can then finish the argument to bound the suboptimality error $e_t$. We have:
\begin{align}
\lyapunov_{t+1} \leq
\stepsize e_t + \lyapunov_{t+1}
&\leq (1-\frac{\stepsize\strongconvex}{2})\lyapunov_t  +(1 - \contraction)^{t+1} (a_0 + A \tilde e_0)  \, .
\end{align}

We have two linearly contracting terms.
The sum contracts linearly with the worst rate between the two (the smallest geometric rate factor).
If we define $\contraction^* := \nu \min(\contraction, \stepsize \strongconvex / 2)$, with $0 <\nu < 1$,\footnote{$\nu$ is introduced to circumvent the problematic case where $\contraction$ and  $\stepsize \strongconvex / 2$ are too close together.} then we get:
\begin{align}
\stepsize e_t + \lyapunov_{t+1}
&\leq (1-\frac{\stepsize\strongconvex}{2})^{t+1}\lyapunov_0 + (1 - \contraction^*)^{t+1} \frac{a_0 + A \tilde e_0}{1 -\eta}
\\
\stepsize e_t &\leq (1 - \contraction^*)^{t+1} \big(\lyapunov_0 + \frac{1}{1 -\eta} (a_0 + A  \tilde e_0) \big) \, ,
\end{align}
where $\eta := \frac{1-M}{1-\contraction^*}$ with $M :=\max(\contraction, \stepsize\strongconvex/2)$. Our geometric rate factor is thus $\contraction^*$ (see Appendix~\ref{apxB:th2}).

\section{Asynchronous Parallel \SVRG\ with the ``After Read'' Labeling}\label{svrg}
\paragraph{\textit{\ASAGA\ vs. asynchronous \SVRG.}}
There are several scenarios in which \ASAGA\ can be practically advantageous over its closely related cousin, asynchronous \SVRG\ (note though that ``asynchronous'' \SVRG\ still requires one synchronization step per epoch to compute a full gradient).

First, while \SAGA\ trades memory for less computation, in the case of generalized linear models the memory cost can be reduced to $\mathcal{O}(n)$, compared to $\mathcal{O}(d)$ for $\SVRG$~\citep{svrg}.
This is of course also true for their asynchronous counterparts.

Second, as \ASAGA\ does not require any synchronization steps, it is better suited to heterogeneous computing environments (where cores have different clock speeds or are shared with other applications).

Finally, \ASAGA\ does not require knowing the condition number $\kappa$ for optimal convergence in the sparse regime.
It is thus adaptive to local strong convexity, whereas \SVRG\ is not.
Indeed, \SVRG\ and its asynchronous variant require setting an additional hyper-parameter -- the epoch size $m$ -- which needs to be at least $\Omega(\kappa)$ for convergence but yields a slower effective convergence rate than \ASAGA\ if it is set much bigger than $\kappa$.
\SVRG\ thus requires tuning this additional hyper-parameter or running the risk of either slower convergence (if the epoch size chosen is much bigger than the condition number) or even not converging at all (if $m$ is chosen to be much smaller than $\kappa$).\footnote{Note that as \SAGA\ (and contrary to the original \SVRG) the \SVRG\ variant from~\citet{qsaga} does not require knowledge of $\kappa$ and is thus adaptive to local strong convexity, which carries over to its asynchronous adaptation that we analyze in Section~\ref{apx:SVRGext}.}

\paragraph{\textit{Motivation for analyzing asynchronous \SVRG.}}
Despite the advantages that we have just listed, in the case of complex models, the storage cost of \SAGA\ may become too expensive for practical use.
\SVRG~\citep{svrg} trades off more computation for less storage and does not suffer from this drawback. It can thus be applied to cases where \SAGA\ cannot~\citep[e.g. deep learning models, see][]{nonconvex}.

Another advantage of \KROMAGNON\ is that the historical gradient term $f'(\tilde x)$ is fixed during an epoch, while its \ASAGA\ equivalent, $\bar \alpha$, has to be updated at each iteration, either by recomputing if from the $\hat \alpha$ -- which is costly -- or by updating a maintained quantity -- which is cheaper but may ultimately result in introducing some bias in the update (see Appendix~\ref{apx:Bias} for more details on this subtle issue).

It is thus worthwhile to carry out the analysis of \KROMAGNON~\citep{mania}\footnote{The speedup analysis presented in~\citet{mania} is not fully satisfactory as it does not achieve state-of-the-art convergence results for either \SVRG\ or \KROMAGNON. Furthermore, we are able to remove their uniform gradient bound assumption, which is inconsistent with strong convexity.}, the asynchronous parallel version of \SVRG, although it has to be noted that since \SVRG\ requires regularly computing batch gradients, \KROMAGNON\ will present regular synchronization steps as well as coordinated computation -- making it less attractive for the asynchronous parallel setting.

We first extend our \ASAGA\ analysis to analyze the convergence of a variant of \SVRG\ presented in~\citet{qsaga}, obtaining exactly the same bounds. 
This variant improves upon the initial algorithm because it does not require tuning the epoch size hyperparameter and is thus adaptive to local strong convexity (see Section~\ref{ssec:svrgalgos}). 
Furthermore, it allows for a cleaner analysis where -- contrary to \SVRG\ -- we do not have to replace the final parameters of an epoch by one of its random iterates.

Then, using our ``after read'' labeling, we are also able to derive a convergence and speedup proof for \KROMAGNON, with comparable results to our \ASAGA\ analysis.
In particular, we prove that as for \ASAGA\, in the ``well-conditioned'' regime \KROMAGNON\ can achieve a linear speedup even without sparsity assumptions.

\subsection{\SVRG\ Algorithms}\label{ssec:svrgalgos}
We start by describing the original \SVRG\ algorithm, the variant given in~\citet{qsaga} and the sparse asynchronous parallel adaptation, \KROMAGNON.

\paragraph{\textit{Original \SVRG\ algorithm.}} 
The standard \SVRG\ algorithm~\citep{svrg} is very similar to \SAGA.
The main difference is that instead of maintaining a table of historical gradients, \SVRG\ uses a ``reference'' batch gradient $f'(\tilde{x})$, updated at regular intervals (typically every $m$ iterations, where $m$ is a hyper-parameter).
\SVRG\ is thus an epoch-based algorithm, where at the beginning of every epoch a reference iterate $\tilde{x}$ is chosen and its gradient is computed.
Then, at every iteration in the epoch, the algorithm samples uniformly at random an index $i \in \{1,\ldots, n\}$, and then executes the following update on~$x$:
\begin{equation}\label{eq:SVRGAupdate}
x^{+} = x - \stepsize \big(f'_i(x) - f'_i(\tilde x) + f'(\tilde x)\big) \,.
\end{equation}
As for \SAGA\ the update direction is unbiased ($\Econd x^{+} = x - \stepsize f'(x)$) and it can be proven~\citep[see][]{svrg} that under a reasonable condition on $\stepsize$ and $m$ (the epoch size), the update has vanishing variance, which enables the algorithm to converge linearly with a constant step size.

\paragraph{\textit{Hofmann's \SVRG\ variant.}}
\citet{qsaga} introduce a variant where the size of the epoch is a random variable.
At each iteration $t$, a first Bernoulli random variable $B_t$ with $p = \nicefrac{1}{n}$ is sampled.
If $B_t = 1$, then the algorithm updates the reference iterate, $\tilde{x} = x_t$ and computes its full gradient as its new ``reference gradient''.
If $B_t = 0$, the algorithm executes the normal \SVRG\ inner update.
Note that this variant is adaptive to local strong convexity, as it does not require the inner loop epoch size~$m = \Omega(\kappa)$ as a hyperparameter.
In that respect it is closer to \SAGA\ than the original \SVRG\ algorithm which is not adaptive.

\paragraph{\textit{\KROMAGNON.}}
\KROMAGNON, introduced in~\citet{mania} is obtained by using the same sparse update technique as Sparse \SAGA, and then running the resulting algorithm in parallel (see Algorithm~\ref{alg:kromagnon}).

\begin{figure*}[ttt!]
\centering
 \renewcommand{\algorithmicloop}{\textbf{for $i=1..m$ do in parallel (asynchronously)}} 
 \renewcommand{\algorithmicendloop}{\algorithmicend\ \textbf{parallel loop}}
 \begin{minipage}[t]{0.62\textwidth}
    \begin{algorithm}[H]
      \caption{\KROMAGNON~\citep{mania}}
      \label{alg:kromagnon}
      \begin{algorithmic}[1]
	    \STATE Initialize shared $x$ and $x_0$
   	    \WHILE{True}
     	    \STATE Compute in parallel $g = f'(x_0)$ (synchronously)

	    \LOOP
	      \STATE \textcolor{blue}{Sample $i$} uniformly in $\{1,...,n\}$
  	      \STATE Let $S_i$ be $f_i$'s support
 	      \STATE $[\hat x]_{S_i} = $ inconsistent read of $x$ on $S_i$
	      \STATE $[\delta x]_{S_i} = - \gamma ([f'_i(\hat x_t) - f'_i(x_0)]_{S_i} + D_i [g]_{S_i})$
	      \FOR{$v$ {\bfseries in} $S_i$}
	        \STATE $[x]_v = [x]_v + [\delta x]_v$  \hfill // atomic
	     \ENDFOR
	   \ENDLOOP
	   \STATE $x_0 = x$
   	    \ENDWHILE
      \end{algorithmic}
    \end{algorithm}
 \end{minipage}
 \hfill
  \begin{minipage}[t]{0.62\textwidth}
  	  \renewcommand{\algorithmicloop}{\textbf{while $s=0$ do in parallel (asynchronously)}} 
  	  \renewcommand{\algorithmicendloop}{\algorithmicend\ \textbf{parallel loop}}
  	\begin{algorithm}[H]
  		\caption{\AHSVRG}
  		\label{alg:ahsvrg}
  		\begin{algorithmic}[1]
  			\STATE Initialize shared $x$, $s$ and $x_0$
  			\WHILE{True}
	  			\STATE Compute in parallel $g = f'(x_0)$ (synchronously)
	  			\STATE $s = 0$
	  			\LOOP
		  			\STATE Sample $B$ with $p = \nicefrac{1}{n}$
		  			\IF{$B = 1$}
			  			\STATE $s = 1$
		  			\ELSE
			  			\STATE \textcolor{blue}{Sample $i$} uniformly in $\{1,...,n\}$
			  			\STATE Let $S_i$ be $f_i$'s support
			  			\STATE $[\hat x]_{S_i} = $ inconsistent read of $x$ on $S_i$
			  			\STATE $[\delta x]_{S_i} = - \gamma ([f'_i(\hat x_t) - f'_i(x_0)]_{S_i} + D_i [g]_{S_i})$
			  			\FOR{$v$ {\bfseries in} $S_i$}
				  			\STATE $[x]_v = [x]_v + [\delta x]_v$  \hfill // atomic
			  			\ENDFOR
		  			\ENDIF
	  			\ENDLOOP
	  			\STATE $x_0 = x$
  			\ENDWHILE
  		\end{algorithmic}
  	\end{algorithm}
  \end{minipage}
\end{figure*}

\subsection{Extension to the \SVRG\ Variant from~\citet{qsaga}}\label{apx:SVRGext}
We introduce \AHSVRG\ -- a sparse asynchronous parallel version for the \SVRG\ variant from~\citet{qsaga} -- in Algorithm~\ref{alg:ahsvrg}.
Every core runs stochastic updates independently as long as they are all sampling inner updates, and coordinate whenever one of them decides to do a batch gradient computation.
The one difficulty of this approach is that each core needs to be able to communicate to every other core that they should stop doing inner updates and start computing a synchronized batch gradient instead.

To this end, we introduce a new shared variable, $s$, which represents the ``state'' of the computation.
This variable is checked by each core $c$ before each update.
If $s = 1$, then another core has called for a batch gradient computation and core $c$ starts computing its allocated part of this computation.
If $s = 0$, core $c$ proceeds to sample a first random variable.
Then it either samples and performs an inner update and keeps going, or it samples a full gradient computation, in which case it updates $s$ to $1$ and starts computing its allocated part of the computation.
Once a full gradient is computed, $s$ is set to $0$ once again and every core resume their loop.

Our \ASAGA\ convergence and speedup proofs can easily be adapted to accommodate \AHSVRG\ since it is closer to \SAGA\ than the initial \SVRG\ algorithm.
To prove convergence, all one has to do is to modify Lemma~\ref{lma:suboptgt} very slightly (the only difference is that the $(t -2\overlap -u -1)_+$ exponent is replaced by $(t - \overlap - u - 1)_+$ and the rest of the proof can be used as is).
The justification for this small tweak is that the batch steps in \SVRG\ are fully synchronized. More details can be found in Appendix~\ref{apxB:ahsvrg}.

\subsection{Fast Convergence and Speedup Rates for \KROMAGNON}
We now state our main theoretical results. We give a detailed outline of the proof in Section~\ref{sec:proofSVRG} and its full details in Appendix~\ref{apx:SVRG}.
\begin{theorem}[Convergence guarantee and rate of \KROMAGNON]\label{thm:SVRG}
	Suppose the step size $\stepsize$ and epoch size $m$ are chosen such that the following condition holds:
	\begin{equation}
		0 < \theta := \frac{\frac{1}{\strongconvex\stepsize m} + 2 \lipschitz (1 + 2 \sqrt{\sparsity} \overlap) (\stepsize + \overlap \strongconvex \stepsize^2)}{1 - 2 \lipschitz (1 + 2 \sqrt{\sparsity} \overlap) (\stepsize + \overlap \strongconvex \stepsize^2)} < 1 \,.
	\end{equation}
	Then the inconsistent read iterates of \KROMAGNON\ converge in expectation at a geometric rate, i.e.
	\begin{equation}
		\E f(\tilde x_k) - f(x^*) \leq \theta^t (f(x_0) -f(x^*)) \, ,
	\end{equation}
	where $\tilde x_k$ is the initial iterate for epoch $k$, which is obtained by choosing uniformly at random among the inconsistent read iterates from the previous epoch.
\end{theorem}
This result is similar to the theorem given in the original \SVRG\ paper~\citep{svrg}.
Indeed, if we remove the asynchronous part (i.e. if we set $\overlap = 0$), we get exactly the same rate and condition.
It also has the same form as the one given in~\citet{smola}, which was derived for dense asynchronous \SVRG\ in the easier setting of consistent read and writes (and in the flawed ``after write'' framework), and gives essentially the same conditions on $\stepsize$ and $m$.

In the canonical example presented in most \SVRG\ papers, with $\kappa = n$, $m=\mathcal{O}(n)$ and $\stepsize = \nicefrac{1}{10\lipschitz}$, \SVRG\ obtains a convergence rate of $0.5$.
\citet{smola} get the same rate by setting $\stepsize = \nicefrac{1}{20 \max(1, \sqrt{\sparsity} \overlap) \lipschitz}$ and $m = \mathcal{O}\big(n (1 + \sqrt{\sparsity} \overlap)\big)$.
Following the same line of reasoning (setting $\stepsize = \nicefrac{1}{20 \max(1, \sqrt{\sparsity} \overlap) \lipschitz}$, $\overlap = \mathcal{O}(n)$ and $\theta = 0.5$ and computing the resulting condition on $m$), these values for $\stepsize$ and $m$ also give us a convergence rate of $0.5$.
Therefore, as in~\citet{smola}, when $\kappa = n$ we get a linear speedup for $\overlap < \nicefrac{1}{\sqrt{\sparsity}}$ (which can be as big as $\sqrt{n}$ in the degenerate case where no data points share any feature with each other).
Note that this is the same speedup condition as \ASAGA\ in this regime.

\SVRG\ theorems are usually similar to Theorem~\ref{thm:SVRG}, which does not give an optimal step size or epoch size.
This makes the analysis of a parallel speedup difficult, prompting authors to compare rates in specific cases with most parameters fixed, as we have just done.
In order to investigate the speedup and step size conditions more precisely and thus derive a more general theorem, we now give \SVRG\ and \KROMAGNON\ results modeled on Theorem~\ref{thm:convergence}.

\begin{corollary}[Convergence guarantee and rate for serial \SVRG]\label{thm:SVRGconvergence}
	Let $\stepsize = \frac{a}{4\lipschitz}$ for any $a\leq \frac{1}{4}$ and $m = \frac{32 \kappa}{a}$.
	Then \SVRG\ converges geometrically in expectation with a rate factor per gradient computation of at least $\rho(a) = \frac{1}{4} \min\big\{\frac{1}{n}, \frac{a}{64\kappa}\big\}$, i.e.
	\begin{equation}\label{eq:svrgcontraction}
		\E f(\tilde x_k)-f(x^*) \leq (1-\rho)^{k(2m + n)} \,  (f(x_0) -f(x^*)) \qquad \forall k \geq 0 \,.
	\end{equation}
\end{corollary}
Due to \SVRG's special structure, we cannot write $\E f(x_t)-f(x^*) \leq (1-\rho)^t \,  (f(x_0) -f(x^*))$ for all $t \geq 0$.
However, expressing the convergence properties of this algorithm in terms of a rate factor per gradient computation (of which there are $2m + n$ per epoch) makes it easier to compare convergence rates, either to similar algorithms such as \SAGA\ or to its parallel variant \KROMAGNON\ -- and thus to study the speedup obtained by parallelizing \SVRG.

Compared to \SAGA, this result is very close.
The main difference is that the additional hyper-parameter $m$ has to be set and requires knowledge of $\strongconvex$.
This illustrates the fact that \SVRG\ is not adaptive to local strong convexity, whereas both \SAGA\ and Hofmann's \SVRG\ are.

\begin{corollary}[Simplified convergence guarantee and rate for \KROMAGNON]\label{thm:kromagnon}
	Let
	\begin{equation}
	\begin{aligned}\label{eq:conditionSVRG}
	a^*(\overlap) = \frac{1}{4 (1 + 2\sqrt{\sparsity} \overlap) (1 + \frac{\overlap}{16 \kappa})}
	\,.
	\end{aligned}
	\end{equation}
	For any step size $\stepsize = \frac{a}{4\lipschitz}$ with $a \leq a^*(\overlap)$ and $m = \frac{32 \kappa}{a}$, \KROMAGNON\ converges geometrically in expectation with a rate factor per gradient computation of at least $\rho(a) = \frac{1}{4} \min\big\{\frac{1}{n}, \frac{a}{64\kappa}\big\}$, i.e.
	\begin{equation}
	\E f(\tilde x_k)-f(x^*) \leq (1-\rho)^{k(2m + n)} \, (f(x_0) -f(x^*)) \qquad \forall k \geq 0 \,.
	\end{equation}
\end{corollary}
This result is again quite close to Corollary~\ref{thm:SVRGconvergence}	derived in the serial case.
We see that the maximum step size is divided by an additional $(1 + 2 \overlap \sqrt{\sparsity})$ term, while the convergence rate is the same.
Comparing the rates and the maximum allowable step sizes in both settings give us the sufficient condition on $\overlap$ to get a linear speedup.

\begin{corollary}[Speedup condition]\label{thm:bigdataSVRG}
Suppose $\overlap \leq \mathcal{O}(n)$ and $\overlap \leq \mathcal{O}({\scriptstyle \frac{1}{\sqrt{\sparsity}}} \max\{1,\frac{n}{\kappa} \})$. If $n \geq \kappa$, also suppose $\overlap \leq \sqrt{n \sparsity^{-\nicefrac{1}{2}}}$. Then using the step size $\stepsize = \nicefrac{a^*(\overlap)}{L}\,$ from~\eqref{eq:conditionSVRG}, \KROMAGNON\ converges geometrically with rate factor $\Omega( \min\{\frac{1}{n}, \frac{1}{\kappa}\})$ (similar to \SVRG), and is thus linearly faster than its sequential counterpart up to a constant factor.
\end{corollary}
This result is almost the same as \ASAGA, with the additional condition that $\overlap \leq \mathcal{O}(\sqrt{n})$ in the well-conditioned regime.
We see that in this regime \KROMAGNON\ can also get the same rate as \SVRG\ even without sparsity, which had not been observed in previous work.

Furthermore, one has to note that $\overlap$ is generally smaller for \KROMAGNON\ than for \ASAGA\ since it is reset to $0$ at the beginning of each new epoch (where all cores are synchronized once more).

\paragraph{\textit{Comparison to related work.}}
\begin{itemize}[topsep=1mm, itemsep=-1mm]
	\item Corollary~\ref{thm:SVRGconvergence} provides a rate of convergence \emph{per gradient computation} for \SVRG, contrary to most of the literature on this algorithm~\citep[including the seminal paper][]{svrg}. 
	This result allows for easy comparison with \SAGA\ and other algorithms (in contrast,~\citealt{s2gd} is more involved).
	\item In contrast to the \SVRG\ analysis from~\citet[Thm. 2]{smola}, our proof technique handles inconsistent reads and a non-uniform processing speed across $f_i$'s. 
	While Theorem~\ref{thm:SVRG} is similar to theirs, Corollary~\ref{thm:SVRGconvergence} and~\ref{thm:kromagnon} are more precise results. 
	They enable a finer analysis of the speedup conditions (Corollary~\ref{thm:bigdataSVRG}) -- including the possible speedup without sparsity regime.
	\item In contrast to the \KROMAGNON\ analysis from~\citet[Thm. 14]{mania}, Theorem~\ref{thm:SVRG} gives a better dependence on the condition number in the rate ($1/\kappa$ vs. $1/\kappa^2$ for them) and on the sparsity (they get $\overlap \leq \mathcal{O}(\sparsity^{\nicefrac{-1}{3}})$), while we remove their gradient bound assumption.
	Our results are state-of-the-art for \SVRG\ (contrary to theirs) and so our speedup comparison is more meaningful.
	Finally, Theorem~\ref{thm:SVRG} gives convergence guarantees on $\hat{x}_t$ \emph{during} the algorithm, whereas they only bound the error for the ``last'' iterate $x_T$.
\end{itemize}

\subsection{Proof of Theorem~\ref{thm:SVRG}}\label{sec:proofSVRG}
We now give a detailed outline of the proof. Its full details can be found in Appendix~\ref{apx:SVRG}.

Our proof technique begins as our \ASAGA\ analysis.
In particular, Properties~\ref{independence},~\ref{prop:unbiased},~\ref{eventconst} are also verified for \KROMAGNON\footnote{Note that similarly to \ASAGA, the \KROMAGNON\ algorithm which we analyze reads the parameters first and then samples. This is necessary in order for Property~\ref{independence} to be verified at $r=t$, although not practical when it comes to actual implementation.}, and as in our \ASAGA\ analysis, we make Assumption~\ref{boundedoverlap} (bounded overlaps).
Consequently, the basic recursive contraction inequality~\eqref{eq:RecursiveIneq1} and Lemma~\ref{lma:1} also hold.
However, when we derive the equivalent of Lemma~\ref{lma:suboptgt}, we get a slightly different form, which prompts a difference in the rest of the proof technique.

\begin{lemma} [Suboptimality bound on $\E \|g_t\|^2$]\label{lma:suboptgtSVRG}
	\begin{equation}\label{gtboundsvrg}
		\E\|g_t\|^2
		\leq 4\lipschitz e_t
		+ 4\lipschitz \tilde e_k
		\qquad \forall k \geq 0, km\leq t \leq (k+1)m\, ,
	\end{equation}
	where $\tilde e_k :=  \E f(\tilde x_k) - f(x^*)$ and $\tilde x_k$ is the initial iterate for epoch $k$.
\end{lemma}
We give the proof in Appendix~\ref{apx:SVRGlemma}. To derive both terms, we use the same technique as for the first term of Lemma~\ref{lma:suboptgt}.
Although this is a much simpler result than Lemma~\ref{lma:suboptgt} in the case of \ASAGA, two key differences prevent us from reusing the same Lyapunov function proof technique.
First, the $e_0$ term in Lemma~\ref{lma:suboptgt} is replaced by $\tilde e_k$ which depends on the epoch number.
Second, this term is not multiplied by a geometrically decreasing quantity, which means the $-2\stepsize e_0$ term is not sufficient to cancel out all of the $e_0$ terms coming from subsequent inequalities.
To solve this issue, we go to more traditional \SVRG\ techniques.

The rest of the proof is as follows:
\begin{enumerate}
	\item By substituting Lemma~\ref{lma:suboptgtSVRG} into Lemma~\ref{lma:1}, we get a master contraction inequality~\eqref{eq:masterSVRG} in terms of $a_{t+1}$, $a_t$ and $e_u, u\leq t$.
	\item As in~\citet{svrg}, we sum the master contraction inequality over a whole epoch, and then use the same randomization trick~\eqref{eq:randomtrick} to relate $(e_t)_{km \leq t \leq(k+1)m -1}$ to $\tilde e_k$.
	\item We thus obtain a contraction inequality between $\tilde e_k$ and $\tilde e_{k-1}$, which finishes the proof for Theorem~\ref{thm:SVRG}.
	\item We then only have to derive the conditions on $\stepsize, \overlap$ and $m$ under which we contractions and compare convergence rates to finish the proofs for Corollary~\ref{thm:SVRGconvergence}, Corollary~\ref{thm:kromagnon} and Corollary~\ref{thm:bigdataSVRG}.
\end{enumerate}

We list the key points below with their proof sketch, and give the detailed proof in Appendix~\ref{apx:SVRG}.

\paragraph{\textit{Master inequality.}}
As in our \ASAGA\ analysis, we apply~\eqref{gtboundsvrg} to the result of Lemma~\ref{lma:1}, which gives us that for all $k \geq 0, km\leq  t \leq (k+1)m -1$ (see Appendix~\ref{apx:SVRGmaster}):

\begin{equation}\label{eq:masterSVRG}
a_{t+1}
\leq (1 - \frac{\stepsize\strongconvex}{2})a_t
+ (4\lipschitz \stepsize^2 C_1 -2\stepsize) e_t
+ 4\lipschitz \stepsize^2 C_2 \sum_{u = \max(km, t-\overlap)}^{t-1} e_u
+ (4\lipschitz \stepsize^2 C_1 + 4 \lipschitz \stepsize^2\overlap C_2) \tilde e_k\,.
\end{equation}

\paragraph{\textit{Contraction inequality.}}
As we previously mentioned, the term in $\tilde e_k$ is not multiplied by a geometrically decreasing factor, so using the same Lyapunov function as for \ASAGA\ cannot work.
Instead, we apply the same method as in the original \SVRG\ paper~\citep{svrg}: we sum the master contraction inequality over a whole epoch.
This gives us (see Appendix~\ref{apx:SVRGContraction}):

\begin{equation}\label{eq:svrgsum}
a_{(k+1)m} \leq a_{km}
+ (4\lipschitz \stepsize^2 C_1 + 4\lipschitz \stepsize^2 \overlap C_2 -2\stepsize) \sum_{t=km}^{(k+1)m-1} e_t
+ m (4\lipschitz \stepsize^2 C_1 + 4 \lipschitz \stepsize^2\overlap C_2) \tilde e_k \,.
\end{equation}

To cancel out the $\tilde e_k$ term, we only have one negative term on the right-hand side of~\eqref{eq:svrgsum}: $-2 \stepsize \sum_{t=km}^{(k+1)m-1} e_t$.
This means we need to relate $\sum_{t=km}^{(k+1)m-1} e_t$ to $\tilde e_k$.
We can do it using the same randomization trick as in~\citet{svrg}: instead of choosing the last iterate of the $k^\mathrm{th}$ epoch as $\tilde x_k$, we pick one of the iterates of the epoch uniformly at random.
This means we get:
\begin{equation}\label{eq:randomtrick}
\tilde e_k = \E f(\tilde x_k) - f(x^*) = \frac{1}{m} \sum_{t=(k-1)m}^{km-1} e_t
\end{equation}
We now have: $\sum_{t=km}^{(k+1)m-1} e_t = m \tilde e_{k+1}$.
Combined with the fact that $a_{km} \leq \frac{2}{\strongconvex} \tilde e_k$ and that we can remove the positive $a_{(k+1)m}$ term from the left-hand-side of~\eqref{eq:svrgsum}, this gives us our final recursion inequality:

\begin{equation}\label{eq:svrgdetailedcontraction}
\big(2\stepsize m - 4\lipschitz \stepsize^2 C_1 m - 4\lipschitz \stepsize^2 \overlap C_2 m\big) \tilde e_{k+1}
\leq
\big(\frac{2}{\strongconvex} + 4\lipschitz \stepsize^2 C_1 m + 4 \lipschitz \stepsize^2\overlap C_2 m\big) \tilde e_k
\end{equation}
Replacing $C_1$ and $C_2$ by their values (defined in~\ref{eq:C12defs}) in~\eqref{eq:svrgdetailedcontraction} directly leads to Theorem~\ref{thm:SVRG}.

\section{\Hogwild\ Analysis}\label{sec:SGD}
In order to show that our improved ``after read'' perturbed iterate framework can be used to revisit the analysis of other optimization routines with correct proofs that do not assume homogeneous computation, we now provide the analysis of the \Hogwild\ algorithm (i.e. asynchronous parallel constant step size \SGD) first introduced in~\citet{hogwild}.

We start by describing \Hogwild\ in Algorithm~\ref{alg:hogwild}, and then give our theoretical convergence and speedup results and their proofs.
Note that our framework allows us to easily remove the classical bounded gradient assumption, which is used in one form or another in most of the literature~\citep{hogwild, taming, mania} -- although it is inconsistent with strong convexity in the unconstrained regime.
This allows for better bounds where the uniform bound on $\|f'_i(x)\|^2$ is replaced by its variance at the optimum.

\begin{figure*}[ttt!]
	\centering
 \begin{minipage}[t]{0.49\textwidth}
 \vspace{\dimexpr\ht\strutbox-\topskip}
   \begin{algorithm}[H]
     \caption{\Hogwild}
     \label{alg:hogwild}
     \begin{algorithmic}[1]
	   \STATE Initialize shared variable $x$
	   \LOOP
	      \STATE $\hat x = $ inconsistent read of $x$
	      \STATE \textcolor{blue}{Sample $i$}  uniformly in $\{1,...,n\}$
	      \STATE Let $S_i$ be $f_i$'s support
    	      \STATE $[\delta x]_{S_i} := -\stepsize f'_i(\hat x)$
        	  \FOR{$v$ {\bfseries in} $S_i$}
        		 \STATE $[x]_v \leftarrow [x]_v + [\delta x]_v$      \hfill // atomic
          \ENDFOR
	   \ENDLOOP
	  \end{algorithmic}
    \end{algorithm}
 \end{minipage}
\end{figure*}

\subsection{Theoretical Results}
We now state the theoretical results of our analysis of \Hogwild\ with inconsistent reads and writes in the ``after read framework''.
We give an outline of the proof in Section~\ref{sec:proofSGD} and its full details in Appendix~\ref{apx:SGD}.
We start with a useful definition.
\begin{definition} Let $\sigma^2 = \Econd\|f'_i(x^*)\|^2$ be the variance of the gradient estimator at the optimum.
\end{definition}
For reference, we start by giving the rate of convergence of serial \SGD~\citep[see e.g.][]{schmidt2014sgd}.
\begin{theorem}[Convergence guarantee and rate of \SGD]\label{thm:convergenceSGDserial}
	Let $a \leq \frac{1}{2}$.
	Then for any step size $\stepsize = \frac{a}{\lipschitz}$, \SGD\ converges in expectation to $b$-accuracy at a geometric rate of at least: $\contraction(a) = \nicefrac{a}{\kappa},$ i.e., $\E \|x_t - x^*\|^2 \leq (1-\rho)^t \|x_0 - x^*\|^2 + b$, where $b = 2 \frac{\stepsize \sigma^2}{\strongconvex}$.
\end{theorem}
As \SGD\ only converges linearly up to a ball around the optimum, to make sure we reach $\epsilon$-accuracy, it is necessary that $\frac{2 \stepsize \sigma^2}{\strongconvex} \leq \epsilon$, i.e. $\stepsize \leq \frac{\epsilon \strongconvex}{2\sigma^2}$.
All told, in order to get linear convergence to $\epsilon$-accuracy, serial \SGD\ requires $\stepsize \leq \min \big\{\frac{1}{2\lipschitz}, \frac{\epsilon \strongconvex}{2\sigma^2}\big\}$. The proof can be found in Appendix~\ref{apx:sgdsgd}.

\begin{theorem}[Convergence guarantee and rate of \Hogwild]\label{thm:convergenceSGD}
Let
\begin{equation}\label{eq:conditionSGD}
a^*(\overlap) := \frac{1}{5 \left(1+ 2 \overlap  \sqrt \sparsity \right) \xi(\kappa, \sparsity, \overlap)} \quad
\begin{aligned} &\text{where } \xi(\kappa, \sparsity, \overlap) := \sqrt{1 + \frac{1}{2 \kappa}  \min\{\frac{1}{\sqrt{\sparsity}}, \overlap\} } \\
& \text{\small{(note that $\xi(\kappa, \sparsity, \overlap) \approx 1$ unless $\kappa < \nicefrac{1}{\sqrt{\sparsity}}  \,\, (\leq \sqrt{n})$)}}.
\end{aligned}
\end{equation}
For any step size $\stepsize = \frac{a}{L}$ with $a \leq \min \big\{a^*(\overlap), \frac{\kappa}{\overlap}\big\}$, the inconsistent read iterates of Algorithm~\ref{alg:hogwild} converge in expectation to $b$-accuracy at a geometric rate of at least: $\contraction(a) = \nicefrac{a}{\kappa},$
i.e., $\E \|\hat x_t - x^*\|^2 \leq (1-\rho)^t (2\|x_0 - x^*\|^2) + b$, where $b = (\frac{8 \stepsize (C_1 + \overlap C_2)}{\strongconvex} + 4 \stepsize^2 C_1 \overlap) \sigma^2$ and $C_1$ and $C_2(\stepsize)$ are defined in~\eqref{eq:C12defs}.
\end{theorem}
Once again this result is quite close to the one obtained for serial \SGD.
Note that we recover this exact condition (up to a small constant factor) if we set $\overlap = 0$, i.e. if we force our asynchronous algorithm to be serial.

The condition $a \leq \frac{\kappa}{\overlap}$ is equivalent to $\stepsize \strongconvex \overlap \leq 1$ and should be thought of as a condition on $\overlap$.
We will see that it is always verified in the regime we are interested in, that is the linear speed-up regime (where more stringent conditions are imposed on $\overlap$).

We now investigate the conditions under which \Hogwild\ is linearly faster than \SGD.
Note that to derive these conditions we need not only compare their respective convergence rates, but also the size of the ball around the optimum to which both algorithms converge.
These quantities are provided in Theorems~\ref{thm:convergenceSGDserial} and~\ref{thm:convergenceSGD}.

\begin{corollary}[Speedup condition]\label{thm:bigdataSGD}
	Suppose $\overlap = \mathcal{O}({\min\{\frac{1}{\sqrt{\sparsity}}}, \kappa\})$.
	Then for any step size $\stepsize \leq \frac{a^*(\overlap)}{\lipschitz} = \mathcal{O}(\frac{1}{\lipschitz})$ (i.e., any allowable step size for \SGD), \Hogwild\ converges geometrically to a ball of radius $r_h = \mathcal{O}(\frac{\stepsize \sigma^2}{\strongconvex})$ with rate factor $\contraction = \frac{\stepsize \strongconvex}{2}$ (similar to \SGD), and is thus linearly faster than its sequential counterpart up to a constant factor.

	Moreover, a universal step size of $\Theta(\frac{1}{L})$ can be used for \Hogwild\ to be adaptive to local strong convexity with a similar rate to \SGD\ (i.e., knowledge of $\kappa$ is not required).
\end{corollary}
If $\stepsize = \mathcal{O}(\nicefrac{1}{\lipschitz})$, \Hogwild\ obtains the same convergence rate as \SGD\ and converges to a ball of equivalent radius.
Since the maximum step size guaranteeing linear convergence for \SGD\ is also $\mathcal{O}(\nicefrac{1}{\lipschitz})$, \Hogwild\ is linearly faster than \SGD\ for any reasonable step size -- under the condition that $\overlap = \mathcal{O}({\min\{\frac{1}{\sqrt{\sparsity}}}, \kappa\})$.
We also remark that since $\stepsize \leq \nicefrac{1}{\lipschitz}$ and $\overlap \leq \kappa$, we have $\stepsize \strongconvex \overlap \leq 1$, which means the condition $a \leq \frac{\kappa}{\overlap}$ is superseded by $a\leq a^*(\overlap)$ in Theorem~\ref{thm:convergenceSGD}.

We note that the condition on $\overlap$ is much more restrictive if the condition number is small than for \ASAGA\ and \KROMAGNON.
This can be explained by the fact that both \SAGA\ and \SVRG\ have a composite rate factor which is not directly proportional to the step size.
As a result, in the well-conditioned setting these algorithms enjoy a range of step sizes that all give the same contraction rate.
This allows their asynchronous variants to use smaller step sizes while maintaining linear speedups.
\SGD, on the other hand, has a rate factor that is directly proportional to its step size, hence the more restrictive condition on $\overlap$.

\paragraph{\textit{Function values results.}}
Our results are derived directly on iterates, that is, we bound the distance between $\hat x_t$ and $x^*$.
We can easily obtain results on function values to bound ${\E f(\hat x_t) - f(x^*)}$ by adapting the classical smoothness inequality:\footnote{See e.g.~\citet{bachandm}.}
${\E f(x_t) - f(x^*) }\leq {\frac{\lipschitz}{2} \E \|x_t - x^*\|^2}$ to the asynchronous parallel setting.

\paragraph{\textit{Convergence to $\epsilon$-accuracy.}}
As noted in~\citet{mania}, for our algorithm to converge to $\epsilon$-accuracy for some $\epsilon > 0$, we require an additional bound on the step size to make sure that the radius of the ball to which we converge is small enough.
For \SGD, this means using a step size $\stepsize = \mathcal{O}(\frac{\epsilon \strongconvex}{\sigma^2})$ (see Appendix~\ref{apx:sgdsgd}).
We can also prove that under the conditions that $\overlap = \mathcal{O}(\frac{1}{\sqrt{\sparsity}})$ and $\stepsize \strongconvex \overlap \leq 1$, \Hogwild\ requires the same bound on the step size to converge to $\epsilon$-accuracy (see Appendix~\ref{apx:SGDcorollary}).

If $\epsilon$ is small enough, the active upper bound on the step size is $\stepsize = \mathcal{O}(\frac{\epsilon \strongconvex}{\sigma^2})$ for both algorithms.
In this regime, we obtain a relaxed condition on $\overlap$ for a linear speedup.
The condition $\overlap \leq \kappa$ which came from comparing maximum allowable step sizes is removed.
Instead, we enforce $\stepsize \strongconvex \overlap \leq 1$, which gives us the weaker condition $\overlap = \mathcal{O}(\frac{\sigma^2}{\epsilon \strongconvex^2})$.
Our condition on the overlap is then: $\overlap = \mathcal{O}({\min\{\frac{1}{\sqrt{\sparsity}}}, \frac{\sigma^2}{\epsilon \strongconvex^2}\})$.
We see that this is similar to the condition obtained by~\citet[Theorem~4]{mania} in their \Hogwild\ analysis, although we have the variance at the optimum $\sigma^2$ instead of a squared global bound on the gradient.

\paragraph{\textit{Comparison to related work.}}
\begin{itemize}[topsep=1mm, itemsep=-1mm]
	\item We give the first convergence analysis for \Hogwild\ with no assumption on a global bound on the gradient ($M$). This allows us to replace the usual dependence in $M^2$ by a term in $\sigma^2$ which is potentially significantly smaller. This means improved upper bounds on the step size and the allowed overlap.
	\item We obtain the same condition on the step size for linear convergence to $\epsilon$-accuracy of \Hogwild\ as previous analysis for serial \SGD~\citep[e.g.][]{srebro} -- given $\overlap \leq \nicefrac{1}{\stepsize \strongconvex}$.
	\item In contrast to the \Hogwild\ analysis from~\citet{hogwild, taming}, our proof technique handles inconsistent reads and a non-uniform processing speed across $f_i$'s.
	Further, Corollary~\ref{thm:bigdataSGD} gives a better dependence on the sparsity than in~\citet{hogwild}, where $\overlap \leq \mathcal{O}(\sparsity^{-\nicefrac{1}{4}})$, and does not require various bounds on the gradient assumptions.
	\item In contrast to the \Hogwild\ analysis from~\citet[Thm. 3]{mania}, removing their gradient bound assumption enables us to get a (potentially) significantly better upper bound condition on $\overlap$ for a linear speedup.
	We also give our convergence guarantee on $\hat{x}_t$ \emph{during} the algorithm, whereas they only bound the error for the ``last'' iterate $x_T$.
\end{itemize}

\subsection{Proof of Theorem~\ref{thm:convergenceSGD} and Corollary~\ref{thm:bigdataSGD}} \label{sec:proofSGD}
Here again, our proof technique begins as our \ASAGA\ analysis, with Properties~\ref{independence},~\ref{prop:unbiased},~\ref{eventconst} also verified for \Hogwild\footnote{Once again, in our analysis the \Hogwild\ algorithms reads the parameters before sampling, so that Property~\ref{independence} is verified for $r=t$.}.
As in our \ASAGA\ analysis, we make Assumption~\ref{boundedoverlap}.
Consequently, the basic recursive contraction inequality~\eqref{eq:RecursiveIneq1} and Lemma~\ref{lma:1} also hold.
As for \KROMAGNON, the proof diverges when we derive the equivalent of Lemma~\ref{lma:suboptgt}.

\begin{lemma} [Suboptimality bound on $\E \|g_t\|^2$]\label{lma:suboptgtSGD}
	For all $t \geq 0$,
	\begin{equation}\label{gtboundsgd}
	\E\|g_t\|^2
	\leq 4\lipschitz e_t
	+ 2 \sigma^2\, .
	\end{equation}
\end{lemma}
We give the proof in Appendix~\ref{apx:SGDlemma}. To derive both terms, we use the same technique as for the first term of Lemma~\ref{lma:suboptgt}.
This result is simpler than both Lemma~\ref{lma:suboptgt} (for \ASAGA) and Lemma~\ref{lma:suboptgtSVRG} (for \KROMAGNON).
The second term in this case does not even vanish as $t$ grows.
This reflects the fact that constant step size \SGD\ does not converge to the optimum but rather to a ball around it.
However, this simpler form allows us to simply unroll the resulting master inequality to get our convergence result.

The rest of the proof is as follows:
\begin{enumerate}
	\item By substituting Lemma~\ref{lma:suboptgtSGD} into Lemma~\ref{lma:1}, we get a master contraction inequality~\eqref{eq:masterSGD} in terms of $a_{t+1}$, $a_t$, $(e_u, u\leq t)$ and $\sigma^2$.
	\item We then unroll this master inequality and cleverly regroup terms to obtain a contraction inequality~\eqref{eq:sgdunroll} between $a_t$, $a_0$ and $\sigma^2$.
	\item By using that $\|\hat x_t - x^*\|^2 \leq 2 a_t + 2\|\hat x_t - x_t\|^2$, we obtain a contraction inequality directly on the ``real'' iterates~\citep[as opposed to the ``virtual'' iterates as in][]{mania}, subject to a maximum step size condition on $\stepsize$. This finishes the proof for Theorem~\ref{thm:convergenceSGD}.
	\item Finally, we only have to derive the conditions on $\stepsize$ and $\overlap$ under which \Hogwild\ converges with a similar convergence rate to a ball with a similar radius than serial \SGD\ to finish the proof for Corollary~\ref{thm:bigdataSGD}.
\end{enumerate}

We list the key points below with their proof sketch, and give the detailed proof in Appendix~\ref{apx:SGD}.

\paragraph{\textit{Master inequality.}}
As in our \ASAGA\ analysis, we plug~\eqref{gtboundsgd} in Lemma~\ref{lma:1}, which gives us that (see Appendix~\ref{apx:SGDtheorem}):

\begin{equation}\label{eq:masterSGD}
a_{t+1}
\leq (1 - \frac{\stepsize\strongconvex}{2})a_t
+ (4\lipschitz \stepsize^2 C_1 -2\stepsize) e_t
+ 4\lipschitz \stepsize^2 C_2 \sum_{u = (t-\overlap)_+}^{t-1} e_u
+ 2 \stepsize^2 \sigma^2 (C_1 + \overlap C_2)\,.
\end{equation}

\paragraph{\textit{Contraction inequality on $x_t$.}}
As we previously mentioned, the term in $\sigma^2$ does not vanish so we cannot use either our \ASAGA\ or our \KROMAGNON\ proof technique.
Instead, we unroll Equation~\eqref{eq:masterSGD} all the way to $t=0$.
This gives us (see Appendix~\ref{apx:SGDtheorem}):

\begin{equation}\label{eq:sgdunroll}
a_{t+1} \leq (1 - \frac{\stepsize \strongconvex}{2})^{t+1} a_{0}
+ (4\lipschitz \stepsize^2 C_1 + 8\lipschitz \stepsize^2 \overlap C_2 -2\stepsize) \sum_{u=0}^{t} (1 - \frac{\stepsize \strongconvex}{2})^{t -u} e_u
+ \frac{4 \stepsize \sigma^2}{\strongconvex} ( C_1 + \overlap C_2) \,.
\end{equation}

\paragraph{\textit{Contraction inequality on $\hat x_t$.}}
We now use that $\|\hat x_t - x^*\|^2 \leq 2 a_t + 2\|\hat x_t - x_t\|^2$ together with our previous bound~\eqref{eq:hatxtbound}. Together with~\eqref{eq:sgdunroll}, we get (see Appendix~\ref{apx:SGDtheorem}):
\begin{align}\label{eq:sgdxhat}
\E \|\hat x_t - x^*\|^2 \leq
(1 - \frac{\stepsize \strongconvex}{2})^{t+1} 2 a_{0}
&+ \Big(\frac{8 \stepsize (C_1 + \overlap C_2)}{\strongconvex} + 4 \stepsize^2 C_1 \overlap \Big) \sigma^2
\nonumber \\
&+ (24\lipschitz \stepsize^2 C_1 + 16\lipschitz \stepsize^2 \overlap C_2 -4\stepsize) \sum_{u=0}^{t} (1 - \frac{\stepsize \strongconvex}{2})^{t -u} e_u\,.
\end{align}
To get our final contraction inequality, we need to safely remove all the $e_u$ terms, so we enforce $16\lipschitz \stepsize^2 C_1 + 16\lipschitz \stepsize^2 \overlap C_2 -4\stepsize \leq 0$.
This leads directly to Theorem~\ref{thm:convergenceSGD}.

\paragraph{\textit{Convergence rate and ball-size comparison.}}
To prove Corollary~\ref{thm:bigdataSGD}, we simply show that under the condition $\overlap = \mathcal{O}({\min\{\frac{1}{\sqrt{\sparsity}}}, \kappa\})$, the biggest allowable step size for \Hogwild\ to converge linearly is $\mathcal{O}(\nicefrac{1}{\lipschitz})$, as is also the case for \SGD; and that the size of the ball to which both algorithms converge is of the same order.
The proof is finished by remarking that for both algorithms, the rates of convergence are directly proportional to the step size.

\section{Empirical Results}\label{results}
We now present the results of our experiments.
We first compare our new sequential algorithm, Sparse \SAGA, to its existing alternative, \SAGA\ with lagged updates and to the original \SAGA\ algorithm as a baseline.
We then move on to our main results, the empirical comparison of \ASAGA , \KROMAGNON\ and \Hogwild.
Finally, we present additional results, including convergence and speedup figures with respect to the number of iteration (i.e. ``theoretical speedups'') and measures on the $\overlap$ constant.

\subsection{Experimental Setup}\label{ImplDetails}
\paragraph{\textit{Models.}}
Although \ASAGA\ can be applied more broadly, we focus on logistic regression, a model of particular practical importance.
The associated objective function takes the following form:
\begin{equation}
\frac{1}{n} \sum_{i=1}^n \log\big(1 + \exp(- b_i a_i^\intercal x)\big) + \frac{\mu}{2} \|x\|^2,
\end{equation}
where $a_i \in \mathbb{R}^d$ and $b_i \in \{-1,+1\}$ are the data samples.

\paragraph{\textit{Data sets.}}
We consider two sparse data sets: RCV1~\citep{RCV1} and URL~\citep{URL}; and a dense one, Covtype~\citep{Covtype}, with statistics listed in the table below.
As in~\citet{SAG}, Covtype is standardized, thus $100\%$ dense.
$\sparsity$ is $\mathcal{O}(1)$ in all data sets, hence not very insightful when relating it to our theoretical results.
Deriving a less coarse sparsity bound remains an open problem.

\begin{table}[ht]
	\caption{Basic data set statistics.}
	\centering
	\begin{tabular}{lcccc}
		\toprule
		{} & $n$ & $d$ & density & $\lipschitz$\\
		\midrule
		{\bf RCV1} & \hfill 697,641 & \hfill 47,236 & \hfill 0.15\% & \hfill 0.25\\
		{\bf URL} & \hfill 2,396,130 & \hfill 3,231,961 & \hfill 0.004\% & \hfill 128.4\\
		{\bf Covtype} & \hfill 581,012 & \hfill 54 & \hfill 100\% & \hfill 48428\\
		\bottomrule
	\end{tabular}
\end{table}

\paragraph{\textit{Hardware and software.}}
Experiments were run on a 40-core machine with 384GB of memory.
All algorithms were implemented in Scala. We chose this high-level language despite its typical 20x slowdown compared to C (when using standard libraries, see Appendix~\ref{scalavsc}) because our primary concern was that the code may easily be reused and extended for research purposes (to this end, we have made all our code available at \url{http://www.di.ens.fr/sierra/research/asaga/}).

\subsection{Implementation Details}
{\bf Regularization.}
Following~\citet{laggedsaga}, the amount of regularization used was set to $\mu=1 / n$.
In each update, we project the gradient of the regularization term (we multiply it by $D_i$ as we also do with the vector $\bar{\alpha}$) to preserve the sparsity pattern while maintaining an unbiased estimate of the gradient.
For squared $\ell_2$, the Sparse \SAGA\ updates becomes:
\begin{equation}
x^+ = x - \gamma (f_i'(x) - \alpha_i + D_i \bar{\alpha} + \mu D_i x).
\end{equation}
{\bf Comparison with the theoretical algorithm.}
The algorithm we used in the experiments is fully detailed in Algorithm~\ref{alg:sagasync}.
There are two differences with Algorithm~\ref{alg:theoretical}.
First, in the implementation we choose $i_t$ at random \textit{before} we read the feature vector $a_{i_t}$.
This enables us to only read the necessary data for a given iteration (i.e.~$[\hat x_t]_{S_i}, [\hat \alpha_i^t], [\bar \alpha^t]_{S_i}$).
Although this violates Property~\ref{independence}, it still performs well in practice.

Second, we maintain $\bar \alpha^t$ in memory.
This saves the cost of recomputing it at every iteration (which we can no longer do since we only read a subset data).
Again, in practice the implemented algorithm enjoys good performance.
But this design choice raises a subtle point: the update is not guaranteed to be unbiased in this setup (see Appendix~\ref{apx:Bias} for more details).

\paragraph{\textit{Step sizes.}}
For each algorithm, we picked the best step size among 10 equally spaced values in a grid, and made sure that the best step size was never at the boundary of this interval. 
For Covtype and RCV1, we used the interval $[\frac{1}{10L}, \frac{10}{L}]$, whereas for URL we used  the interval $[\frac{1}{L}, \frac{100}{L}]$ as it admitted larger step sizes.
It turns out that the best step size was fairly constant for different number of cores for both \ASAGA\ and \KROMAGNON, and both algorithms had similar best step sizes ($0.7$ for RCV1, $0.05$ for URL and $5 \times 10^{-5}$ for Covtype).

\subsection{Comparison of Sequential Algorithms: Sparse \SAGA\ vs Lagged updates}\label{sec:ssagacomp}

\begin{figure*}
	\includegraphics[width=\linewidth]{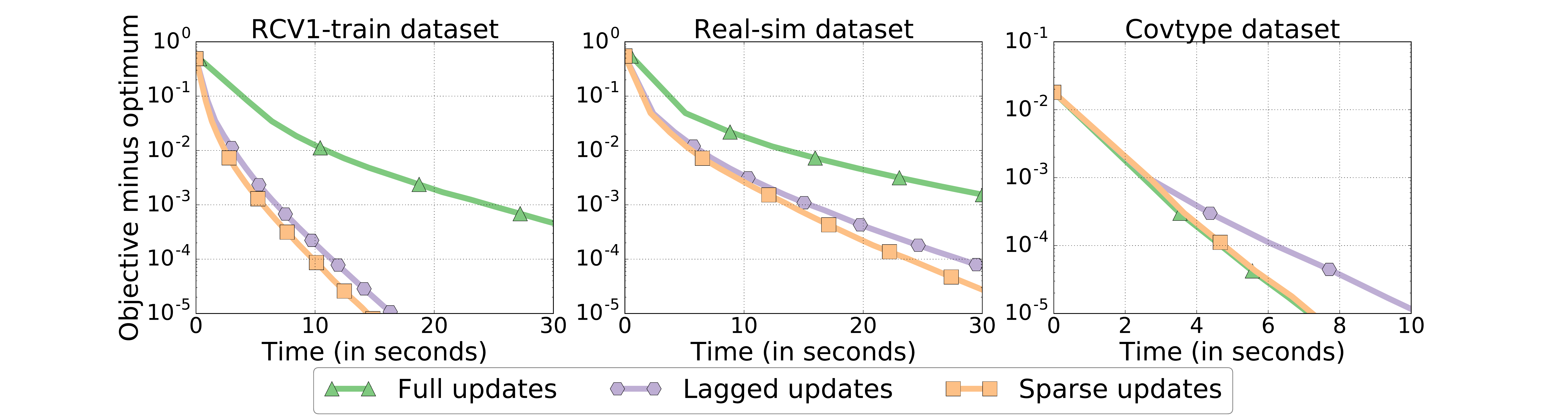}
	\includegraphics[width=\linewidth]{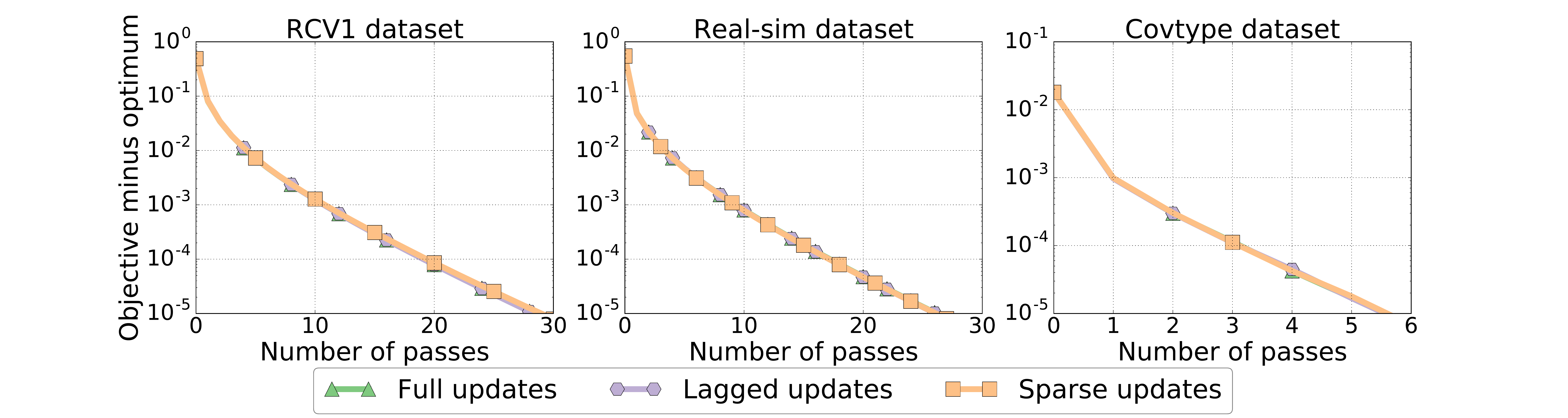}
	\caption{{\bf Lagged vs Sparse \SAGA\ updates}.
		Suboptimality with respect to time for different \SAGA\ update schemes on various data sets. First row: suboptimality as a function of time. Second row: suboptimality as a the number of passes over the data set.
		For sparse data sets (RCV1 and Real-sim), lagged and sparse updates have a lower cost per iteration which result in faster convergence.}\label{fig:fig_1}
\end{figure*}

We compare the Sparse \SAGA\ variant proposed in Section~\ref{scs:sparse_saga} to two other approaches: the naive (i.e., dense) update scheme and the lagged updates implementation described in~\citet{SAGA}.
Note that we use different datasets from the parallel experiments, including a subset of the RCV1 data set and the Realsim data set (see description in Appendix~\ref{apx:datasets}).
Figure~\ref{fig:fig_1} reveals that sparse and lagged updates have a lower cost per iteration than their dense counterpart, resulting in faster convergence for sparse data sets.
Furthermore, while the two approaches had similar convergence in terms of number of iterations, the Sparse \SAGA\ scheme is slightly faster in terms of runtime (and as previously pointed out, sparse updates are better adapted for the asynchronous setting).
For the dense data set (Covtype), the three approaches exhibit similar performance.

\subsection{\ASAGA\ vs. \KROMAGNON\ vs. \Hogwild}\label{sec:asynccomp}
\begin{figure*}[ttt!]
	\centering
	\begin{subfigure}[t]{0.48\linewidth}
		\centering
		\includegraphics[width = 1.05\linewidth]{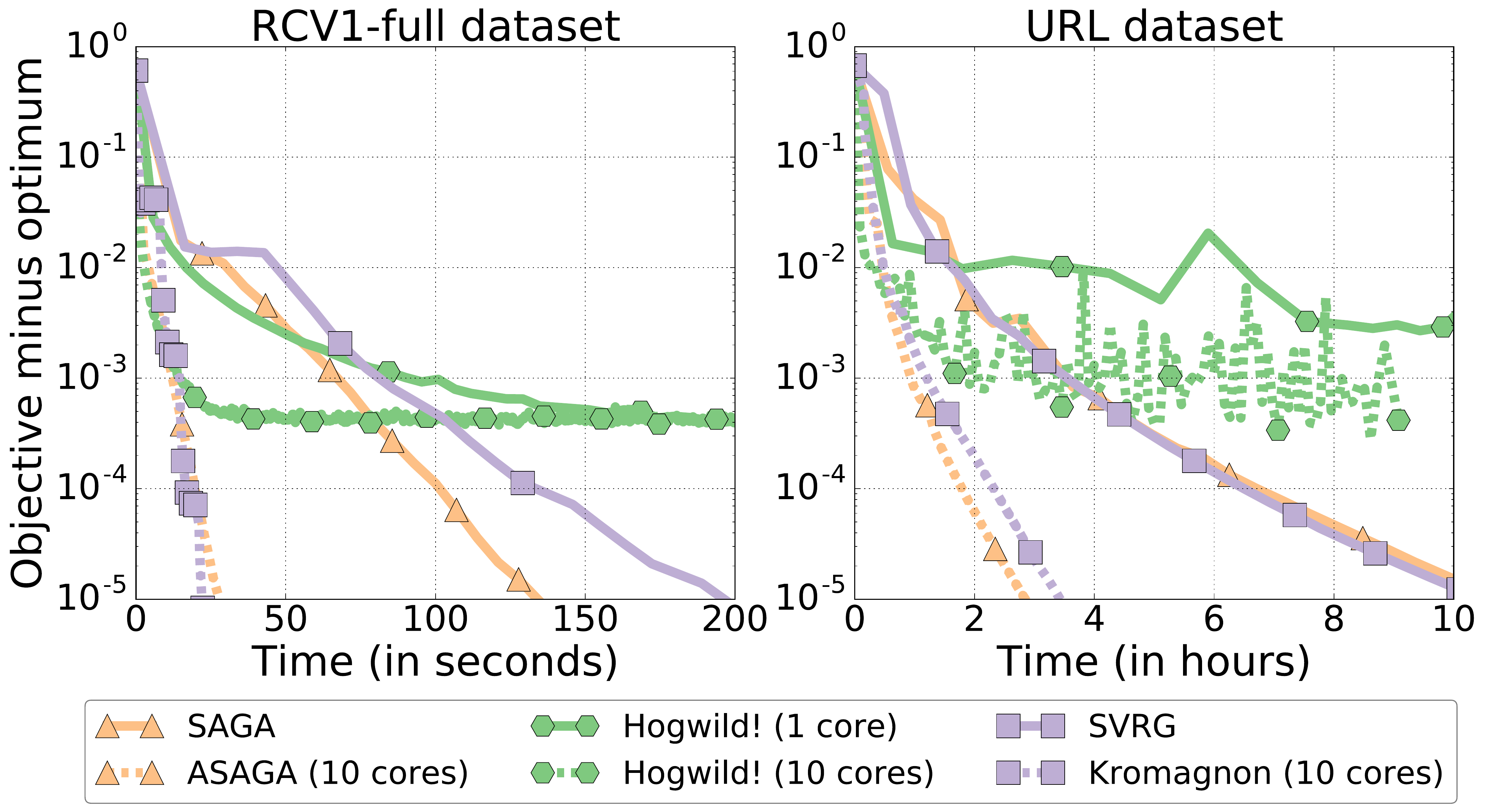}
		\caption{Suboptimality as a function of time.}
		\label{fig:fig_2}
	\end{subfigure}
	\hfill
	\begin{subfigure}[t]{0.48\linewidth}
		\centering
		\includegraphics[width = 1\linewidth]{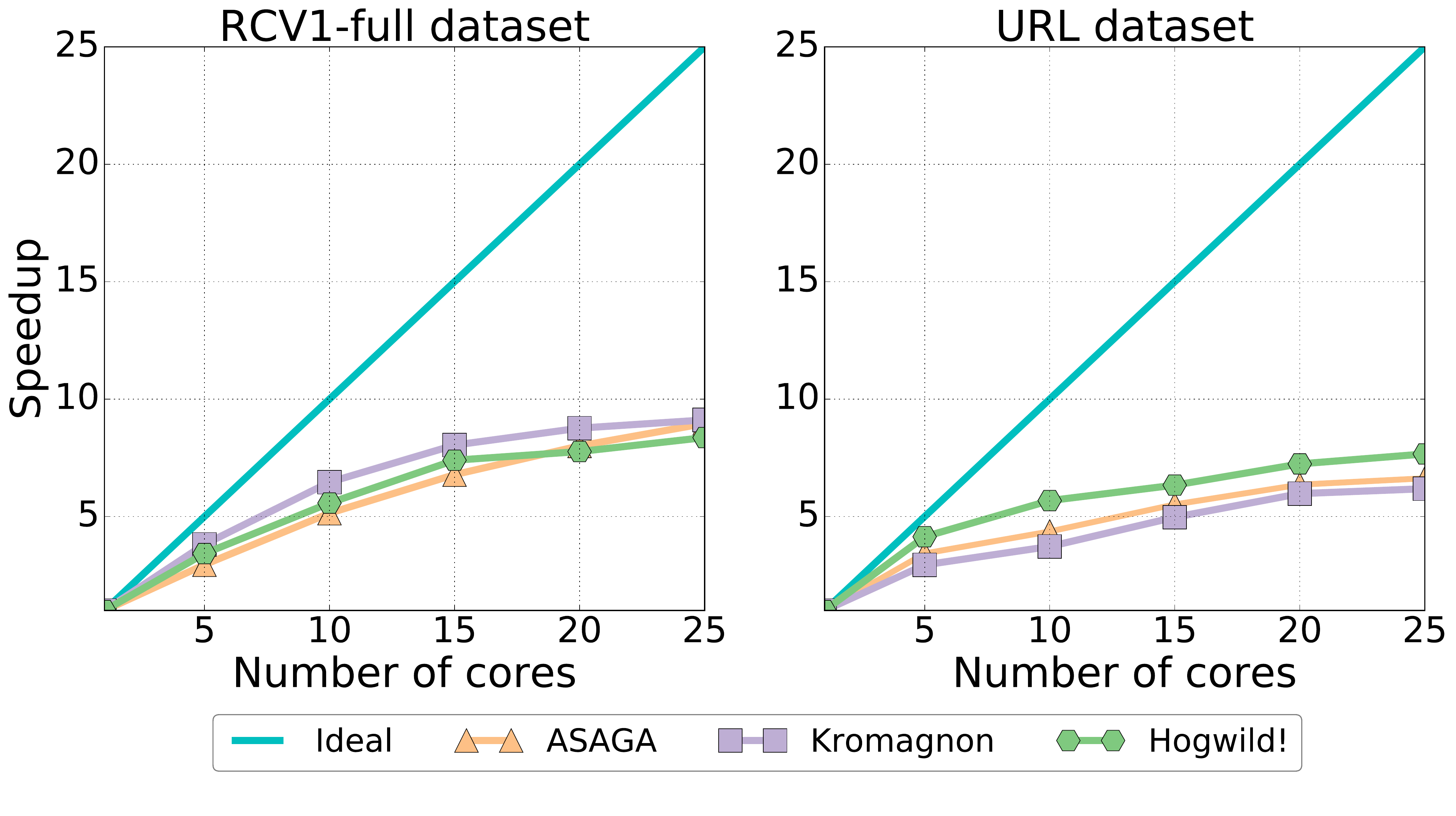}
		\caption{Speedup as a function of the number of cores}
		\label{fig:fig_3}
	\end{subfigure}
	\caption{ {\bf Convergence and speedup for asynchronous stochastic gradient descent methods}.
		We display results for RCV1 and URL. Results for Covtype can be found in Section~\ref{apx:speedup}. }
\end{figure*}
We compare three different asynchronous variants of stochastic gradient methods on the aforementioned data sets: \ASAGA, presented in this work, \KROMAGNON, the asynchronous sparse \SVRG\ method described in~\citet{mania} and \Hogwild~\citep{hogwild}.
Each method had its step size chosen so as to give the fastest convergence (up to a suboptimality of $10^{-3}$ in the special case of \Hogwild).
The results can be seen in Figure~\ref{fig:fig_2}: for each method we consider its asynchronous version with both one (hence sequential) and ten processors.
This figure reveals that the asynchronous version offers a significant speedup over its sequential counterpart.

We then examine the speedup relative to the increase in the number of cores.
The speedup is measured as time to achieve a suboptimality of $10^{-5}$ ($10^{-3}$ for \Hogwild) with one core divided by time to achieve the same suboptimality with several cores, averaged over 3 runs.
Again, we choose step size leading to fastest convergence\footnote{Although we performed grid search on a large interval, we observed that the best step size was fairly constant for different number of cores, and similar for \ASAGA\ and \KROMAGNON.} (see Appendix~\ref{scalavsc} for information about the step sizes).
Results are displayed in Figure~\ref{fig:fig_3}.

As predicted by our theory, we observe linear ``theoretical'' speedups (i.e. in terms of number of iterations, see Section~\ref{apx:speedup}).
However, with respect to running time, the speedups seem to taper off after $20$ cores.
This phenomenon can be explained by the fact that our hardware model is by necessity a simplification of reality.
As noted in~\citet{duchi}, in a modern machine there is no such thing as \textit{shared memory}.
Each core has its own levels of cache (L1, L2, L3) in addition to RAM.
These faster pools of memory are fully leveraged when using a single core.
Unfortunately, as soon as several cores start writing to common locations, cache coherency protocols have to be deployed to ensure that the information is consistent across cores.
These protocols come with computational overheads.
As more and more cores are used, the shared information goes lower and lower in the memory stack, and the overheads get more and more costly.
It may be the case that on much bigger data sets, where the cache memory is unlikely to provide benefits even for a single core (since sampling items repeatedly becomes rare), the running time speedups actually improve.
More experimentation is needed to quantify these effects and potentially increase performance.

\subsection{Effect of Sparsity}
\begin{figure}
	\center \includegraphics[width=0.7\linewidth]{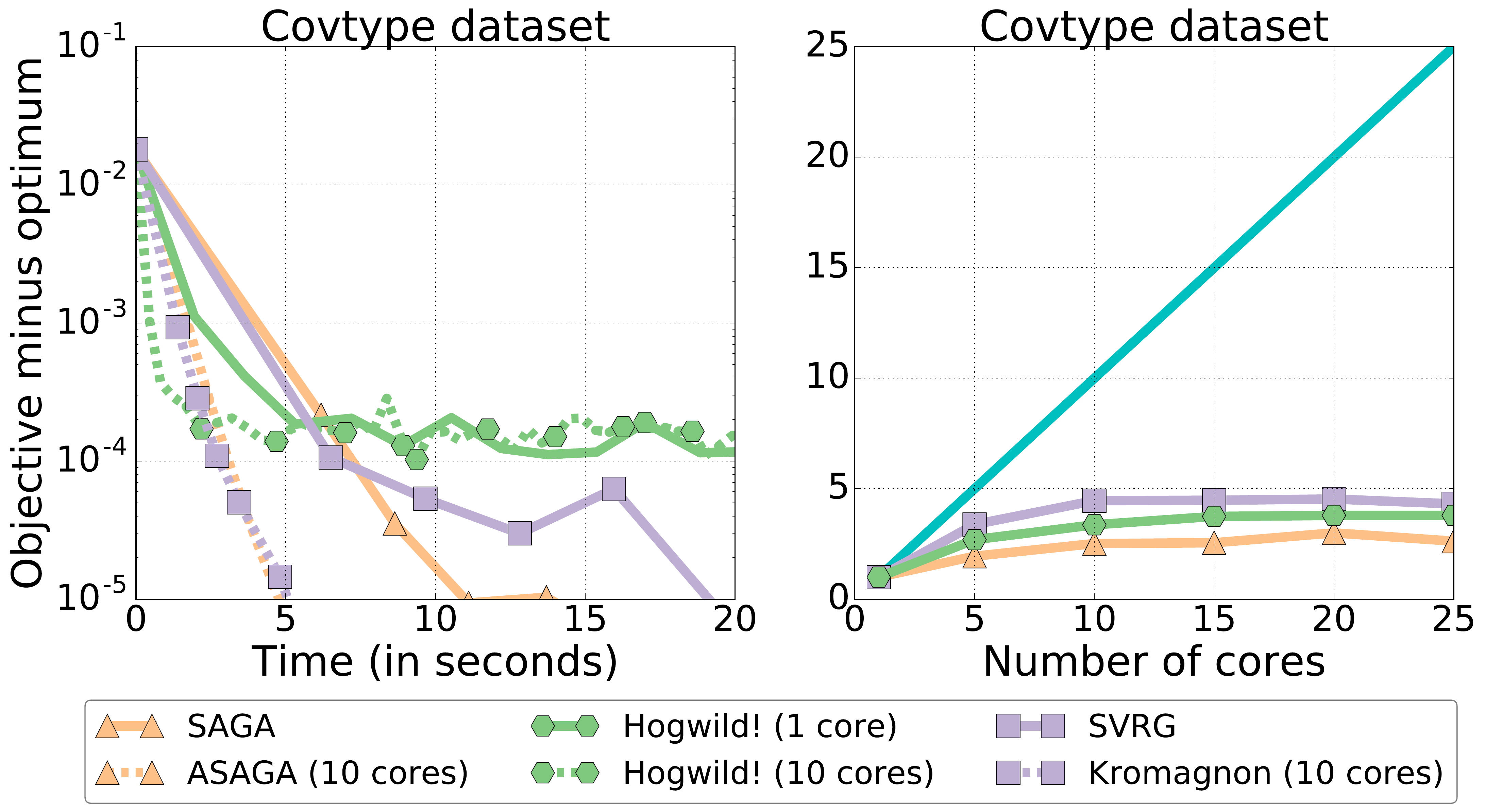}
	\caption{Comparison on the Covtype data set. Left: suboptimality. Right: speedup. The number of cores in the legend only refers to the left plot. }\label{fig:covtype}
\end{figure}

Sparsity plays an important role in our theoretical results, where we find that while it is necessary in the ``ill-conditioned'' regime to get linear speedups, it is not in the ``well-conditioned'' regime.
We confront this to real-life experiments by comparing the convergence and speedup performance of our three asynchronous algorithms on the Covtype data set, which is fully dense after standardization.
The results appear in Figure~\ref{fig:covtype}.

While we still see a significant improvement in speed when increasing the number of cores, this improvement is smaller than the one we observe for sparser data sets.
The speedups we observe are consequently smaller, and taper off earlier than on our other data sets.
However, since the observed ``theoretical'' speedup is linear (see Section~\ref{apx:speedup}), we can attribute this worse performance to higher hardware overhead.
This is expected because each update is fully dense and thus the shared parameters are much more heavily contended for than in our sparse datasets.

One thing we notice when computing the $\sparsity$ constant for our data sets is that it often fails to capture the full sparsity distribution, being essentially a maximum: for all three data sets, we obtain $\sparsity = \mathcal{O}(1)$.
This means that $\sparsity$ can be quite big even for very sparse data sets.
Deriving a less coarse bound remains an open problem.

\subsection{Theoretical Speedups}\label{apx:speedup}
In the previous experimental sections, we have shown experimental speedup results where suboptimality was a function of the running time.
This measure encompasses both theoretical algorithmic properties and hardware overheads (such as contention of shared memory) which are not taken into account in our analysis.

\begin{figure*}
	\includegraphics[width=\linewidth]{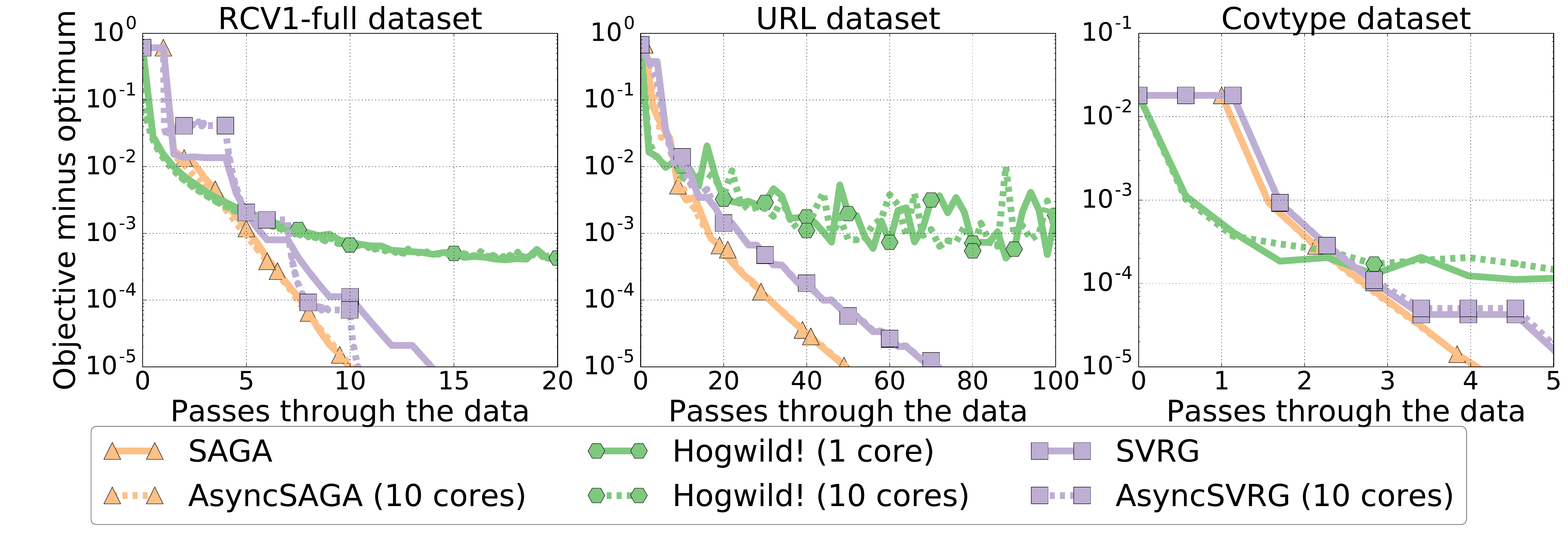}
	\caption{{\bf Theoretical speedups}.
		Suboptimality with respect to number of iterations for \ASAGA, \KROMAGNON\ and \Hogwild\ with 1 and 10 cores.
		Curves almost coincide, which means the theoretical speedup is almost the number of cores $p$, hence linear.}\label{fig:theoretical_speedups}
\end{figure*}

In order to isolate these two effects, we now plot our convergence experiments where suboptimality is a function of the number of iterations; thus, we abstract away any potential hardware overhead.\footnote{To do so, we implement a global counter which is sparsely updated (every $100$ iterations for example) in order not to modify the asynchrony of the system.
	This counter is used only for plotting purposes and is not needed otherwise.
}
The experimental results can be seen in Figure~\ref{fig:theoretical_speedups}.

For all three algorithms and all three data sets, the curves for $1$ and $10$ cores almost coincide, which means that we are indeed in the ``theoretical linear speedup'' regime.
Indeed, when we plotted the amount of iterations required to converge to a given accuracy as a function of the number of cores, we obtained straight  horizontal lines for our three algorithms.

The fact that the speedups we observe in running time are less than linear can thus be attributed to various hardware overheads, including shared variable contention -- the compare-and-swap operations are more and more expensive as the number of competing requests increases -- and cache effects as mentioned in Section~\ref{sec:asynccomp}.

\subsection{A Closer Look at the $\overlap$ Constant} \label{apxD}\label{sec:overlap}
\subsubsection{Theory}\label{sssec:overlaptheory}
In the parallel optimization literature, $\overlap$ is often referred to as a proxy for the number of cores.
However, intuitively as well as in practice, it appears that there are a number of other factors that can influence this quantity.
We will now attempt to give a few qualitative arguments as to what these other factors might be and how they relate to $\overlap$.

\paragraph{\textit{Number of cores.}}
The first of these factors is indeed the number of cores.
If we have $p$ cores, $\overlap \geq p - 1$.
Indeed, in the best-case scenario where all cores have exactly the same execution speed for a single iteration, $\overlap = p -1$.

\paragraph{\textit{Length of an iteration.}}
To get more insight into what $\overlap$ really encompasses, let us now try to define the worst-case scenario in the preceding example.
Consider $2$ cores.
In the worst case, one core runs while the other is stuck.
Then the overlap is $t$ for all $t$ and eventually grows to $+\infty$.
If we assume that one core runs twice as fast as the other, then $\overlap = 2$.
If both run at the same speed, $\overlap = 1$.

It appears then that a relevant quantity is $R$, the ratio between the fastest execution time and the slowest execution time for a single iteration.
We have $\overlap \leq (p-1) R$, which can be arbitrarily bigger than $p$.

There are several factors at play in $R$ itself. These include:
\begin{itemize}
	\item the speed of execution of the cores themselves (i.e. clock time).

	\item the data matrix itself. Different support sizes for $f_i$ means different gradient computation times. If one $f_i$ has support of size $n$ while all the others have support of size $1$ for example, $R$ may eventually become very big.

	\item the length of the computation itself. The longer our algorithm runs, the more likely it is to explore the potential corner cases of the data matrix.
\end{itemize}

The overlap is then upper bounded by the number of cores multiplied by the ratio of the maximum iteration time over the minimum iteration time (which is linked to the sparsity distribution of the data matrix).
This is an upper bound, which means that in some cases it will not really be useful.
For example, if one factor has support size $1$ and all others have support size $d$, the probability of the event which corresponds to the upper bound is exponentially small in $d$.
We conjecture that a more useful indicator could be ratio of the maximum iteration time over the expected iteration time.

To sum up this preliminary theoretical exploration, the $\overlap$ term encompasses much more complexity than is usually implied in the literature.
This is reflected in the experiments we ran, where the constant was orders of magnitude bigger than the number of cores.

\subsubsection{Experimental Results}
In order to verify our intuition about the $\overlap$ variable, we ran several experiments on all three data sets, whose characteristics are reminded in Table~\ref{table:1}.
$\delta_l^i$ is the support size of $f_i$.

\begin{table}[ht]
	\caption{Density measures including minimum, average and maximum support size $\delta_l^i$ of the factors.}
	\centering
	\label{dataset-table2}
	\begin{tabular}{lccccccc}
		\toprule
		{} & $n$ & $d$ & density & $\max(\delta_l^i)$ & $\min(\delta_l^i)$ & $\bar \delta_l$ & $\max(\delta_l^i) / \bar \delta_l$\\
		\midrule
		{\bf RCV1} & \hfill 697,641 & \hfill 47,236 & \hfill 0.15\% & \hfill 1,224 & \hfill 4 & \hfill 73.2 & \hfill 16.7\\
		{\bf URL} & \hfill 2,396,130 & \hfill 3,231,961 & \hfill 0.003\% & \hfill 414 & \hfill 16 & \hfill 115.6 & \hfill 3.58 \\
		{\bf Covtype} & \hfill 581,012 & \hfill 54 & \hfill 100\% & \hfill 12 & \hfill 8 & \hfill 11.88 & \hfill 1.01 \\
		\bottomrule
	\end{tabular}
	\label{table:1}
\end{table}

To estimate $\overlap$, we compute the average overlap over $100$ iterations, i.e. the difference in labeling between the end of the hundredth iteration and the start of the first iteration on, divided by $100$.
This quantity is a lower bound on the actual overlap (which is a maximum, not an average).
We then take its maximum observed value.
The reason why we use an average is that computing the overlap requires using a global counter, which we do not want to update every iteration since it would make it a heavily contentious quantity susceptible of artificially changing the asynchrony pattern of our algorithm.

The results we observe are order of magnitude bigger than $p$, indicating that $\overlap$ can indeed not be dismissed as a mere proxy for the number of cores, but has to be more carefully analyzed.

\begin{figure*}
	\includegraphics[width=\linewidth]{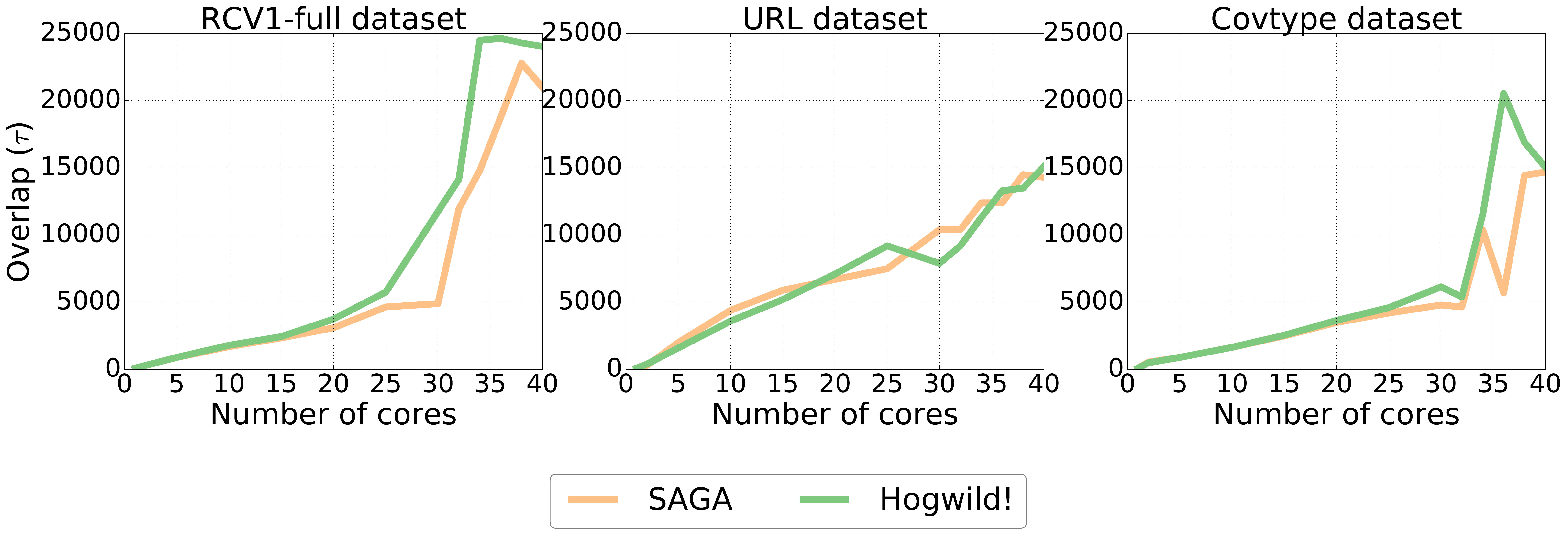}
	\caption{{\bf Overlap}.
		Overlap as a function of the number of cores for both \ASAGA\ and \Hogwild\ on all three data sets.}\label{fig:overlap}
\end{figure*}

First, we plot the maximum observed $\overlap$ as a function of the number of cores (see Figure~\ref{fig:overlap}).
We observe that the relationship does indeed seem to be roughly linear with respect to the number of cores until 30 cores.
After 30 cores, we observe what may be a phase transition where the slope increases significantly.

Second, we measured the maximum observed $\overlap$ as a function of the number of epochs.
We omit the figure since we did not observe any dependency; that is, $\overlap$ does not seem to depend on the number of epochs.
We know that it must depend on the number of iterations (since it cannot be bigger, and is an increasing function with respect to that number for example), but it appears that a stable value is reached quite quickly (before one full epoch is done).

If we allowed the computations to run forever, we would eventually observe an event such that $\overlap$ would reach the upper bound mentioned in Section~\ref{sssec:overlaptheory}, so it may be that $\overlap$ is actually a very slowly increasing function of the number of iterations.

\section{Conclusions and Future Work}
Building on the recently proposed ``perturbed iterate'' framework, we have proposed a novel perspective to clarify an important technical issue present in a large fraction of the recent convergence rate proofs for asynchronous parallel optimization algorithms.
To resolve it, we have introduced a novel ``after read'' framework and demonstrated its usefulness by analyzing three asynchronous parallel incremental optimization algorithms, including \ASAGA, a novel sparse and fully asynchronous variant of the incremental gradient algorithm \SAGA.
Our proof technique accommodates more realistic settings than is usually the case in the literature (such as inconsistent reads and writes and an unbounded gradient); we obtain tighter conditions than in previous work.
In particular, we show that \ASAGA\ is linearly faster than \SAGA\ under mild conditions, and that sparsity is not always necessary to get linear speedups.
Our empirical benchmarks confirm speedups up to 10x.

\citet{laggedsaga} have shown that \SAG\ enjoys much improved performance when combined with non-uniform sampling and line-search.
We have also noticed that our $\sparsityr$ constant (being essentially a maximum) sometimes fails to accurately represent the full sparsity distribution of our data sets.
Finally, while our algorithm can be directly ported to a distributed master-worker architecture, its communication pattern would have to be optimized to avoid prohibitive costs. Limiting communications can be interpreted as artificially increasing the delay, yielding an interesting trade-off between delay influence and communication costs.

These constitute interesting directions for future analysis, as well as a further exploration of the $\overlap$ term, which we have shown encompasses more complexity than previously thought.

\acks{We would like to thank Xinghao Pan for sharing with us their implementation of \KROMAGNON, as well as Alberto Chiappa for spotting a typo in the proof.
This work was partially supported by a Google Research Award and the MSR-Inria Joint Center.
FP acknowledges financial support from the chaire {\em \'Economie des nouvelles donn\'ees} with the {\em data science} joint research initiative with the {\em fonds AXA pour la recherche}.
}

\clearpage
\appendix

\paragraph{\textit{Appendix Outline:}}
\begin{itemize}
\item In Appendix~\ref{apxA}, we adapt the proof from~\citet{qsaga} to prove Theorem~\ref{th1}, our convergence result for serial Sparse \SAGA.
\item In Appendix~\ref{apxB}, we give the complete details for the proof of convergence for~\ASAGA\ (Theorem~\ref{thm:convergence}) as well as its linear speedup regimes  (Corollary~\ref{thm:bigdata}).
\item In Appendix~\ref{apx:SVRG}, we give the full details for the proof of convergence for~\KROMAGNON\ (Theorem~\ref{thm:SVRG}) as well as a simpler convergence result for both \SVRG\ (Corollary~\ref{thm:SVRGconvergence}) and \KROMAGNON\ (Corollary~\ref{thm:kromagnon}) and finally the latter's linear speedup regimes (Corollary~\ref{thm:bigdataSVRG})
\item In Appendix~\ref{apx:SGD}, we give the full details for the proof of convergence for~\Hogwild\ (Theorem~\ref{thm:convergenceSGD}) as well as its linear speedup regimes (Corollary~\ref{thm:bigdataSGD}).
\item In Appendix~\ref{apxC}, we explain why adapting the lagged updates implementation of \SAGA\ to the asynchronous setting is difficult.
\item In Appendix~\ref{apxE}, we give additional details about the data sets and our implementation.
\end{itemize}

\section{Sparse \SAGA\ -- Proof of Theorem~\ref{th1}}\label{apxA}
\paragraph{\textit{Proof sketch for~\citet{qsaga}.}}
As we will heavily reuse the proof technique from~\citet{qsaga}, we start by giving its sketch.

First, the authors combine classical strong convexity and Lipschitz inequalities to derive the following inequality~\citep[Lemma~1]{qsaga}:
\begin{align}\label{eq:sparseappendix}
\Econd \|x^{+} \! - \!x^*\|^2 \leq 
&(1 \! - \! \stepsize\strongconvex) \|x \! -\! x^*\|^2 
+ 2\stepsize^2 \Econd \|\alpha_i - f'_i(x^*)\|^2
+ (4 \stepsize^2 \lipschitz-2\stepsize)\big(f(x) - f(x^*)\big).
\end{align}
This gives a contraction term, as well as two additional terms; $2\stepsize^2 \Econd \|\alpha_i - f'_i(x^*)\|^2$ is a positive variance term, but $(4 \stepsize^2 \lipschitz-2\stepsize)\big(f(x) - f(x^*)\big)$ is a negative suboptimality term (provided $\stepsize$ is small enough).
The suboptimality term can then be used to cancel the variance one.

Second, the authors use a classical smoothness upper bound to control the variance term and relate it to the suboptimality.
However, since the $\alpha_i$ are partial gradients computed at previous time steps, the upper bounds of the variance involve suboptimality at previous time steps, which are not directly relatable to the current suboptimality.

Third, to circumvent this issue, a Lyapunov function is defined to encompass both current and past terms.
To finish the proof,~\citet{qsaga} show that the Lyapunov function is a contraction.

	\paragraph{\textit{Proof outline.}}\label{sparseoutlineappendix}
	Fortunately, we can reuse most of the proof from~\citet{qsaga} to show that Sparse \SAGA\ converges at the same rate as regular \SAGA.
	In fact, once we establish that \citet[Lemma~1]{qsaga} is still verified we are done.
	
	To prove this, we show that the gradient estimator is unbiased, and then derive close variants of equations $(6)$ and $(9)$ in their paper, which we remind the reader of here:
	\begin{align}
	\Econd \|f'_i(x) - \bar \alpha_i\|^2 &\leq 2\Econd \|f'_i(x) - f'_i(x^*)\|^2 \!+\! 2 \Econd \|\bar \alpha_i - f'_i(x^*)\|^2
	\tag*{\citet[Eq. 6]{qsaga}}
	\\
	\Econd \|\bar \alpha_i - f'_i(x^*)\|^2 &\leq \Econd \|\alpha_i - f'_i(x^*)\|^2 \, .
	\tag*{\citet[Eq. 9]{qsaga}}
	\end{align} 
	
	\paragraph{\textit{Unbiased gradient estimator.}}
	We first show that the update estimator is unbiased.
	The estimator is unbiased if:
	\begin{align}\label{eq:apxBias}
	\Econd D_i \bar \alpha = \Econd \alpha_i = \frac{1}{n}\sum_{i=1}^n \alpha_i \, .
	\end{align}	
	We have:
	\begin{align*}
	\Econd D_i \bar \alpha 
	= \frac{1}{n} \sum_{i=1}^n D_i \bar \alpha
	= \frac{1}{n} \sum_{i=1}^n P_{S_i} D \bar \alpha
	= \frac{1}{n} \sum_{i=1}^n \sum_{v \in S_i} \frac{[\bar \alpha]_v e_v}{p_v}
	= \sum_{v=1}^{d} \left( \sum_{i\, | \, v \in S_i} 1 \right) \frac{[\bar \alpha]_v e_v}{n p_v}  \, ,
	\end{align*}
	where $e_v$ is the vector whose only nonzero component is the $v$ component which is equal to $1$.
	
	By definition, $\sum_{i|v \in S_i} 1 = n p_v,$ which gives us Equation~\eqref{eq:apxBias}.
	
	\paragraph{\textit{Deriving~\citet[Equation 6]{qsaga}.}}
	We define $\bar \alpha_i := \alpha_i - D_i\bar \alpha$ (contrary to~\citet{qsaga} where the authors define $\bar \alpha_i := \alpha_i - \bar \alpha$ since they do not concern themselves with sparsity).
	Using the inequality $\|a+b\|^2 \leq 2 \|a\|^2 + 2 \|b\|^2$, we get:
	\begin{align}\label{eq:apxhof6}
	\Econd \|f'_i(x) - \bar \alpha_i\|^2 
	\leq 2\Econd \|f'_i(x) -f'_i(x^*)\|^2 + 2\Econd \|\bar\alpha_i -f'_i(x^*)\|^2,
	\end{align}
	which is our equivalent to~\citet[Eq. 6]{qsaga}, where only our definition of $\bar \alpha_i$ differs.
	
	\paragraph{\textit{Deriving~\citet[Equation 9]{qsaga}.}}
	We want to prove~\citet[Eq. 9]{qsaga}:
	\begin{align}
	\Econd \|\bar \alpha_i -f'_i(x^*)\|^2
	\leq \Econd \|\alpha_i -f'_i(x^*)\|^2 .
	\end{align}
	
	We have:
	\begin{align} \label{eq:alphaiBarVariance}
	\Econd \|\bar \alpha_i -f'_i(x^*)\|^2
	&= \Econd \|\alpha_i -f'_i(x^*)\|^2
	-2\Econd \langle \alpha_i - f'_i(x^*), D_i \bar \alpha \rangle + \Econd\|D_i\bar \alpha\|^2 .
	\end{align}
	
	Let $D_{\neg i} := P_{S_i^c} D$; we then have the orthogonal decomposition $D \alpha = D_i \alpha + D_{\neg i} \alpha$ with $D_i \alpha \perp D_{\neg i} \alpha$, as they have disjoint support. We now use the orthogonality of $D_{\neg i} \alpha$ with any vector with support in $S_i$ to simplify the expression~\eqref{eq:alphaiBarVariance} as follows:
	\begin{align}
	\Econd \langle \alpha_i - f'_i(x^*), D_i\bar \alpha \rangle
	&= \Econd \langle \alpha_i - f'_i(x^*), D_i \bar \alpha + D_{\neg i} \bar \alpha \rangle 
	\tag*{$(\alpha_i - f'_i(x^*) \perp D_{\neg i} \bar \alpha)$}
	\nonumber \\
	&= \Econd \langle \alpha_i - f'_i(x^*), D \bar \alpha \rangle
	\nonumber \\
	&= \langle \Econd \big(\alpha_i - f'_i(x^*)\big), D \bar \alpha \rangle
	\nonumber \\
	&= \langle \Econd \alpha_i, D \bar \alpha \rangle
	\tag*{($f'(x^*) = 0$)}
	\nonumber \\
	&=\bar \alpha^\intercal D \bar \alpha \,  .
	\end{align}
	
	Similarly, 
	\begin{align}
	\Econd \|D_i\bar \alpha\|^2 
	&= \Econd\langle D_i\bar \alpha, D_i\bar \alpha \rangle
	\nonumber \\
	&= \Econd\langle D_i\bar \alpha, D \bar \alpha \rangle
	\tag*{($D_i \bar \alpha \perp D_{\neg i} \bar \alpha$)}
	\nonumber \\
	&= \langle \Econd D_i\bar \alpha, D \bar \alpha \rangle
	\nonumber \\
	&= \bar \alpha^\intercal D \bar \alpha  \, .
	\end{align}
	
	Putting it all together,
	\begin{align}
	\Econd \|\bar \alpha_i -f'_i(x^*)\|^2 
	= \Econd \|\alpha_i -f'_i(x^*)\|^2 - \bar \alpha^\intercal D \bar \alpha
	\leq \Econd \|\alpha_i -f'_i(x^*)\|^2 .
	\end{align}
	
	This is our version of~\citet[Equation 9]{qsaga}, which finishes the proof of~\citet[Lemma~1]{qsaga}. 
	The rest of the proof from~\citet{qsaga} can then be reused without modification to obtain Theorem~\ref{th1}.
	\hfill\BlackBox

\section{\ASAGA\ -- Proof of Theorem~\ref{thm:convergence} and Corollary~\ref{thm:bigdata}}\label{apxB}
\subsection{Initial Recursive Inequality Derivation} \label{app:RecurDerivation}
We start by proving Equation~\eqref{eq:RecursiveIneq1}.
Let $g_t := g(\hat x_t, \hat \alpha^t, i_t)$. From~\eqref{eq:PIupdate}, we get:
\begin{align*}
\|x_{t+1} - x^*\|^2 \nonumber
= \|x_t -\stepsize g_t -x^*\|^2
&= \|x_t -x^*\|^2 + \stepsize^2 \|g_t\|^2 -2\stepsize\langle x_t -x^*,  g_t\rangle
\\
&= \|x_t -x^*\|^2 + \stepsize^2 \|g_t\|^2
	- 2 \stepsize\langle \hat x_t -x^*,  g_t\rangle +2\stepsize\langle \hat x_t -x_t,  g_t\rangle .
\end{align*}

In order to prove Equation~\eqref{eq:RecursiveIneq1}, we need to bound the $- 2 \stepsize\langle \hat x_t -x^*,  g_t\rangle$ term.
Thanks to Property~\ref{independence}, we can write:
\begin{align*}
\E \langle \hat x_t -x^*,  g_t\rangle
= \E \langle \hat x_t -x^*, \Econd g_t \rangle
= \E \langle \hat x_t -x^*,  f'(\hat x_t)\rangle  \, .
\end{align*}

We can now use a classical strong convexity bound as well as a squared triangle inequality to get:
\begin{align}
- \langle \hat x_t -x^*,  f'(\hat x_t)\rangle &\leq - \big(f(\hat x_t) -f(x^*)\big) -\frac{\strongconvex}{2}\|\hat x_t - x^*\|^2
\tag*{(Strong convexity bound)} \nonumber \\
- \|\hat x_t - x^*\|^2 &\leq \|\hat x_t - x_t\|^2 - \frac{1}{2} \|x_t - x^*\|^2
\tag*{($\|a+b\|^2 \leq 2 \|a\|^2 + 2 \|b\|^2$)} \nonumber \\
- 2 \stepsize \E \langle \hat x_t -x^*,  g_t\rangle &\leq
	- \frac{\stepsize \strongconvex}{2} \E \|x_t - x^*\|^2
	+ \stepsize \strongconvex \E \|\hat x_t - x_t\|^2
	-2 \stepsize \big(\E f(\hat x_t) - f(x^*)\big)  \, .
\end{align}
Putting it all together, we get the initial recursive inequality~\eqref{eq:RecursiveIneq1}, rewritten here explicitly:
\begin{align}
a_{t+1} \leq
	(1 -\frac{\stepsize \strongconvex}{2}) a_t
	+ \stepsize^2 \E \|g_t\|^2
	+ \stepsize\strongconvex \E\|\hat x_t - x_t\|^2
	+ 2\stepsize \E \langle \hat x_t -x_t,  g_t\rangle
	-2\stepsize e_t  \, ,
\end{align}
where $a_t := \E \|x_t - x^*\|^2$ and $e_t := \E f(\hat x_t) - f(x^*)$.

\subsection{Proof of Lemma~\ref{lma:1} (inequality in terms of $g_t := g(\hat x_{t}, \hat \alpha^t, i_{t})$)} \label{apxB:lma1}
To prove Lemma~\ref{lma:1}, we now bound both $\E\|\hat x_t - x_t\|^2$ and $\E\langle \hat x_t -x_t,  g_t\rangle$ with respect to $(\E\|g_u\|^2)_{u\leq t}$.

\paragraph{\textit{Bounding $\E\langle \hat x_t -x_t,  g_t\rangle$ in terms of $g_u$.}}
\begin{align}
\frac{1}{\stepsize} \E\langle \hat x_t - x_t, g_t \rangle
&= \sum_{u=(t - \overlap)_+}^{t-1} \E \langle G_u^t g_u, g_t \rangle
\tag*{(by Equation~\eqref{eq:async})} \nonumber \\
&\leq \sum_{u=(t - \overlap)_+}^{t-1} \E | \langle g_u, g_t \rangle |
\tag*{($G_u^t$ diagonal matrices with terms in $\{0, 1\}$)}\nonumber \\
&\leq \sum_{u=(t - \overlap)_+}^{t-1} \frac{\sqrt{\sparsity}}{2}(\E\|g_{u}\|^2 + \E\|g_{t}\|^2)
\tag*{(by Proposition~\ref{prop:1})}\nonumber \\
&\leq \frac{\sqrt{\sparsity}}{2} \sum_{u=(t - \overlap)_+}^{t-1}\E\|g_{u}\|^2 + \frac{\sqrt{\sparsity}\overlap}{2}\E\|g_{t}\|^2 .
\end{align}

\paragraph{\textit{Bounding $\E\|\hat x_t - x_t\|^2$ with respect to $g_u$}}
Thanks to the expansion for $\hat x_t - x_t$~\eqref{eq:async}, we get:
\begin{align*}
\|\hat x_t - x_t\|^2
\leq \stepsize^2 \sum_{u, v=(t -\overlap)_+}^{t-1}|\langle G_u^t g_{u}, G_v^t g_{v}\rangle |
\leq \stepsize^2 \sum_{u=(t -\overlap)_+}^{t-1}\|g_{u}\|^2
	+ \stepsize^2 \sum_{\substack{u, v=(t-\overlap)_+ \\u\neq v}}^{t-1} |\langle G_u^t g_{u}, G_v^t g_{v}\rangle |  \, .
\end{align*}
Using~\eqref{sparseproduct} from Proposition~\ref{prop:1}, we have that for $u \neq v$:

\begin{equation} \label{supersparse}
\E |\langle G_u^t g_{u}, G_v^t g_{v}\rangle |
\leq \E |\langle g_{u}, g_{v}\rangle |
\leq \frac{\sqrt{\sparsity}}{2}(\E\|g_{u}\|^2 + \E\|g_{v}\|^2) \, .
\end{equation}
By taking the expectation and using~\eqref{supersparse}, we get:
\begin{align}\label{eq:hatxtbound}
\E\|\hat x_t - x_t\|^2
&\leq \stepsize^2 \sum_{u=(t-\overlap)_+}^{t-1}\E\|g_{u}\|^2
	+ \stepsize^2 \sqrt{\sparsity}(\overlap-1)_+ \sum_{u=(t-\overlap)_+}^{t-1}\E\|g_{u}\|^2
\nonumber \\
&= \stepsize^2 \big(1+\sqrt{\sparsity}(\overlap-1)_+ \big)\sum_{u=(t-\overlap)_+}^{t-1}\E\|g_{u}\|^2 
\nonumber \\
&\leq \stepsize^2 \big(1+\sqrt{\sparsity}\overlap \big)\sum_{u=(t-\overlap)_+}^{t-1}\E\|g_{u}\|^2 .
\end{align}
We can now rewrite~\eqref{eq:RecursiveIneq1} in terms of $\E\|g_t\| ^2$, which finishes the proof for Lemma~\ref{lma:1} (by introducing $C_1$ and $C_2$ as specified by~\ref{eq:C12defs} in Lemma~\ref{lma:1}):
\begin{align}
a_{t+1} &\leq
	(1 - \frac{\stepsize\strongconvex}{2}) a_t
	- 2\stepsize e_t
	+ \stepsize^2 \E\|g_t\|^2
	+ \stepsize^3 \strongconvex(1+\sqrt{\sparsity}\overlap)\sum_{u=(t-\overlap)_+}^{t-1}\E\|g_{u}\|^2
\nonumber \\
	&\quad + \stepsize^2 \sqrt{\sparsity}\sum_{u=(t-\overlap)_+}^{t-1}\E\|g_{u}\|^2
	+ \stepsize^2 \sqrt{\sparsity}\overlap\E\|g_t\|^2
\nonumber \\
&\leq (1 - \frac{\stepsize\strongconvex}{2}) a_t
	- 2\stepsize e_t
	+ \stepsize^2 C_1 \E\|g_t\|^2
	+ \stepsize^2 C_2 \sum_{u=(t-\overlap)_+}^{t-1}\E\|g_{u}\|^2 .
\end{align}
\hfill\BlackBox

\subsection{Proof of Lemma~\ref{lma:suboptgt} (suboptimality bound on $\E \|g_t\|^2$)} \label{apxB:lma2}
We now derive our bound on $g_t$ with respect to suboptimality.
From Appendix~\ref{apxA}, we know that:
\begin{align}
\E\|g_t\|^2
&\leq 2 \E \|f'_{i_t}(\hat x_t)-f'_{i_t}(x^*)\|^2
	+ 2 \E \|\hat \alpha_{i_t}^t - f'_{i_t}(x^*)\|^2 \label{eq:classicsaga}
\\
\E \|f'_{i_t}(\hat x_t)-f'_{i_t}(x^*)\|^2
&\leq 2\lipschitz\big(\E f(\hat x_t) - f(x^*)\big)
= 2\lipschitz e_t  \, . \label{eq:classicsaga2}
\end{align}
\textbf{N. B.: In the following, $i_t$ is a random variable picked uniformly at random in $\{1,...,n\}$, whereas $i$ is a fixed constant.}

We still have to handle the $\E \|\hat \alpha_{i_t}^t - f'_{i_t}(x^*)\|^2$ term and express it in terms of past suboptimalities.
We know from our definition of $t$ that $i_t$ and $\hat x_u$ are independent $\forall u<t$.
Given the ``after read'' global ordering, $\Econd$ -- the expectation on $i_t$ conditioned on $\hat x_t$ and all ``past" $\hat x_u$ and $i_u$ -- is well defined, and we can rewrite our quantity as:
\begin{align*}
\E \|\hat \alpha_{i_t}^t - f'_{i_t}(x^*)\|^2
= \E \big( \Econd \|\hat \alpha_{i_t}^t - f'_{i_t}(x^*)\|^2 \big)
&= \E \frac{1}{n} \sum_{i=1}^n  \|\hat \alpha_i^t - f'_i(x^*)\|^2
\\
&=  \frac{1}{n} \sum_{i=1}^n \E \|\hat \alpha_i^t - f'_i(x^*)\|^2 .
\end{align*}

Now, with $i$ fixed, let $u_{i,l}^t$ be the time of the iterate last used to write the $[\hat \alpha_{i}^t]_l$ quantity, i.e. $[\hat \alpha_{i}^t]_l = [f'_{i}(\hat x_{u_{i, l}^t})]_l$.
We know\footnote{In the case where $u=0$, one would have to replace the partial gradient with $\alpha_i^0$. We omit this special case here for clarity of exposition.} that $0 \leq u_{i,l}^t \leq t - 1$.
To use this information, we first need to split $\hat \alpha_i$ along its dimensions to handle the possible inconsistencies among them:
\begin{align*}
\E \|\hat \alpha_i^t - f'_i(x^*)\|^2
=
\E \sum_{l=1}^d \big([\hat \alpha_i^t]_l-[f'_i(x^*)]_l\big)^2
=
\sum_{l=1}^d \E \Big[ \big([\hat \alpha_i^t]_l-[f'_i(x^*)]_l\big)^2 \Big] .
\end{align*}
This gives us:
\begin{align}
\E \|\hat \alpha_i^t - f'_i(x^*)\|^2
&=
\sum_{l=1}^d \E \Big[ \big(f'_{i}(\hat x_{u_{i, l}^t})_l-f'_i(x^*)_l\big)^2 \Big] \notag \\
&=
\sum_{l=1}^d\E \Big[\sum_{u=0}^{t-1} \ind_{\{u_{i, l}^t = u\}} \big(f'_i(\hat x_u)_l-f'_i(x^*)_l\big)^2 \Big] \notag \\
&=
\sum_{u=0}^{t-1} \sum_{l=1}^d\E \Big[\ind_{\{u_{i, l}^t = u\}} \big(f'_i(\hat x_u)_l-f'_i(x^*)_l\big)^2 \Big]
\label{eq:IndicatorsAppearance}.
\end{align}

We will now rewrite the indicator so as to obtain independent events from the rest of the equality.
This will enable us to distribute the expectation.
Suppose $u>0$ ($u=0$ is a special case which we will handle afterwards). $\{u_{i, l}^t = u\}$ requires two things:
\begin{enumerate}
\item at time $u$, $i$ was picked uniformly at random,
\item (roughly) $i$ was not picked again between $u$ and $t$.
\end{enumerate}
We need to refine both conditions because we have to account for possible collisions due to asynchrony.
We know from our definition of $\overlap$ that the $t^\mathrm{th}$ iteration finishes before at $t + \overlap + 1$, but it may still be unfinished by time $t + \overlap$.
This means that we can only be sure that an update selecting $i$ at time $v$ has been written to memory at time $t$ if $v \leq t -\overlap -1$.
Later updates may not have been written yet at time $t$.
Similarly, updates before $v = u + \overlap +1$ may be overwritten by the $u^\mathrm{th}$ update so we cannot infer that they did not select $i$. From this discussion, we conclude that $u_{i, l}^t = u$ implies that $i_v \neq i$ for all $v$ between $u+\overlap+1$ and $t-\overlap-1$, though it can still happen that $i_v = i$ for $v$ outside this range.

Using the fact that $i_u$ and $i_v$ are independent for $v \neq u$, we can thus upper bound the indicator function appearing in~\eqref{eq:IndicatorsAppearance} as follows:
\begin{align}\label{eq:indicatrices}
\ind_{\{u_{i,l}^t = u\}}
\leq \ind_{\{i_u=i\}}
	 \ind_{\{i_v \neq i\ \forall v\ \text{s.t.}\ u+\overlap+1 \leq v \leq t-\overlap-1\}} .
\end{align}
This gives us:
\begin{align}
\E \Big[ &\ind_{\{u_{i,l}^t = u\}} \big(f'_i(\hat x_u)_l-f'_i(x^*)_l\big)^2 \Big]
\nonumber \\
&\leq \E \Big[\ind_{\{i_u=i\}}
	\ind_{\{i_v \neq i\ \forall v\ \text{s.t.}\ u+\overlap+1 \leq v \leq t-\overlap-1\}} \big(f'_i(\hat x_u)_l-f'_i(x^*)_l\big)^2\Big]
\nonumber \\
&\leq P\{i_u=i\}
	P\{i_v \neq i\ \forall v\ \text{s.t.}\ u+\overlap+1 \leq v \leq t-\overlap-1\}
	\E\big(f'_i(\hat x_u)_l-f'_i(x^*)_l\big)^2
\tag*{($i_v \perp\!\!\! \perp \hat x_u, \forall v \geq u$)} \nonumber \\
&\leq \frac{1}{n}(1 - \frac{1}{n})^{(t-2\overlap-u -1)_+}
	\E\big(f'_i(\hat x_u)_l-f'_i(x^*)_l\big)^2 .
\end{align}
Note that the third line used the crucial independence assumption $i_v \perp\!\!\! \perp \hat x_u, \forall v \geq u$ arising from our ``After Read'' ordering.
Summing over all dimensions $l$, we then get:
\begin{align}
\E \Big[ \ind_{\{u_{i,l}^t = u\}} \|f'_i(\hat x_u)-f'_i(x^*)\|^2 \Big]
\leq \frac{1}{n}(1 - \frac{1}{n})^{(t-2\overlap-u-1)_+}
	\E \|f'_i(\hat x_u)-f'_i(x^*)\|^2 .
\end{align}
So now:
\begin{align}
\E \|\hat \alpha_{i_t}^t - f'_{i_t}(x^*)\|^2 - \lambda \tilde e_0
&\leq \frac{1}{n}\sum_{i=1}^n \sum_{u=1}^{t-1} \frac{1}{n}(1 - \frac{1}{n})^{(t-2\overlap-u-1)_+} \E \|f'_i(\hat x_{u}) - f'_i(x^*)\|^2
\nonumber \\
&= \sum_{u=1}^{t-1} \frac{1}{n}(1 - \frac{1}{n})^{(t-2\overlap-u -1)_+} \frac{1}{n}\sum_{i=1}^n \E \|f'_i(\hat x_{u}) - f'_i(x^*)\|^2
\nonumber \\
&= \sum_{u=1}^{t-1} \frac{1}{n}(1 - \frac{1}{n})^{(t-2\overlap-u -1)_+} \E \Big(\Econd \|f'_{i_u}(\hat x_{u}) - f'_{i_u}(x^*)\|^2 \Big)
\tag*{($i_u \perp\!\!\!\perp \hat x_u$)} \nonumber \\
&\leq \frac{2\lipschitz}{n} \sum_{u=1}^{t-1} (1 - \frac{1}{n})^{(t-2\overlap-u -1)_+} e_u
\tag*{(by Equation~\ref{eq:classicsaga2})} \nonumber \\
&= \frac{2\lipschitz}{n} \sum_{u=1}^{(t-2\overlap -1)_+} (1 - \frac{1}{n})^{t-2\overlap-u -1} e_u
	+ \frac{2\lipschitz}{n} \sum_{u=\max(1, t-2\overlap)}^{t-1} e_u  \, .\label{eq:51}
\end{align}

Note that we have excluded $\tilde e_0$ from our formula, using a generic $\lambda$ multiplier.
We need to treat the case $u=0$ differently to bound $\ind_{\{u_{i,l}^t = u\}}$.
Because all our initial $\alpha_i$ are initialized to a fixed $\alpha_i^0$, $\{u_i^t = 0\}$ just means that $i$ has not been picked between $0$ and $t-\overlap -1$, i.e. $\{i_v \neq i\ \forall\ v\ \text{s.t.}\ 0 \leq v \leq t - \overlap -1\}$.
This means that the $\ind_{\{i_u=i\}}$ term in~\eqref{eq:indicatrices} disappears and thus we lose a $\frac{1}{n}$ factor compared to the case where $u>1$.

Let us now evaluate $\lambda$.
We have:
\begin{align}\label{eq:52}
\E \Big[\ind_{\{u_i^t = 0\}} \|\alpha_i^0-f'_i(x^*)\|^2 \Big]
&\leq \E \Big[\ind_{\{i_v \neq i\ \forall\ v\ \text{s.t.}\ 0 \leq v \leq t - \overlap -1\}} \|\alpha_i^0-f'_i(x^*)\|^2 \Big]
\nonumber \\
&\leq P\{i_v \neq i\ \forall\ v\ \text{s.t.}\ 0 \leq v \leq t - \overlap -1\} \E \|\alpha_i^0-f'_i(x^*)\|^2
\nonumber \\
&\leq (1 - \frac{1}{n})^{(t-\overlap)_+} \E \|\alpha_i^0-f'_i(x^*)\|^2 .
\end{align}
Plugging~\eqref{eq:51} and~\eqref{eq:52} into~\eqref{eq:classicsaga}, we get Lemma~\ref{lma:suboptgt}:
\begin{align}\label{eq:53}
\E\|g_t\|^2
\leq 4\lipschitz e_t
	+ \frac{4\lipschitz}{n} \sum_{u=1}^{t-1} (1 - \frac{1}{n})^{(t-2\overlap-u -1)_+} e_u
	+ 4\lipschitz (1 - \frac{1}{n})^{(t-\overlap)_+} \tilde e_0 \,  ,
\end{align}
where we have introduced $\tilde e_0 = \frac{1}{2\lipschitz} \E\|\alpha_i^0 - f'_i(x^*)\|^2$.
Note that in the original \SAGA\ algorithm, a batch gradient is computed to set the $\alpha_i^0 = f'_i(x_0)$.
In this setting, we can write Lemma~\ref{lma:suboptgt} using only $\tilde e_0 \leq e_0$ thanks to~\eqref{eq:classicsaga2}.
In the more general setting where we initialize all $\alpha_i^0$ to a fixed quantity, we cannot use~\eqref{eq:classicsaga2} to bound $\E\|\alpha_i^0 - f'_i(x^*)\|^2$ which means that we have to introduce $\tilde e_0$.

\subsection{Lemma~\ref{lma:suboptgt} for \AHSVRG}\label{apxB:ahsvrg}
In the simpler case of \AHSVRG\ as described in~\ref{apx:SVRGext}, we have a slight variation of \eqref{eq:IndicatorsAppearance}:
\begin{equation}
\E \|f'_i(x_0^t) - f'_i(x^*)\|^2
=
\sum_{u=0}^{t-1} \E \| \ind_{\{u_i^t = u\}} \big(f'_i(\hat x_u) - f'_i(x^*)\big) \|^2
\label{eq:IndicatorsAppearanceSVRG}.
\end{equation}
Note that there is no sum over dimensions in this case because the full gradient computations and writes are synchronized (so the reference gradient is consistent).

As in Section~\ref{apxB:lma2}, we can upper bound the indicator $\ind_{\{u_i^t = u\}}$.
Now, $\{u_{i, l}^t = u\}$ requires two things: first, the next $B$ variable sampled after the $u^\mathrm{th}$ update, $\tilde B_u$\footnote{We introduce this quantity because the iterations where full gradients are computed do not receive a time label since they do not correspond to updates to the iterates.}, was $1$; second, $B$ was $0$ for every update between $u$ and $t$ (roughly).
Since the batch step is fully synchronized, we do not have to worry about updates from the past overwriting the reference gradient (and the iterates $x_u$ where we compute the gradient contains all past updates because we have waited for every core to finish its current update).

However, updating the state variable $s$ to $1$ once a $B = 1$ variable is sampled is not atomic.
So it is possible to have iterations with time label bigger than $u$ and that still use an older reference gradient for their update\footnote{Conceivably, another core could start a new iteration, draw $B = 1$ and try to update $s$ to $1$ themselves. This is not an issue since the operation of updating $s$ to $1$ is idempotent. Only one reference gradient would be computed in this case.}.
Fortunately, we can consider the state update as any regular update to shared parameters.
As such, Assumption~\ref{boundedoverlap} applies to it.
This means that we can be certain that the reference gradient has been updated for iterations with time label $v \geq u + \overlap + 1$.

This gives us:
\begin{align}\label{lmaAHSVRG}
\E \| \ind_{\{u_i^t = u\}} \big(f'_i(\hat x_u) - f'_i(x^*)\big) \|^2
&\leq \E \Big[\ind_{\{\tilde B_u=1\}}
\ind_{\{B_v = 0 \ \forall v\ \text{s.t.}\ u + 1 \leq v \leq t-\overlap-1\}} \|f'_i(\hat x_u) - f'_i(x^*)\|^2 \Big]
\nonumber \\
&\leq \frac{1}{n} \big(1 - \frac{1}{n}\big)^{(t - \overlap - u - 1)_+} \E \|f'_i(\hat x_u) - f'_i(x^*)\|^2 \, .
\end{align}
This proves Lemma~\ref{lma:suboptgt} for \AHSVRG\ (while we actually have a slightly better exponent, ($t - \overlap - u - 1)_+$, we can upperbound it by the term in Lemma~\ref{lma:suboptgt}).
Armed with this result, we can finish the proof of Theorem~\ref{thm:convergence} for \AHSVRG\ in exactly the same manner as for \ASAGA.
By remarking that the cost to get to iteration $t$ (including computing reference batch gradients) is the same in the sequential and parallel version, we see that our analysis for Corollary~\ref{thm:bigdata} for \ASAGA\ also applies for \AHSVRG, so both algorithms obey the same convergence and speedup results.

\subsection{Master Inequality Derivation}\label{apxB:master}
Now, if we combine the bound on $\E\|g_{t}\|^2$ which we just derived (i.e. Lemma~\ref{lma:suboptgt}) with Lemma~\ref{lma:1}, we get:
\begin{equation}
\begin{aligned}
a_{t+1}
\leq &(1 - \frac{\stepsize\strongconvex}{2}) a_t
	- 2\stepsize e_t
\\
	&+ 4\lipschitz\stepsize^2 C_1e_t
	+ \frac{4\lipschitz\stepsize^2 C_1}{n} \sum_{u=1}^{t-1} (1 - \frac{1}{n})^{(t-2\overlap-u -1)_+} e_u
	+ 4\lipschitz \stepsize^2 C_1(1 - \frac{1}{n})^{(t-\overlap)_+} \tilde e_0
\\
	&+4\lipschitz\stepsize^2 C_2\sum_{u=(t-\overlap)_+}^{t-1} e_u
	+4\lipschitz\stepsize^2 C_2 \sum_{u=(t-\overlap)_+}^{t-1} (1 - \frac{1}{n})^{(u - \overlap)_+} \tilde e_0
\\
	&+ \frac{4\lipschitz\stepsize^2 C_2}{n} \sum_{u=(t-\overlap)_+}^{t-1} \sum_{v=1}^{u-1} (1-\frac{1}{n})^{(u - 2\overlap - v -1)_+}e_v  \, .
\end{aligned}
\end{equation}
If we define $H_t := \sum_{u=1}^{t-1} (1 - \frac{1}{n})^{(t-2\overlap-u-1)_+} e_u$, then we get:
\begin{equation}\label{eq:master}
\begin{aligned}
a_{t+1}
\leq &(1 - \frac{\stepsize\strongconvex}{2}) a_t
	- 2\stepsize e_t
\\
	&+ 4\lipschitz\stepsize^2 C_1 \big(e_t  + (1 - \frac{1}{n})^{(t-\overlap)_+} \tilde e_0 \big)
	+ \frac{4\lipschitz\stepsize^2 C_1}{n} H_t
\\
	&+4\lipschitz\stepsize^2 C_2\sum_{u=(t-\overlap)_+}^{t-1} (e_u +  (1 - \frac{1}{n})^{(u - \overlap)_+} \tilde e_0 \big)
	+\frac{4\lipschitz\stepsize^2 C_2}{n} \sum_{u=(t-\overlap)_+}^{t-1} H_u  \, ,
\end{aligned}
\end{equation}
which is the master inequality~\eqref{master}.

\subsection{Lyapunov Function and Associated Recursive Inequality}\label{apxB:lyapunov}
We define $\lyapunov_t := \sum_{u=0}^t (1-\contraction)^{t-u}a_u$ for some target contraction rate $\contraction < 1$ to be defined later.
We have:
\begin{align}
\lyapunov_{t+1}
&= (1 - \contraction)^{t+1}a_0
	+ \sum_{u=1}^{t+1}(1 - \contraction)^{t+1-u}a_u
= (1 - \contraction)^{t+1}a_0
	+ \sum_{u=0}^t(1 - \contraction)^{t-u}a_{u+1} \,  .
\end{align}
We now use our new bound on $a_{t+1}$,~\eqref{eq:master}:
\begin{align}\label{eq:rutMaster}
\lyapunov_{t+1}
&\leq (1 - \contraction)^{t+1}a_0
	+ \sum_{u = 0}^t(1 - \contraction)^{t-u} \Big[
		(1 - \frac{\stepsize\strongconvex}{2}) a_u
		- 2\stepsize e_u
		+ 4\lipschitz\stepsize^2 C_1 \big(e_u  + (1 - \frac{1}{n})^{(u-\overlap)_+} \tilde e_0 \big)
\notag \\
		&\qquad \qquad \qquad \qquad \qquad \qquad \qquad
		+ \frac{4\lipschitz\stepsize^2 C_1}{n} H_u
		+\frac{4\lipschitz\stepsize^2 C_2}{n} \sum_{v=(u-\overlap)_+}^{u-1} H_v
\notag \\
		&\qquad \qquad \qquad \qquad \qquad \qquad \qquad
		+4\lipschitz\stepsize^2 C_2\sum_{v=(u-\overlap)_+}^{u-1} (e_v +  (1 - \frac{1}{n})^{(v - \overlap)_+} \tilde e_0 \big)
	\Big]
\notag \\
&\leq (1 - \contraction)^{t+1}a_0
	+ (1-\frac{\stepsize\strongconvex}{2})\lyapunov_t
\notag \\
	&\qquad \qquad + \sum_{u = 0}^t(1 - \contraction)^{t-u} \Big[
		- 2\stepsize e_u
		+ 4\lipschitz\stepsize^2 C_1 \big(e_u  + (1 - \frac{1}{n})^{(u-\overlap)_+} \tilde e_0 \big)
\notag \\
		&\qquad \qquad \qquad \qquad \qquad \qquad
		+ \frac{4\lipschitz\stepsize^2 C_1}{n} H_u
		+\frac{4\lipschitz\stepsize^2 C_2}{n} \sum_{v=(u-\overlap)_+}^{u-1} H_v
\notag \\
		&\qquad \qquad \qquad \qquad \qquad \qquad
		+4\lipschitz\stepsize^2 C_2\sum_{v=(u-\overlap)_+}^{u-1} (e_v +  (1 - \frac{1}{n})^{(v - \overlap)_+} \tilde e_0 \big)
	\Big] .
\end{align}

We can now rearrange the sums to expose a simple sum of $e_u$ multiplied by factors $r_u^t$:
\begin{align}\label{apx:Lyapunov}
\lyapunov_{t+1} \leq (1 - \contraction)^{t+1}a_0 + (1-\frac{\stepsize\strongconvex}{2})\lyapunov_t + \sum_{u=1}^t r_u^t e_u + r_0^t \tilde e_0  \, .
\end{align}

\subsection{Proof of Lemma~\ref{lma:3} (sufficient condition for convergence for \ASAGA)}\label{apxB:lma3}
We want to make explicit what conditions on $\contraction$ and $\stepsize$ are necessary to ensure that $r_u^t$ is negative for all $u \geq 1$.
Since each $e_u$ is positive, we will then be able to safely drop the sum term from the inequality.
The $r_0^t$ term is a bit trickier and is handled separately.
Indeed, trying to enforce that $r_0^t$ is negative results in a significantly worse condition on $\stepsize$ and eventually a convergence rate smaller by a factor of $n$ than our final result.
Instead, we handle this term directly in the Lyapunov function.

\paragraph{\textit{Computation of $r_u^t$.}}
Let's now make the multiplying factor explicit.
We assume $u \geq 1$.

We split $r_u^t$ into five parts coming from~\eqref{eq:rutMaster}:
\begin{itemize}
\item $r_1$, the part coming from the $-2\stepsize e_u$ terms;
\item $r_2$, coming from $4\lipschitz\stepsize^2 C_1 e_u$;
\item $r_3$, coming from $\frac{4\lipschitz\stepsize^2 C_1}{n} H_u$;
\item $r_4$, coming from $4\lipschitz\stepsize^2 C_2\sum_{v=(u-\overlap)_+}^{u-1} e_v$;
\item $r_5$, coming from $\frac{4\lipschitz\stepsize^2 C_2}{n} \sum_{v=(u-\overlap)_+}^{u-1} H_v$.
\end{itemize}
$r_1$ is easy to derive. Each of these terms appears only in one inequality.
So for $u$ at time $t$, the term is:
\begin{equation}\label{eq:r1}
r_1 = -2 \stepsize (1- \contraction)^{t-u} .
\end{equation}
For much the same reasons, $r_2$ is also easy to derive and is:
\begin{equation}\label{eq:r2}
r_2 = 4\lipschitz\stepsize^2 C_1(1- \contraction)^{t-u} .
\end{equation}
$r_3$ is a bit trickier, because for a given $v > 0$ there are several $H_u$ which contain $e_v$.
The key insight is that we can rewrite our double sum in the following manner:
\begin{align}
\sum_{u=0}^t(1 - \contraction)^{t-u} &\sum_{v=1}^{u-1} (1-\frac{1}{n})^{(u - 2\overlap -v -1)_+} e_v
\nonumber \\
&= \sum_{v=1}^{t-1} e_v \sum_{u=v+1}^{t}(1 - \contraction)^{t-u}  (1-\frac{1}{n})^{(u - 2\overlap -v -1)_+}
\nonumber \\
&\leq \sum_{v=1}^{t-1} e_v \Big[
		\sum_{u=v+1}^{\min(t, v + 2\overlap)}(1 - \contraction)^{t-u}
		+ \sum_{u=v + 2 \overlap +1}^{t}(1 - \contraction)^{t-u}  (1-\frac{1}{n})^{u - 2\overlap -v -1}
	\Big]
\nonumber \\
&\leq \sum_{v=1}^{t-1} e_v \Big[
		2\overlap (1 - \contraction)^{t-v-2 \overlap}
		+ (1 - \contraction)^{t-v-2\overlap-1} \sum_{u=v + 2 \overlap +1}^{t}q^{u - 2\overlap -v -1}
	\Big]
\nonumber \\
&\leq \sum_{v=1}^{t-1} (1 - \contraction)^{t-v} e_v (1 - \contraction)^{-2\overlap -1} \big[
		2\overlap
		+ \frac{1}{1-q}
	\big],
\end{align}
where we have defined:
\begin{equation} \label{eq:qDefinition}
q :=\frac{1-1/n}{1-\contraction}, \quad \text{with the assumption $\contraction < \frac{1}{n}$}  \, .
\end{equation}
Note that we have bounded the $\min(t, v+2\overlap)$ term by $v+2\overlap$ in the first sub-sum, effectively adding more positive terms.

This gives us that at time $t$, for $u$:
\begin{equation}\label{eq:r3}
r_3 \leq \frac{4\lipschitz\stepsize^2 C_1}{n}(1 - \contraction)^{t-u} (1 - \contraction)^{-2\overlap -1} \big[2\overlap + \frac{1}{1-q}\big] .
\end{equation}
For $r_4$ we use the same trick:
\begin{align}
\sum_{u=0}^t (1-\contraction)^{t-u} \sum_{v=(u-\overlap)_+}^{u-1} e_v
&= \sum_{v=0}^{t-1} e_v  \sum_{u=v+1}^{\min(t, v+\overlap)} (1-\contraction)^{t-u}
\nonumber \\
&\leq \sum_{v=0}^{t-1} e_v  \sum_{u=v+1}^{v+\overlap} (1-\contraction)^{t-u}
\leq \sum_{v=0}^{t-1} e_v \overlap (1-\contraction)^{t-v - \overlap} .
\end{align}
This gives us that at time $t$, for $u$:
\begin{equation}\label{eq:r4}
r_4 \leq 4\lipschitz\stepsize^2 C_2(1 - \contraction)^{t-u} \overlap (1-\contraction)^{-\overlap} \,  .
\end{equation}

Finally we compute $r_5$ which is the most complicated term.
Indeed, to find the factor of $e_w$ for a given $w > 0$, one has to compute a triple sum, $\sum_{u = 0}^t(1 - \contraction)^{t-u} \sum_{v=(u-\overlap)_+}^{u-1} H_v$.
We start by computing the factor of $e_w$ in the inner double sum, $\sum_{v=(u-\overlap)_+}^{u-1} H_v$.

\begin{align}
\sum_{v=(u-\overlap)_+}^{u-1} \sum_{w=1}^{v-1}(1-\frac{1}{n})^{(v -2\overlap -w -1)_+} e_w
= \sum_{w=1}^{u-2} e_w \sum_{v=\max(w+1, u-\overlap)}^{u-1} (1-\frac{1}{n})^{(v -2\overlap -w -1)_+}  \, .
\end{align}
Now there are at most $\overlap$ terms for each $e_w$.
If $w \leq u - 3\overlap -1$, then the exponent is positive in every term and it is always bigger than $u -3\overlap -1 -w$, which means we can bound the sum by $\overlap (1-\frac{1}{n})^{u -3\overlap -1 -w}$.
Otherwise we can simply bound the sum by $\overlap$. We get:
\begin{align}
\sum_{v=(u-\overlap)_+}^{u-1} H_v
\leq \sum_{w=1}^{u-2}
	\big[
		\ind_{\{u -3\overlap \leq w \leq u -2\}}\overlap
		+ \ind_{\{w \leq u -3\overlap -1\}}\overlap (1-\frac{1}{n})^{u -3\overlap -1 -w}
	\big] e_w  \, .
\end{align}
This means that for $w$ at time $t$:
\begin{align}\label{eq:r5}
r_5
&\leq \frac{4\lipschitz\stepsize^2 C_2}{n} \sum_{u=0}^t (1 -\contraction)^{t-u}
	\big[
		\ind_{\{u -3\overlap \leq w \leq u -2\}}\overlap
		+ \ind_{\{w \leq u -3\overlap -1\}}\overlap (1-\frac{1}{n})^{u -3\overlap -1 -w}
	\big]
\nonumber \\
&\leq \frac{4\lipschitz\stepsize^2 C_2}{n}
	\Big[
		\sum_{u=w+2}^{\min(t, w +3\overlap)} \overlap (1 -\contraction)^{t-u}
		+ \sum_{u=w +3\overlap +1}^t \overlap (1-\frac{1}{n})^{u -3\overlap -1 -w} (1 -\contraction)^{t-u}
	\Big]
\nonumber \\
&\leq \frac{4\lipschitz\stepsize^2 C_2}{n} \overlap
	\Big[
		(1 -\contraction)^{t-w}(1 -\contraction)^{-3\overlap} 3\overlap
\nonumber \\ &\qquad \qquad \qquad
		+ (1 -\contraction)^{t-w}(1 -\contraction)^{-1 -3\overlap}\sum_{u=w +3\overlap +1}^t (1-\frac{1}{n})^{u -3\overlap -1 -w} (1 -\contraction)^{-u +3\overlap +1 +w}
	\Big]
\nonumber \\
&\leq \frac{4\lipschitz\stepsize^2 C_2}{n} \overlap (1 -\contraction)^{t-w}(1 -\contraction)^{-3\overlap-1}
	\big(
		3\overlap
		+ \frac{1}{(1-q)}
	\big) \,  .
\end{align}
By combining the five terms together (\ref{eq:r1},~\ref{eq:r2},~\ref{eq:r3},~\ref{eq:r4} and~\ref{eq:r5}), we get that $\forall u$ s.t. $1 \leq u \leq t$:
\begin{equation}\label{eq:rut}
\begin{aligned}
r_u^t \leq (1- \contraction)^{t-u}
	\Big[
		&-2 \stepsize
		+4\lipschitz\stepsize^2 C_1
		+\frac{4\lipschitz\stepsize^2 C_1}{n}
			(1 - \contraction)^{-2\overlap -1} \big(
				2\overlap
				+ \frac{1}{1-q}
			\big)
\\
		&+4\lipschitz\stepsize^2 C_2 \overlap(1-\contraction)^{-\overlap}
		+\frac{4\lipschitz\stepsize^2 C_2}{n} \overlap (1 -\contraction)^{-3\overlap -1}
			\big(
				3\overlap
				+ \frac{1}{1-q}
		\big)
	\Big] .
\end{aligned}
\end{equation}

\paragraph{\textit{Computation of $r_0^t$.}}
Recall that we treat the $\tilde e_0$ term separately in Section~\ref{apxB:lma2}.
The initialization of \SAGA\ creates an initial synchronization, which means that the contribution of $\tilde e_0$ in our bound on $\E\|g_t\|^2$~\eqref{eq:53} is roughly $n$ times bigger than the contribution of any $e_u$ for $1 < u < t$.\footnote{This is explained in details right before~\eqref{eq:52}.}
In order to safely handle this term in our Lyapunov inequality, we only need to prove that it is bounded by a reasonable constant.
Here again, we split $r_0^t$ in five contributions coming from~\eqref{eq:rutMaster}:
\begin{itemize}
\item $r_1$, the part coming from the $-2\stepsize e_u$ terms;
\item $r_2$, coming from $4\lipschitz\stepsize^2 C_1 e_u$;
\item $r_3$, coming from $4\lipschitz\stepsize^2 C_1 (1 -\frac{1}{n})^{(u -\overlap)_+} \tilde e_0$;
\item $r_4$, coming from $4\lipschitz\stepsize^2 C_2\sum_{v=(u-\overlap)_+}^{u-1} e_v$;
\item $r_5$, coming from $4\lipschitz\stepsize^2 C_2\sum_{v=(u-\overlap)_+}^{u-1} (1 -\frac{1}{n})^{(v -\overlap)_+} \tilde e_0$.
\end{itemize}
Note that there is no $\tilde e_0$ in $H_t$, which is why we can safely ignore these terms here.

We have $r_1 = -2\stepsize (1 -\contraction)^t$ and $r_2=4\lipschitz\stepsize^2 C_1 (1 -\contraction)^t$.

Let us compute $r_3$.
\begin{align}
\sum_{u=0}^t (1 -\contraction)^{t -u} &(1 -\frac{1}{n})^{(u -\overlap)_+}
\nonumber \\
&= \sum_{u=0}^{\min(t, \overlap)} (1 -\contraction)^{t -u}
	+ \sum_{u=\overlap +1}^t (1 -\contraction)^{t -u} (1 -\frac{1}{n})^{u -\overlap}
\nonumber \\
&\leq (\overlap + 1) (1 -\contraction)^{t -\overlap}
	+ (1 -\contraction)^{t -\overlap} \sum_{u=\overlap +1}^t (1 -\contraction)^{\overlap -u} (1 -\frac{1}{n})^{u -\overlap}
\nonumber \\
&\leq (1 -\contraction)^t (1 -\contraction)^{-\overlap}
	\big(
		\overlap + 1 + \frac{1}{1 -q}
	\big)  \, .
\end{align}
This gives us:
\begin{align}
r_3 \leq (1 -\contraction)^t 4\lipschitz\stepsize^2 C_1 (1 -\contraction)^{-\overlap} \big(\overlap + 1 + \frac{1}{1 -q} \big)  \, .
\end{align}

We have already computed $r_4$ for $u>0$ and the computation is exactly the same for $u=0$. $r_4 \leq (1 - \contraction)^t 4\lipschitz\stepsize^2 C_2 \frac{\overlap}{1-\contraction}$ .

Finally we compute $r_5$.
\begin{align}
\sum_{u=0}^t (1 -\contraction)^{t -u} &\sum_{v=(u -\overlap)_+}^{u -1} (1 -\frac{1}{n})^{(v -\overlap)_+}
\nonumber \\
&=\sum_{v=1}^{t -1} \sum_{u=v +1}^{\min(t, v +\overlap)} (1 -\contraction)^{t -u} (1 -\frac{1}{n})^{(v -\overlap)_+}
\nonumber \\
&\leq \sum_{v=1}^{\min(t-1, \overlap)} \sum_{u=v +1}^{v +\overlap} (1 -\contraction)^{t -u}
	+ \sum_{v=\overlap + 1}^{t -1} \sum_{u=v +1}^{\min(t, v +\overlap)} (1 -\contraction)^{t -u} (1 -\frac{1}{n})^{v -\overlap}
\nonumber \\
&\leq \overlap^2 (1 -\contraction)^{t -2\overlap}
	+ \sum_{v=\overlap +1}^{t -1} (1 -\frac{1}{n})^{v -\overlap} \overlap (1 -\contraction)^{t -v -\overlap}
\nonumber \\
&\leq \overlap^2 (1 -\contraction)^{t -2\overlap}
	+ \overlap (1 -\contraction)^t (1 -\contraction)^{-2\overlap}\sum_{v=\overlap +1}^{t -1} (1 -\frac{1}{n})^{v -\overlap} \overlap (1 -\contraction)^{-v +\overlap}
\nonumber  \\
&\leq (1 -\contraction)^t (1 -\contraction)^{-2\overlap}\big(\overlap^2 + \overlap \frac{1}{1 -q}\big)  \, .
\end{align}
Which means:
\begin{align}
r_5 \leq (1 -\contraction)^t 4\lipschitz\stepsize^2 C_2 (1 -\contraction)^{-2\overlap}\big(\overlap^2 + \overlap \frac{1}{1 -q}\big) .
\end{align}
Putting it all together, we get that: $\forall t \geq 0$,
\begin{equation} \label{r0t}
\begin{aligned}
r_0^t \leq (1 -\contraction)^t
	\Big[    
	\Big( &-2 \stepsize + 
			4\lipschitz\stepsize^2 C_1+4\lipschitz\stepsize^2 C_2 \frac{\overlap}{1-\contraction} \Big) \frac{e_0}{\tilde e_0} \\
			&+4\lipschitz\stepsize^2 C_1 (1 -\contraction)^{-\overlap} \big(\overlap + 1 + \frac{1}{1 -q} \big)
		+4\lipschitz\stepsize^2 C_2 \overlap (1 -\contraction)^{-2\overlap}\big(\overlap + \frac{1}{1 -q}\big)
	\Big] .
\end{aligned}
\end{equation}

\paragraph{\textit{Sufficient condition for convergence.}}
We need all $r_u^t, u \geq 1$ to be negative so we can safely drop them from~\eqref{apx:Lyapunov}.
Note that for every $u$, this is the same condition.
We will reduce that condition to a second-order polynomial sign condition.
We also remark that since $\stepsize \geq 0$, we can upper bound our terms in $\stepsize$ and $\stepsize^2$ in this upcoming polynomial, which will give us sufficient conditions for convergence.

Now, recall that $C_2(\stepsize)$ (as defined in~\eqref{eq:C12defs}) depends on $\stepsize$. We thus need to expand it once more to find our conditions.
We have:
\begin{align*}
C_1 &= 1 + \sqrt{\sparsity}\overlap ; \qquad
C_2 = \sqrt{\sparsity} + \stepsize\strongconvex C_1  \, .
\end{align*}

Dividing the bracket in~\eqref{eq:rut} by~$\stepsize$ and rearranging as a second degree polynomial, we get the condition:
\begin{align}
4\lipschitz &\Bigg(
	C_1
	+ \frac{C_1}{n}(1-\contraction)^{-2\overlap -1} \Big[2\overlap + \frac{1}{1-q}\Big]
	+ \Big[ \frac{\sqrt{\sparsity}\overlap}{(1-\contraction)^{\overlap}} + \frac{\sqrt{\sparsity}\overlap}{n}(1 -\contraction)^{-3\overlap -1}(3\overlap + \frac{1}{1 -q}) \Big]
	\Bigg) \stepsize
\nonumber \\
	&+ 8\strongconvex C_1 \lipschitz \overlap \Big[(1- \contraction)^{-\overlap} + \frac{1}{n}(1 -\contraction)^{-3\overlap -1}(3\overlap + \frac{1}{1 -q})
	\Big] \stepsize^2
	+ 2
\leq 0  \, . \label{eq:ConvergenceCondition}
\end{align}
The discriminant of this polynomial is always positive, so $\stepsize$ needs to be between its two roots.
The smallest is negative, so the condition is not relevant to our case (where $\stepsize > 0$).
By solving analytically for the positive root~$\phi$, we get an upper bound condition on~$\stepsize$ that can be used for any overlap~$\overlap$ and guarantee convergence.
Unfortunately, for large~$\overlap$, the upper bound becomes exponentially small because of the presence of~$\overlap$ in the exponent in~\eqref{eq:ConvergenceCondition}.
More specifically, by using the bound $1/(1-\contraction) \leq \exp(2\contraction)$\footnote{This bound can be derived from the inequality $(1-x/2) \geq \exp(-x)$ which is valid for $0 \leq x \leq 1.59$.} and thus $(1-\contraction)^{-\overlap} \leq \exp(2 \overlap \contraction)$ in~\eqref{eq:ConvergenceCondition}, we would obtain factors of the form $\exp(\overlap/n)$ in the denominator for the root~$\phi$ (recall that $\contraction < 1/n$).

Our Lemma~\ref{lma:3} is derived instead under the assumption that $\overlap \leq \mathcal{O}(n)$, with the constants chosen in order to make the condition~\eqref{eq:ConvergenceCondition} more interpretable and to relate our convergence result with the standard SAGA convergence (see Theorem~\ref{th1}). As explained in Section~\ref{apxD}, the assumption that $\overlap \leq \mathcal{O}(n)$ appears reasonable in practice.
First, by using Bernoulli's inequality, we have:
\begin{equation} \label{eq:Bernouilli}
(1 - \contraction)^{k \overlap} \geq 1 - k\overlap \contraction \qquad \textnormal{for integers} \quad k\overlap \geq 0 \, .
\end{equation}
To get manageable constants, we make the following slightly more restrictive assumptions on the target rate~$\contraction$\footnote{Note that we already expected $\contraction < 1/n$.} and overlap~$\overlap$:\footnote{This bound on $\overlap$ is reasonable in practice, see Appendix~\ref{apxD}.}
\begin{align}
\contraction &\leq \frac{1}{4n} \label{eq:assumptionContraction}
\\
\overlap &\leq \frac{n}{10} \, . \label{eq:assumptionOverlap}
\end{align}
We then have:
\begin{align}
\frac{1}{1-q} &\leq \frac{4n}{3} & &
\\
\frac{1}{1-\contraction} &\leq \frac{4}{3} & &
\\
k\overlap\contraction &\leq \frac{3}{40} & &
\text{for $1 \leq k \leq 3$}
\\
(1 -\contraction)^{-k\overlap} &\leq \frac{1}{1-k\overlap \contraction} \leq \frac{40}{37} & &
\text{for $1 \leq k \leq 3$ and by using~\eqref{eq:Bernouilli}}.
\end{align}
We can now upper bound loosely the three terms in brackets appearing in~\eqref{eq:ConvergenceCondition} as follows:
\begin{align}
(1-\contraction)^{-2\overlap -1} \big[2\overlap + \frac{1}{1-q}\big] &\leq 3n \label{eq:Simp1}
\\
\sqrt{\sparsity}\overlap (1-\contraction)^{-\overlap} + \frac{\sqrt{\sparsity}\overlap}{n}(1 -\contraction)^{-3\overlap -1}(3\overlap + \frac{1}{1 -q})
\leq 4 \sqrt{\sparsity}\overlap &\leq 4C_1
\\
(1-\contraction)^{-\overlap} + \frac{1}{n}(1 -\contraction)^{-3\overlap -1}(3\overlap + \frac{1}{1 -q}) &\leq 4 \, .
\label{eq:Simp2}
\end{align}
By plugging~\eqref{eq:Simp1}--\eqref{eq:Simp2} into~\eqref{eq:ConvergenceCondition}, we get the simpler sufficient condition on~$\stepsize$:
\begin{align}
-1 + 16\lipschitz C_1 \stepsize + 16\lipschitz C_1\strongconvex\overlap\stepsize^2 \leq 0 \, .
\end{align}
The positive root $\phi$ is:
\begin{align}\label{eq:phi}
\phi = \frac{16\lipschitz C_1(\sqrt{1 + \frac{\strongconvex\overlap}{4\lipschitz C_1}} -1)}{32 \lipschitz C_1\strongconvex\overlap}
= \frac{\sqrt{1 + \frac{\strongconvex\overlap}{4\lipschitz C_1}} -1}{2\strongconvex\overlap}  \, .
\end{align}
We simplify it further by using the inequality:\footnote{This inequality can be derived by using the concavity property $f(y) \leq f(x) + (y-x) f'(x)$ on the differentiable concave function $f(x)=\sqrt{x}$ with $y=1$.}
\begin{equation}\label{eq:concavesqrt}
\sqrt{x} - 1 \geq \frac{x - 1}{2 \sqrt{x}}  \qquad \forall x > 0 \, .
\end{equation}
Using~\eqref{eq:concavesqrt} in~\eqref{eq:phi}, and recalling that $\kappa := \lipschitz / \strongconvex$, we get:
\begin{align}
\phi \geq \frac{1}{16\lipschitz C_1 \sqrt{1 + \frac{\overlap}{4\kappa C_1}}} \, .
\end{align}
Since $\frac{\overlap}{C_1} = \frac{\overlap}{1 + \sqrt{\sparsity}\overlap} \leq \min(\overlap, \frac{1}{\sqrt{\sparsity}})$, we get that a sufficient condition on our step size is:
\begin{equation}
\stepsize \leq \frac{1}{16\lipschitz (1 + \sqrt{\sparsity} \overlap) \sqrt{1 + \frac{1}{4\kappa} \min(\overlap, \frac{1}{\sqrt{\sparsity}})}}  \, .
\end{equation}

Subject to our conditions on $\stepsize$, $\contraction$ and $\overlap$, we then have that: $r_u^t \leq 0\ \text{for all}\ u\ \text{s.t.}\ 1 \leq u \leq t$.
This means we can rewrite~\eqref{apx:Lyapunov} as:
\begin{align}\label{eq:lya2}
\lyapunov_{t+1} \leq (1 - \contraction)^{t+1}a_0 + (1-\frac{\stepsize\strongconvex}{2})\lyapunov_t + r_0^t \tilde e_0  \, .
\end{align}

Now, we could finish the proof from this inequality, but it would only give us a convergence result in terms of $a_t = \E \|x_t - x^*\|^2$. A better result would be in terms of the suboptimality at $\hat x_t$ (because $\hat x_t$ is a real quantity in the algorithm whereas $x_t$ is virtual). Fortunately, to get such a result, we can easily adapt~\eqref{eq:lya2}.

We make $e_t$ appear on the left side of~\eqref{eq:lya2}, by adding $\stepsize$ to $r_t^t$ in~\eqref{apx:Lyapunov}:\footnote{We could use any multiplier from $0$ to $2\stepsize$, but choose $\stepsize$ for simplicity. For this reason and because our analysis of the $r_t^t$ term was loose, we could derive a tighter bound, but it does not change the leading terms.}
\begin{align}\label{eq:newlyapunov}
\stepsize e_t + \lyapunov_{t+1} \leq (1 - \contraction)^{t+1}a_0 + (1-\frac{\stepsize\strongconvex}{2})\lyapunov_t + \sum_{u=1}^{t-1} r_u^t e_u + r_0^t \tilde e_0 + (r_t^t + \stepsize)e_t .
\end{align}

We now require the stronger property that $\stepsize +r_t^t \leq 0$, which translates to replacing $-2\stepsize$ with $-\stepsize$ in~\eqref{eq:rut}:
\begin{equation}
\begin{aligned}
0 \geq
	\Big[
		&- \stepsize
		+4\lipschitz\stepsize^2 C_1
		+\frac{4\lipschitz\stepsize^2 C_1}{n}
			(1 - \contraction)^{-2\overlap -1} \big(
				2\overlap
				+ \frac{1}{1-q}
			\big)
\\
		&+4\lipschitz\stepsize^2  C_2 \overlap (1 - \contraction)^{-\overlap}
		+\frac{4\lipschitz\stepsize^2  C_2}{n} \overlap (1 -\contraction)^{-3\overlap}
			\big(
				3\overlap
				+ \frac{1}{1-q}
		\big)
	\Big] .
\end{aligned}
\end{equation}

We can easily derive a new stronger condition on $\stepsize$ under which we can drop all the $e_u, u > 0$ terms in~\eqref{eq:newlyapunov}:
\begin{align}
\stepsize \leq \stepsize^* = \frac{1}{32\lipschitz (1 + \sqrt{\sparsity} \overlap) \sqrt{1 + \frac{1}{8\kappa} \min(\overlap, \frac{1}{\sqrt{\sparsity}})}},
\end{align}
and thus under which we get:
\begin{align}\label{eq:lyapu3}
\stepsize  e_t + \lyapunov_{t+1} \leq (1 - \contraction)^{t+1}a_0 + (1-\frac{\stepsize\strongconvex}{2})\lyapunov_t + r_0^t \tilde e_0 .
\end{align}
This finishes the proof of Lemma~\ref{lma:3}.
\hfill\BlackBox

\subsection{Proof of Theorem~\ref{thm:convergence} (convergence guarantee and rate of \ASAGA)}\label{apxB:th2}
\paragraph{\textit{End of Lyapunov convergence.}}
We continue with the assumptions of Lemma~\ref{lma:3} which gave us~\eqref{eq:lyapu3}.
Thanks to~\eqref{r0t}, we can also rewrite $r_0^t \leq (1 -\contraction)^{t+1} A$ where $A$ is a constant which depends on $n$, $\sparsity$, $\stepsize$ and $\lipschitz$ but is finite and crucially does not depend on $t$.
In fact, by reusing similar arguments as in~\ref{apxB:lma3}, we can show the loose bound $A \leq \stepsize n$ under the assumptions of Lemma~\ref{lma:3} (including $\stepsize \leq \stepsize^*$).\footnote{In particular, note that $e_0$ does not appear in the definition of $A$ because it turns out that the parenthesis group multiplying $e_0$ in~\eqref{r0t} is negative. Indeed, it contains less positive terms than~\eqref{eq:rut} which we showed to be negative under the assumptions from Lemma~\ref{lma:3}.}
We then have:
\begin{align}
\lyapunov_{t+1}  \leq \stepsize e_t + \lyapunov_{t+1}
&\leq (1-\frac{\stepsize\strongconvex}{2})\lyapunov_t  +(1 - \contraction)^{t+1} (a_0 + A \tilde e_0)
\nonumber \\
&\leq (1-\frac{\stepsize\strongconvex}{2})^{t+1}\lyapunov_0 + (a_0 + A \tilde e_0) \sum_{k=0}^{t+1} (1 - \contraction)^{t+1 -k} (1-\frac{\stepsize\strongconvex}{2})^k .
\end{align}

We have two linearly contracting terms.
The sum contracts linearly with the minimum geometric rate factor between $\stepsize \strongconvex / 2$ and $\contraction$.
If we define $m := \min(\contraction, \stepsize\strongconvex/2)$, $M := \max(\contraction, \stepsize\strongconvex/2)$ and $\contraction^* := \nu m$ with $0 <\nu < 1$,\footnote{$\nu$ is introduced to circumvent the problematic case where $\contraction$ and $\stepsize \strongconvex / 2$ are too close together, which does not prevent the geometric convergence, but makes the constant $\frac{1}{1-\eta}$ potentially very big (in the case both terms are equal, the sum even becomes an annoying linear term in t).} we then get:\footnote{Note that if $m \neq \contraction$, we can perform the index change $t+1-k \rightarrow k$ to get the sum.}
\begin{align}
\stepsize e_t \leq
\stepsize e_t + \lyapunov_{t+1}
&\leq (1-\frac{\stepsize\strongconvex}{2})^{t+1}\lyapunov_0 + (a_0 + A \tilde e_0) \sum_{k=0}^{t+1} (1 - m)^{t+1-k} (1-M)^k
\nonumber \\
&\leq (1-\frac{\stepsize\strongconvex}{2})^{t+1}\lyapunov_0 + (a_0 + A \tilde e_0) \sum_{k=0}^{t+1} (1 - \contraction^*)^{t+1-k} (1-M)^k
\nonumber \\
&\leq (1-\frac{\stepsize\strongconvex}{2})^{t+1}\lyapunov_0 + (a_0 + A \tilde e_0) (1 - \contraction^*)^{t+1} \sum_{k=0}^{t+1} (1 - \contraction^*)^{-k} (1-M)^k
\nonumber \\
&\leq (1-\frac{\stepsize\strongconvex}{2})^{t+1}\lyapunov_0 + (1 - \contraction^*)^{t+1} \frac{1}{1 -\eta} (a_0 + A \tilde e_0)
\nonumber \\
&\leq (1 - \contraction^*)^{t+1} \big(a_0 + \frac{1}{1 -\eta} (a_0 + A \tilde e_0) \big) ,
\end{align}
where $\eta := \frac{1-M}{1-\contraction^*}$.
We have $\frac{1}{1 - \eta} = \frac{1 - \contraction^*}{M - \contraction^*}$.

By taking $\nu = \frac{4}{5}$ and setting $\contraction = \frac{1}{4n}$ -- its maximal value allowed by the assumptions of Lemma~\ref{lma:3} -- we get $M \geq \frac{1}{4n}$ and $\contraction^* \leq \frac{1}{5n}$, which means $\frac{1}{1 - \eta} \leq 20n$.
All told, using $A \leq \stepsize n$, we get:
\begin{equation}
e_t \leq (1 - \contraction^*)^{t+1} \tilde C_0,
\end{equation}
where:
\begin{equation}
\tilde C_0 := \frac{21n}{\stepsize}\Big(\|x_0 - x^*\|^2 + \stepsize \frac{n}{2L}\E\|\alpha_i^0 - f'_i(x^*)\|^2\Big) \, .
\end{equation}

Since we set $\contraction = \frac{1}{4n}, \nu = \frac{4}{5}$, we have $\nu \contraction = \frac{1}{5n}$.
Using a step size $\stepsize = \frac{a}{\lipschitz}$ as in Theorem~\ref{thm:convergence}, we get $\nu \frac{\stepsize \strongconvex}{2} = \frac{2a}{5 \kappa}$.
We thus obtain a geometric rate of $\contraction^* = \min \{\frac{1}{5n}, a\frac{2}{5\kappa}\}$, which we simplified to $\frac{1}{5} \min \{\frac{1}{n}, a\frac{1}{\kappa}\}$ in Theorem~\ref{thm:convergence}, finishing the proof.
We also observe that $\tilde C_0 \leq \frac{60n}{\stepsize} C_0$, with $C_0$ defined in Theorem~\ref{th1}.
\hfill\BlackBox

\subsection{Proof of Corollary~\ref{thm:bigdata} (speedup regimes for \ASAGA)} \label{apx:proofASAGAspeedup}
Referring to~\citet{qsaga} and our own Theorem~\ref{th1}, the geometric rate factor of \SAGA\ is $\frac{1}{5}\min\{\frac{1}{n}, \frac{a}{\kappa}\}$ for a step size of $\stepsize = \frac{a}{5L}$. We start by proving the first part of the corollary which considers the step size $\stepsize = \frac{a}{L}$ with $a = a^*(\overlap)$.
We distinguish between two regimes to study the parallel speedup our algorithm obtains and to derive a condition on $\overlap$ for which we have a linear speedup.

\paragraph{\textit{Well-conditioned regime.}}
In this regime, $n > \kappa$ and the geometric rate factor of sequential \SAGA\ is $\frac{1}{5n}$.
To get a linear speedup (up to a constant factor), we need to enforce $\contraction^* = \Omega(\frac{1}{n})$.
We recall that $\contraction^* = \min \{\frac{1}{5n}, a\frac{1}{5\kappa}\}$.

We already have $\frac{1}{5n} = \Omega(\frac{1}{n})$. This means that we need $\overlap$ to verify $\frac{a^*(\overlap)}{5\kappa} = \Omega(\frac{1}{n})$, where $a^*(\overlap) = \frac{1}{32 \left(1+ \overlap  \sqrt \sparsity \right) \xi(\kappa, \sparsity, \overlap)}$ according to Theorem~\ref{thm:convergence}.
Recall that $\xi(\kappa, \sparsity, \overlap) := \sqrt{1 + \frac{1}{8 \kappa}  \min\{\frac{1}{\sqrt{\sparsity}}, \overlap\} }$.
Up to a constant factor, this means we can give the following sufficient condition:

\begin{equation}
\frac{1}{\kappa \left(1+ \overlap  \sqrt \sparsity \right) \xi(\kappa, \sparsity, \overlap)}
= \Omega \Big(\frac{1}{n}\Big)
\end{equation}
i.e.
\begin{equation} \label{eq:SufficientCondition}
\left(1+ \overlap  \sqrt \sparsity \right) \xi(\kappa, \sparsity, \overlap)
= \mathcal{O} \Big( \frac{n}{\kappa} \Big) \, .
\end{equation}

We now consider two alternatives, depending on whether $\kappa$ is bigger than $\frac{1}{\sqrt{\sparsity}}$ or not. If $\kappa \geq \frac{1}{\sqrt{\sparsity}}$, then $\xi(\kappa, \sparsity, \overlap) < 2$ and we can rewrite the sufficient condition~\eqref{eq:SufficientCondition} as:

\begin{align}
\overlap = \mathcal{O}(1) \frac{n}{\kappa\sqrt{\sparsity}}.
\end{align}

In the alternative case, $\kappa \leq \frac{1}{\sqrt{\sparsity}}$.
Since $a^*(\overlap)$ is decreasing in $\overlap$, we can suppose $\overlap \geq \frac{1}{\sqrt{\sparsity}}$ without loss of generality and thus $\xi(\kappa, \sparsity, \overlap) = \sqrt{1 + \frac{1}{8 \kappa \sqrt{\sparsity}}}$.
We can then rewrite the sufficient condition~\eqref{eq:SufficientCondition} as:
\begin{align}\label{eq:sufcondcor3}
\frac{\overlap \sqrt{\sparsity}}{\sqrt{\kappa}\sqrt[4]{\sparsity}} &= \mathcal{O}(\frac{n}{\kappa}) \,;
\nonumber \\
\overlap &= \mathcal{O}(1)\frac{n}{\sqrt{\kappa}\sqrt[4]{\sparsity}} \, .
\end{align}

We observe that since we have supposed that $\kappa \leq \frac{1}{\sqrt{\sparsity}}$, we have $\sqrt{\kappa \sqrt{\sparsity}} \leq \kappa \sqrt{\sparsity} \leq 1$, which means that our initial assumption that $\overlap < \frac{n}{10}$ is stronger than condition~\eqref{eq:sufcondcor3}.

We can now combine both cases to get the following sufficient condition for the geometric rate factor of \ASAGA\ to be the same order as sequential \SAGA\ when $n > \kappa$:
\begin{align}
\overlap = \mathcal{O}(1) \frac{n}{\kappa\sqrt{\sparsity}}; \quad
\overlap = \mathcal{O}(n) \, .
\end{align}

\paragraph{\textit{Ill-conditioned regime.}}
In this regime, $\kappa > n$ and the geometric rate factor of sequential \SAGA\ is $a \frac{1}{\kappa}$.
Here, to obtain a linear speedup, we need $\contraction^* = \mathcal{O}(\frac{1}{\kappa})$.
Since $\frac{1}{n} > \frac{1}{\kappa}$, all we require is that $\frac{a^*(\overlap)}{\kappa} = \Omega(\frac{1}{\kappa})$ where $a^*(\overlap) = \frac{1}{32 \left(1+ \overlap  \sqrt \sparsity \right) \xi(\kappa, \sparsity, \overlap)}$, which reduces to $a^*(\overlap) = \Omega(1)$.

We can give the following sufficient condition:
\begin{align}
\frac{1}{\left(1+ \overlap  \sqrt \sparsity \right) \xi(\kappa, \sparsity, \overlap)} = \Omega(1)
\end{align}

Using that $\frac{1}{n} \leq \sparsity \leq 1$ and that $\kappa > n$, we get that $\xi(\kappa, \sparsity, \overlap) \leq 2$, which means our sufficient condition becomes:

\begin{align}
\overlap \sqrt{\sparsity} &= \mathcal{O}(1)
\nonumber \\
\overlap &= \frac{\mathcal{O}(1)}{\sqrt{\sparsity}} .
\end{align}
This finishes the proof for the first part of Corollary~\ref{thm:illcondition}.

\paragraph{\textit{Universal step size.}} 
If $\overlap = \mathcal{O}(\frac{1}{\sqrt{\sparsity}})$, then $\xi(\kappa, \sparsity, \overlap) = \mathcal{O}(1)$ and $(1+\overlap \sqrt{\sparsity}) = \mathcal{O}(1)$, and thus $a^*(\overlap) = \Omega(1)$ (for any $n$ and $\kappa$). This means that the universal step size $\stepsize = \Theta(1/L)$ satisfies $\stepsize \leq a^*(\overlap)$ for any $\kappa$, giving the same rate factor $\Omega( \min\{\frac{1}{n}, \frac{1}{\kappa}\})$ that sequential \SAGA\ has, completing the proof for the second part of Corollary~\ref{thm:illcondition}.
\hfill\BlackBox

\section{\KROMAGNON\ -- Proof of Theorem~\ref{thm:SVRG} and Corollary~\ref{thm:bigdataSVRG}}\label{apx:SVRG}
\subsection{Proof of Lemma~\ref{lma:suboptgtSVRG} (suboptimality bound on $\E \|g_t\|^2$)}\label{apx:SVRGlemma}
\citet[Lemma~9]{mania}, tells us that for serial sparse \SVRG\, we have for all $km \leq t \leq (k+1)m -1$:
\begin{equation}
\E\|g_t\|^2
\leq 2 \E \|f'_{i_t}(\hat x_t)-f'_{i_t}(x^*)\|^2
+ 2 \E \|f'_{i_t}(\hat x_k) - f'_{i_t}(x^*)\|^2 .
\end{equation}
This remains true in the case of \KROMAGNON.
We can use~\citet[Equations 7 and 8]{qsaga} to bound both terms in the following manner:
\begin{equation}
\E\|g_t\|^2
\leq 4 \lipschitz (\E f(\hat x_t)-f(x^*))
+ 4 \lipschitz (\E f(\tilde x_k) - f(x^*))
\leq 4\lipschitz e_t + 4 \lipschitz \tilde e_k \,.
\end{equation}
\hfill\BlackBox

\subsection{Proof of Theorem~\ref{thm:SVRG} (convergence guarantee and rate of \KROMAGNON)}
\paragraph{\textit{Master inequality derivation.}}\label{apx:SVRGmaster}
As in our \ASAGA\ analysis, we plug Lemma~\ref{lma:suboptgtSVRG} in Lemma~\ref{lma:1}, which gives us that for all $k \geq 0, km\leq  t \leq (k+1)m -1$:

\begin{equation}
a_{t+1}
\leq (1 - \frac{\stepsize\strongconvex}{2})a_t
+ \stepsize^2 C_1 (4\lipschitz e_t + 4 \lipschitz \tilde e_k)
+ \stepsize^2 C_2 \sum_{u = \max(km, t - \overlap)}^{t-1} (4\lipschitz e_t + 4 \lipschitz \tilde e_k)
- 2\stepsize e_t \,.
\end{equation}
By grouping the $\tilde e_k$ and the $e_t$ terms we get our master inequality~\eqref{eq:masterSVRG}:
\begin{equation*}
a_{t+1}
\leq (1 - \frac{\stepsize\strongconvex}{2})a_t
+ (4\lipschitz \stepsize^2 C_1 -2\stepsize) e_t
+ 4\lipschitz \stepsize^2 C_2 \sum_{u = \max(km, t-\overlap)}^{t-1} e_u
+ (4\lipschitz \stepsize^2 C_1 + 4 \lipschitz \stepsize^2\overlap C_2) \tilde e_k \,.
\end{equation*}
\paragraph{\textit{Contraction inequality derivation.}}\label{apx:SVRGContraction}
We now adopt the same method as in the original \SVRG\ paper~\citep{svrg}; we sum the master inequality over a whole epoch and then we use the randomization trick:
\begin{equation}\label{eq:randomtrickapx}
\tilde e_k = \E f(\tilde x_k) - f(x^*) = \frac{1}{m} \sum_{t=(k-1)m}^{km-1} e_t
\end{equation}
This gives us:
\begin{equation}\label{eq:SVRGwhatever}
\begin{aligned}
\sum_{t=km +1}^{(k+1)m} a_t
\leq
&(1 - \frac{\stepsize\strongconvex}{2})\sum_{t=km}^{(k+1)m-1} a_t
+ (4\lipschitz \stepsize^2 C_1 -2\stepsize) \sum_{t=km}^{(k+1)m-1} e_t
\\
&+ 4\lipschitz \stepsize^2 C_2 \sum_{t=km}^{(k+1)m-1} \sum_{u = \max(km, t-\overlap)}^{t-1} e_u
+ m (4\lipschitz \stepsize^2 C_1 + 4 \lipschitz \stepsize^2\overlap C_2) \tilde e_k \,.
\end{aligned}
\end{equation}
Since $1 - \frac{\stepsize\strongconvex}{2} < 1$, we can upper bound it by $1$, and then remove all the telescoping terms from~\eqref{eq:SVRGwhatever}.
We also have:
\begin{equation}
\begin{aligned}
\sum_{t=km}^{(k+1)m-1} \sum_{u = \max(km, t-\overlap)}^{t-1} e_u =
\sum_{u=km}^{(k+1)m-2} \sum_{t=u+1}^{\min((k+1)m-1, u + \overlap)} e_u &\leq
\overlap \sum_{u=km}^{(k+1)m-2} e_u
\\&\leq
\overlap \sum_{u=km}^{(k+1)m-1} e_u \,.
\end{aligned}
\end{equation}
All told:
\begin{equation}
a_{(k+1)m} \leq
a_{km}
+(4\lipschitz \stepsize^2 C_1 + 4\lipschitz \stepsize^2 \overlap C_2 -2\stepsize) \sum_{t=km}^{(k+1)m-1} e_t
+ m (4\lipschitz \stepsize^2 C_1 + 4 \lipschitz \stepsize^2\overlap C_2) \tilde e_k \,.
\end{equation}
Now we use the randomization trick~\eqref{eq:randomtrickapx}:
\begin{equation}\label{eq:SVRGrand}
a_{(k+1)m} \leq
a_{km}
+(4\lipschitz \stepsize^2 C_1 + 4\lipschitz \stepsize^2 \overlap C_2 -2\stepsize)  m \tilde e_{k+1}
+ m (4\lipschitz \stepsize^2 C_1 + 4 \lipschitz \stepsize^2\overlap C_2) \tilde e_k \,.
\end{equation}
Finally, in order to get a recursive inequality in $\tilde e_k$, we can remove the positive $a_{(k+1)m}$ term from the left-hand side of~\eqref{eq:SVRGrand} and bound the $a_{km}$ term on the right-hand side by $\nicefrac{2e_{km}}{\strongconvex}$ using a standard strong convexity inequality.
We get our final contraction inequality~\eqref{eq:svrgdetailedcontraction}:
\begin{equation*}
(2\stepsize m - 4m\lipschitz \stepsize^2 C_1 - 4m\lipschitz \stepsize^2 \overlap C_2) \tilde e_{k+1}
\leq
\big(\frac{2}{\strongconvex} + 4m\lipschitz \stepsize^2 C_1 + 4 m \lipschitz \stepsize^2\overlap C_2\big) \tilde e_k \,.
\end{equation*}
\hfill\BlackBox

\subsection{Proof of Corollary~\ref{thm:SVRGconvergence}, Corollary~\ref{thm:kromagnon} and Corollary~\ref{thm:bigdataSVRG} (speedup regimes)}

\paragraph{\textit{A simpler result for \SVRG.}}
The standard convergence rate for serial \SVRG\ is given by:
\begin{equation}
\theta := \frac{\frac{1}{\strongconvex\stepsize m} + 2 \lipschitz \stepsize}{1 - 2 \lipschitz \stepsize}\,.
\end{equation}
If we define $a$ such that $\stepsize = \nicefrac{a}{4 \lipschitz}$, we obtain:
\begin{equation}
\theta = \frac{\frac{4 \kappa}{a m} + \frac{a}{2}}{1 - \frac{a}{2}} \,.
\end{equation}
Now, since we need $\theta \leq 1$, we see that we require $a \leq 1$.
The optimal value of the denominator is then $1$ (when $a = 0$), whereas the worst case value is $\nicefrac{1}{2}$ ($a = 1$).
We can thus upper bound $\theta$ by replacing the denominator with $\nicefrac{1}{2}$, while satisfied that we do not lose more than a factor of $2$.
This gives us:
\begin{equation}
\theta \leq \frac{8 \kappa}{a m} + a \,.
\end{equation}
Enforcing $\theta \leq \nicefrac{1}{2}$ can be done easily by choosing $a \leq \nicefrac{1}{4}$ and $m = \nicefrac{32 \kappa}{a}$.
Now, to be able to compare algorithms easily, we want to frame our result in terms of rate factor per gradient computation $\contraction$, such that~\eqref{eq:svrgcontraction} is verified:
\begin{equation*}
\E f(\tilde x_k)-f(x^*) \leq (1-\rho)^{k(2m + n)} \,  (\E f(x_0) -f(x^*)) \qquad \forall k \geq 0 \,.
\end{equation*}
We define $\contraction_b := 1 - \theta$.
We want to estimate $\contraction$ such that $(1-\rho)^{2m + n} = 1 - \contraction_b$.
We get that $\contraction = 1 - (1 - \contraction_b)^{\frac{1}{2m + n}}$.
Using Bernoulli's inequality, we get:
\begin{equation}
\contraction \geq \frac{\contraction_b}{2m + n} \geq \frac{1}{2(2m + n)} \geq \frac{1}{4} \min\big\{\frac{1}{2m}, \frac{1}{n}\big\} \geq \frac{1}{4} \min\big\{\frac{a}{64\kappa}, \frac{1}{n}\big\} \, .
\end{equation}
This finishes the proof for Corollary~\ref{thm:SVRGconvergence}.
\hfill \BlackBox

\paragraph{\textit{A simpler result for \KROMAGNON.}}
We also define $a$ such that $\stepsize = \nicefrac{a}{4 \lipschitz}$.
Theorem~\ref{thm:SVRG} tells us that:
\begin{equation}
\theta = \frac{\frac{4 \kappa}{a m} + \frac{a}{2} C_3 (1 + \frac{\overlap}{16\kappa})}{1 - \frac{a}{2} C_3 (1 + \frac{\overlap}{16\kappa})}\,.
\end{equation}
We can once again upper bound $\theta$ by removing its denominator at a reasonable worst-case cost of a factor of $2$:
\begin{equation}
\theta \leq \frac{8 \kappa}{a m} + a C_3 (1 + \frac{\overlap}{16\kappa}) \,.
\end{equation}
Now, to enforce $\theta \leq \nicefrac{1}{2}$, we can choose $a \leq \frac{1}{4C_3 (1 + \frac{\overlap}{16\kappa})}$ and $m = \frac{32 \kappa}{a}$.
We also obtain a rate factor per gradient computation of: $\contraction \geq \frac{1}{4} \min\{\frac{a}{64\kappa}, \frac{1}{n}\}$.
This finishes the proof of Corollary~\ref{thm:kromagnon}.
\hfill \BlackBox

\paragraph{\textit{Speedup conditions.}}
All we have to do now is to compare the rate factors of \SVRG\ and \KROMAGNON\ in different regimes.
Note that while our convergence result hold for any $a \leq \nicefrac{1}{4}$ \SVRG\ (or the slightly more complex expression in the case of \KROMAGNON), the best step size (in terms of number of gradient computations) ensuring $\theta \leq \frac{1}{2}$ is the biggest allowable one -- thus this is the one we use for our comparison.

Suppose we are in the ``well-conditioned'' regime where $n \geq \kappa$. The rate factor of \SVRG\ is $\Omega(\nicefrac{1}{n})$.
To make sure we have a linear speedup, we need the rate factor of \KROMAGNON\ to also be $\Omega(\nicefrac{1}{n})$, which means that:
\begin{equation}\label{eq:svrgspeedupcond}
\frac{1}{256 \kappa C_3 + 16 \overlap C_3} = \Omega(\frac{1}{n})
\end{equation}
Recalling that $C_3 = 1 + 2\overlap \sqrt{\sparsity} $, we can rewrite~\eqref{eq:svrgspeedupcond} as:
\begin{equation}
\kappa = \mathcal{O}(n)\,; \quad \kappa \overlap \sqrt{\sparsity} = \mathcal{O}(n)\,; \quad \overlap = \mathcal{O}(n)\,; \quad \overlap^2 \sqrt{\sparsity} = \mathcal{O}(n) \,.
\end{equation}
We can condense these conditions down to:
\begin{equation}
\overlap = \mathcal{O}(\frac{n}{\kappa \sqrt{\sparsity}})\, ; \qquad \overlap = \mathcal{O}(\sqrt{n \sparsity^{-\nicefrac{1}{2}}}) \,.
\end{equation}
Suppose now we are in the ``ill-conditioned'' regime, where $\kappa \geq n$.
The rate factor of \SVRG\ is now $\Omega(\nicefrac{1}{\kappa})$.
To make sure we have a linear speedup, we need the rate factor of \KROMAGNON\ to also be $\Omega(\nicefrac{1}{\kappa})$, which means that:
\begin{equation}\label{eq:svrgspeedupcondIC}
\frac{1}{256 \kappa C_3 + 16 \overlap C_3} = \Omega(\frac{1}{\kappa})
\end{equation}
We can derive the following sufficient conditions:
\begin{equation}
\overlap = \mathcal{O}(\frac{1}{\sqrt{\sparsity}}) \, ; \qquad \overlap = \mathcal{O}(\kappa) \,.
\end{equation}
Since $\kappa \geq n$, we obtain the conditions of Corollary~\ref{thm:bigdataSVRG} and thus finish its proof.
\hfill \BlackBox

\section{\Hogwild\ -- Proof of Theorem~\ref{thm:convergenceSGD} and Corollary~\ref{thm:bigdataSGD}}\label{apx:SGD}
\subsection{Proof of Lemma~\ref{lma:suboptgtSGD} (suboptimality bound on $\E \|g_t\|^2$)}\label{apx:SGDlemma}
As was the case for proving Lemma~\ref{lma:suboptgt} and Lemma~\ref{lma:suboptgtSVRG}, we simply introduce $f'_i(x^*)$ in $g_t$ to derive our bound.

\begin{align}
\E \|g_t\|^2
= \E \|f'_i(\hat x_t)\|^2
&= \E \|f'_i(\hat x_t) - f'_i(x^*) + f'_i(x^*)\|^2
\nonumber
\\
&\leq 2\E \|f'_i(\hat x_t) - f'_i(x^*)\|^2 + 2\E \|f'_i(x^*)\|^2
\tag*{$(\|a + b\|^2 \leq 2\|a\|^2 + 2\|b\|^2)$}
\\
&\leq 4 \lipschitz e_t + 2 \sigma^2 \,.
\tag*{(\citet{qsaga}, Eq (7) \& (8))}
\end{align}
\hfill\BlackBox
\subsection{Proof of Theorem~\ref{thm:convergenceSGD} (convergence guarantee and rate of \Hogwild)}\label{apx:SGDtheorem}
\paragraph{\textit{Master inequality derivation.}}
Once again, we plug Lemma~\ref{lma:suboptgtSGD} into Lemma~\ref{lma:1} which gives us:
\begin{equation}
a_{t+1}
\leq (1 - \frac{\stepsize\strongconvex}{2})a_t
+ \stepsize^2 C_1 (4\lipschitz e_t + 2\sigma^2)
+ \stepsize^2 C_2 \sum_{u = (t - \overlap)_+}^{t-1} (4\lipschitz e_u + 2\sigma^2)
- 2\stepsize e_t \,.
\end{equation}
By grouping the $e_t$ and the $\sigma^2$ terms we get our master inequality~\eqref{eq:masterSGD}:
\begin{equation*}
a_{t+1}
\leq (1 - \frac{\stepsize\strongconvex}{2})a_t
+ (4\lipschitz \stepsize^2 C_1 - 2 \stepsize) e_t
+ 4\lipschitz \stepsize^2 C_2 \sum_{u = (t - \overlap)_+}^{t-1} e_u
+ 2 \stepsize^2 \sigma^2 (C_1 +  \overlap C_2)\,.
\end{equation*}

\paragraph{\textit{Contraction inequality derivation ($x_t$).}}
We now unroll~\eqref{eq:masterSGD} all the way to $t=0$ to get:
\begin{equation}\label{eq:sgdcontr}
\begin{aligned}
a_{t+1}
\leq (1 - \frac{\stepsize\strongconvex}{2})^{t+1} a_0
&+ \sum_{u = 0}^{t} (1 - \frac{\stepsize\strongconvex}{2})^{t-u} (4\lipschitz \stepsize^2 C_1 - 2 \stepsize) e_u
\\
&+ \sum_{u = 0}^{t} (1 - \frac{\stepsize\strongconvex}{2})^{t-u} 4\lipschitz \stepsize^2 C_2 \sum_{v = (u - \overlap)_+}^{u-1} e_v
\\
&+ \sum_{u = 0}^{t} (1 - \frac{\stepsize\strongconvex}{2})^{t-u} 2 \stepsize^2 \sigma^2 (C_1 +  \overlap C_2) \,.
\end{aligned}
\end{equation}
Now we can simplify these terms as follows:
\begin{align}
\sum_{u = 0}^{t} (1 - \frac{\stepsize\strongconvex}{2})^{t-u}  \sum_{v = (u - \overlap)_+}^{u-1} e_v
&=
\sum_{v = 0}^{t-1} \sum_{u = v +1}^{\min(t, v + \overlap)} (1 - \frac{\stepsize\strongconvex}{2})^{t-u} e_v
\nonumber \\
&=
\sum_{v = 0}^{t-1} (1 - \frac{\stepsize\strongconvex}{2})^{t-v} e_v \sum_{u = v +1}^{\min(t, v + \overlap)} (1 - \frac{\stepsize\strongconvex}{2})^{v-u}
\nonumber \\
&\leq \sum_{v = 0}^{t-1} (1 - \frac{\stepsize\strongconvex}{2})^{t-v} e_v \overlap (1 - \frac{\stepsize\strongconvex}{2})^{-\overlap}
\nonumber \\
&\leq \overlap (1 - \frac{\stepsize\strongconvex}{2})^{-\overlap} \sum_{v = 0}^{t} (1 - \frac{\stepsize\strongconvex}{2})^{t-v} e_v\, . \label{eq:sgdsimp1}
\end{align}
This $(1 - \frac{\stepsize\strongconvex}{2})^{-\overlap}$ term is easily bounded, as we did in~\eqref{eq:Simp1} for \ASAGA.
Using Bernoulli's inequality~\eqref{eq:Bernouilli}, we get that if we assume $\overlap \leq \frac{1}{\stepsize \strongconvex}$:\footnote{While this assumption on $\overlap$ may appear restrictive, it is in fact weaker than the condition for a linear speed-up obtained by our analysis in Corollary~\ref{thm:bigdataSGD}.}
\begin{equation}\label{eq:sgdsimp2}
(1 - \frac{\stepsize\strongconvex}{2})^{-\overlap} \leq 2 \,.
\end{equation}
We note that the last term in~\eqref{eq:sgdcontr} is a geometric sum:
\begin{align}\label{eq:sgdsimp3}
\sum_{u = 0}^{t} (1 - \frac{\stepsize\strongconvex}{2})^{t-u} \sigma^2
&= \frac{2}{\stepsize \strongconvex} \sigma^2 \,.
\end{align}
We plug~\eqref{eq:sgdsimp1}--\eqref{eq:sgdsimp3} in~\eqref{eq:sgdcontr} to obtain~\eqref{eq:sgdunroll}:
\begin{equation*}
a_{t+1} \leq (1 - \frac{\stepsize \strongconvex}{2})^{t+1} a_{0}
+ (4\lipschitz \stepsize^2 C_1 + 8\lipschitz \stepsize^2 \overlap C_2 -2\stepsize) \sum_{u=0}^{t} (1 - \frac{\stepsize \strongconvex}{2})^{t -u} e_u
+ \frac{4 \stepsize \sigma^2}{\strongconvex} ( C_1 + \overlap C_2) \,.
\end{equation*}

\paragraph{\textit{Contraction inequality derivation ($\hat x_t$).}}
We now have an contraction inequality for the convergence of $x_t$ to $x^*$.
However, since this quantity does not exist (except if we fix the number of iterations prior to running the algorithm and then wait for all iterations to be finished -- an unwieldy solution), we rather want to prove that $\hat x_t$ converges to $x^*$.
In order to do this, we use the simple following inequality:
\begin{equation}
\|\hat x_t - x^*\|^2 \leq 2 a_t + 2\|\hat x_t - x_t\|^2\,.
\end{equation}
We already have a contraction bound on the first term~\eqref{eq:sgdunroll}.
For the second term, we combine~\eqref{eq:hatxtbound} with Lemma~\ref{lma:suboptgtSGD} to get:
\begin{equation}\label{eq:thing}
\E \|\hat x_t - x_t\|^2
\leq
4 \lipschitz \stepsize^2 C_1 \sum_{u = (t - \overlap)_+}^{t-1} e_u + 2\stepsize^2 \overlap \sigma^2\,.
\end{equation}
To make it easier to combine with~\eqref{eq:sgdunroll}, we rewrite~\eqref{eq:thing} as:
\begin{equation}\label{eq:xhatxt}
\begin{aligned}
\E \|\hat x_t - x_t\|^2
&\leq
4 \lipschitz \stepsize^2 C_1 (1 - \frac{\stepsize \strongconvex}{2})^{-\overlap} \sum_{u = (t - \overlap)_+}^{t-1} (1 - \frac{\stepsize \strongconvex}{2})^{t-1-u} e_u + 2\stepsize^2 \overlap \sigma^2
\\
&\leq
8 \lipschitz \stepsize^2 C_1 \sum_{u = (t - \overlap)_+}^{t-1} (1 - \frac{\stepsize \strongconvex}{2})^{t-1-u} e_u + 2\stepsize^2 \overlap \sigma^2
\\
&\leq
8 \lipschitz \stepsize^2 C_1 \sum_{u = 0}^{t-1} (1 - \frac{\stepsize \strongconvex}{2})^{t-1-u} e_u + 2\stepsize^2 \overlap \sigma^2\,.
\end{aligned}
\end{equation}
Combining~\eqref{eq:sgdunroll} and~\eqref{eq:xhatxt} gives us~\eqref{eq:sgdxhat}:
\begin{align}
\E \|\hat x_t - x^*\|^2 \leq
(1 - \frac{\stepsize \strongconvex}{2})^{t+1} 2 a_{0}
&+ (\frac{8 \stepsize (C_1 + \overlap C_2)}{\strongconvex} + 4 \stepsize^2 C_1 \overlap) \sigma^2
\nonumber \\
&+ (24\lipschitz \stepsize^2 C_1 + 16\lipschitz \stepsize^2 \overlap C_2 -4\stepsize) \sum_{u=0}^{t} (1 - \frac{\stepsize \strongconvex}{2})^{t -u} e_u\,. \nonumber
\end{align}

\paragraph{\textit{Maximum step size condition on $\stepsize$.}}
To prove Theorem~\ref{thm:convergenceSGD}, we need an inequality of the following type: $\E \|\hat x_t - x_t\|^2 \leq (1-\contraction)^t a + b$.
To give this form to Equation~\eqref{eq:sgdxhat}, we need to remove all the $(e_u, u < t)$ terms from its right-hand side.
To safely do so, we need to enforce that all these terms are negative, hence that:
\begin{equation}
24\lipschitz \stepsize^2 C_1 + 16\lipschitz \stepsize^2 \overlap C_2 -4\stepsize \leq 0 \,.
\end{equation}
Plugging the values of $C_1$ and $C_2$ we get:
\begin{equation}
4\lipschitz \strongconvex \overlap (1 + \sqrt{\sparsity} \overlap) \stepsize^2
+ 6\lipschitz (1 + 2\sqrt{\sparsity} \overlap) \stepsize
-1 \leq 0 \,.
\end{equation}
As in our \ASAGA\ analysis, this reduces to a second-order polynomial sign condition.
We remark again that since $\stepsize \geq 0$, we can upper bound our terms in $\stepsize$ and $\stepsize^2$ in this polynomial, which will still give us sufficient conditions for convergence.
This means if we define $C_3:=1+2\sqrt{\sparsity} \overlap$, a sufficient condition is:
\begin{equation}
4\lipschitz \strongconvex \overlap C_3 \stepsize^2
+ 6\lipschitz C_3 \stepsize
-1 \leq 0\,.
\end{equation}
The discriminant of this polynomial is always positive, so $\stepsize$ needs to be between its two roots.
The smallest is negative, so the condition is not relevant to our case (where $\stepsize > 0$).
By solving analytically for the positive root~$\phi$, we get an upper bound condition on~$\stepsize$ that can be used for any overlap~$\overlap$ and guarantee convergence.
This positive root is:
\begin{equation}
\phi = \frac{3}{4} \frac{\sqrt{1 + \frac{\strongconvex \overlap}{2 \lipschitz C_3}} - 1)}{\lipschitz C_3}\,.
\end{equation}
We simplify it further by using~\eqref{eq:concavesqrt}:
\begin{equation}
\phi \geq \frac{3}{16 \lipschitz C_3 \sqrt{1 + \frac{\overlap}{2\kappa C_3}}}\,.
\end{equation}
This finishes the proof for Theorem~\ref{thm:convergenceSGD}.\hfill\BlackBox

\subsection{Proof of Theorem~\ref{thm:convergenceSGDserial} (convergence result for serial \SGD) }\label{apx:sgdsgd}
In order to analyze Corollary~\ref{thm:bigdataSGD}, we need to derive the maximum allowable step size for serial \SGD.
Note that \SGD\ verifies a simpler contraction inequality than Lemma~\ref{lma:1}.
For all $t \geq 0$:
\begin{equation} \label{eq:recursiveserialsgd}
a_{t+1} \leq
(1 - \stepsize\strongconvex) a_t + \stepsize^2 \E\|g_t\|^2 - 2\stepsize e_t  \, ,
\end{equation}
Here, the contraction factor is $(1 - \stepsize\strongconvex)$ instead of $(1 - \frac{\stepsize \strongconvex}{2})$ because $\hat x_t = x_t$ so there is no need for a triangle inequality to get back $\|x_t - x^*\|^2$ from $\|\hat x_t - x^*\|^2$ after we apply our strong convexity bound in our initial recursive inequality (see Section~\ref{app:RecurDerivation}).
Lemma~\ref{lma:suboptgtSGD} also holds for serial \SGD.
By plugging it into~\eqref{eq:recursiveserialsgd}, we get:
\begin{equation}\label{eq:recsersgd}
a_{t+1}
\leq (1 - \stepsize\strongconvex)a_t
+ (4\lipschitz \stepsize^2 - 2\stepsize) e_t
+ 2 \stepsize^2 \sigma^2\,.
\end{equation}
We then unroll~\eqref{eq:recsersgd} until $t=0$ to get:
\begin{equation}
a_{t+1}
\leq (1 - \stepsize\strongconvex)^{t+1} a_0
+ (4\lipschitz \stepsize^2 - 2\stepsize) \sum_{u = 0}^t (1 - \stepsize\strongconvex)^{t-u} e_u
+ 2 \frac{\stepsize \sigma^2}{\strongconvex}\,.
\end{equation}
To get linear convergence up to a ball around the optimum, we need to remove the $(e_u, 0\leq u \leq t)$ terms from the right-hand side of the equation.
To safely do this, we need these terms to be negative, i.e. $4\lipschitz \stepsize^2 - 2\stepsize \leq 0$.
We can then trivially derive the condition on $\stepsize$ to achieve linear convergence: $\stepsize \leq \frac{1}{2 \lipschitz}$.

We see that if $\stepsize = \nicefrac{a}{\lipschitz}$ with $a \leq \nicefrac{1}{2}$, \SGD\ converges at a geometric rate of at least: $\contraction(a) = \nicefrac{a}{\kappa},$ up to a ball of radius $2 \frac{\stepsize \sigma^2}{\strongconvex}$ around the optimum.
Now, to make sure we reach $\epsilon$-accuracy, we need $\frac{2 \stepsize \sigma^2}{\strongconvex} \leq \epsilon$, i.e. $\stepsize \leq \frac{\epsilon \strongconvex}{2\sigma^2}$.
All told, in order to get linear convergence to $\epsilon$-accuracy, serial \SGD\ requires $\stepsize \leq \min \big\{\frac{1}{2\lipschitz}, \frac{\epsilon \strongconvex}{2\sigma^2}\big\}$.

\subsection{Proof of Corollary~\ref{thm:bigdataSGD} (speedup regimes for \Hogwild)}\label{apx:SGDcorollary}
The convergence rate of both \SGD\ and \Hogwild\ is directly proportional to the step size.
Thus, in order to make sure \Hogwild\ is linearly faster than \SGD\ for any reasonable step size, we need to show that the maximum allowable step size ensuring linear convergence for \Hogwild\ -- given in Theorem~\ref{thm:convergenceSGD} -- is of the same order as the one for \SGD, $\mathcal{O}(\nicefrac{1}{\lipschitz})$.
Recalling that $\stepsize = \frac{a}{\lipschitz}$, we get the following sufficient condition: $a^*(\overlap) = \mathcal{O}(1)$.

Given~\eqref{eq:conditionSGD}, the definition of $a^*(\overlap)$, we require both:
\begin{equation}
\overlap \sqrt{\sparsity} = \mathcal{O}(1) \, ; \qquad\sqrt{1 + \frac{1}{2 \kappa}  \min\{\frac{1}{\sqrt{\sparsity}}, \overlap\} } = \mathcal{O}(1)\,.
\end{equation}
This gives us the final condition on $\overlap$ for a linear speedup: $\overlap = \mathcal{O}({\min\{\frac{1}{\sqrt{\sparsity}}}, \kappa\})$.

To finish the proof of Corollary~\ref{thm:bigdataSGD}, we only have to show that under this condition, the size of the ball is of the same order regardless of the algorithm used.

Using $\stepsize \strongconvex \overlap \leq 1$ and $\overlap \leq \frac{1}{\sqrt{\sparsity}}$, we get that $(\frac{8 \stepsize (C_1 + \overlap C_2)}{\strongconvex} + 4 \stepsize^2 C_1 \overlap) \sigma^2 = \mathcal{O}( \frac{\stepsize \sigma^2}{\strongconvex})$, which finishes the proof of Corollary~\ref{thm:bigdataSGD}.
Note that these two conditions are weaker than $\overlap = \mathcal{O}({\min\{\frac{1}{\sqrt{\sparsity}}}, \kappa\})$, which allows us to get better bounds in case we want to reach $\epsilon$-accuracy with $\frac{\epsilon \strongconvex}{\sigma^2} \ll \frac{1}{2\lipschitz}$.
\hfill\BlackBox

\section{On the Difficulty of Parallel Lagged Updates}\label{apxC}\label{apx:DifficultyLagged}
In the implementation presented in~\citet{laggedsaga}, the dense part ($\bar \alpha$) of the updates is deferred.
Instead of writing dense updates, counters $c_d$ are kept for each coordinate of the parameter vector -- which represent the last time these variables were updated -- as well as the average gradient $\bar \alpha$ for each coordinate.
Then, whenever a component $[\hat x]_d$ is needed (in order to compute a new gradient), we subtract $\stepsize (t-c_d) [\bar \alpha]_d$ from it and $c_d$ is set to $t$.
It is possible to do this without modifying the algorithm because  $[\bar \alpha]_d$ only changes when $[\hat x]_d$ also does.

In the sequential setting, this results in the same iterations as performing the updates in a dense fashion, since the coordinates are only stale when they are not used.
Note that at the end of an execution all counters have to be subtracted at once to get the true final parameter vector (and to bring every $c_d$ counter to the final $t$).

In the parallel setting, several issues arise:
\begin{itemize}
\item two cores might be attempting to correct the lag at the same time.
In which case since updates are done as additions and not replacements (which is necessary to ensure that there are no overwrites), the lag might be corrected multiple times, i.e. overly corrected.
\item we would have to read and write atomically to each $[\hat x_d], c_d, [\bar \alpha]_d$ triplet, which is highly impractical.
\item we would need to have an explicit global counter, which we do not in \ASAGA\ (our global counter~$t$ being used solely for the proof).
\item in the dense setting, updates happen coordinate by coordinate.
So at time $t$ the number of $\bar \alpha$ updates a coordinate has received from a fixed past time $c_d$ is a random variable, which may differs from coordinate to coordinate.
Whereas in the lagged implementation, the multiplier is always $(t-c_d)$ which is a constant (conditional to $c_d$), which means a potentially different $\hat x_t$.
\item the trick used in~\citet{smola} for asynchronous parallel \SVRG\ does not apply here because it relies on the fact that the ``reference'' gradient term in \SVRG\ is constant throughout a whole epoch, which is not the case for \SAGA.
\end{itemize}

All these points mean both that the implementation of such a scheme in the parallel setting would be impractical, and that it would actually yields a different algorithm than the dense version, which would be even harder to analyze.
This is confirmed by~\citet{cyclades}, where the authors tried to implement a parallel version of the lagged updates scheme and had to alter the algorithm to succeed, obtaining an algorithm with suboptimal performance as a result.

\section{Additional Empirical Details} \label{apxE}
\subsection{Detailed Description of Data Sets}\label{apx:datasets}
We run our experiments on four data sets. In every case, we run logistic regression for the purpose of binary classification.

\paragraph{\textit{RCV1 ($n=697,641$, $d=47,236$).}}
The first is the Reuters Corpus Volume I (RCV1) data set~\citep{RCV1}, an archive of over 800,000 manually categorized newswire stories made available by Reuters, Ltd. for research purposes.
The associated task is a binary text categorization.

\paragraph{\textit{URL ($n=2,396,130$, $d=3,231,961$).}}
Our second data set was first introduced in~\citet{URL}. Its associated task is a binary malicious URL detection.
This data set contains more than 2 million URLs obtained at random from Yahoo's directory listing (for the ``benign'' URLs) and from a large Web mail provider (for the ``malicious'' URLs).
The benign to malicious ratio is 2.
Features include lexical information as well as metadata.
This data set was obtained from the libsvmtools project.\footnote{\url{http://www.csie.ntu.edu.tw/~cjlin/libsvmtools/datasets/binary.html}}

\paragraph{\textit{Covertype ($n=581,012$, $d=54$).}}
On our third data set, the associated task is a binary classification problem~\citep[down from 7 classes originally, following the pre-treatment of][]{Covtype}. The features are cartographic variables.
Contrarily to the first two, this is a dense data set.

\paragraph{\textit{Realsim ($n=73,218$, $d=20,958$).}}
We only use our fourth data set for non-parallel experiments and a specific compare-and-swap test.
It constitutes of UseNet articles taken from four discussion groups (simulated auto racing, simulated aviation, real autos, real aviation).

\subsection{Implementation Details}
\paragraph{\textit{Hardware.}}
All experiments were run on a Dell PowerEdge 920 machine with 4 Intel Xeon E7-4830v2 processors with 10 2.2GHz cores each and 384GB 1600 MHz RAM.

\paragraph{\textit{Software.}} \label{scalavsc}
All algorithms were implemented in the Scala language and the software stack consisted of a Linux operating system running Scala 2.11.7 and Java 1.6.

We chose this expressive, high level language for our experimentation despite its typical 20x slower performance compared to C because our primary concern was that the code may easily be reused and extended for research purposes (which is harder to achieve with low level, heavily optimized C code; especially for error prone parallel computing).

As a result our timed experiments exhibit sub-optimal running times, e.g. compared to~\citet{s2gd}.
This is as we expected.
The observed slowdown is both consistent across data sets (roughly 20x) and with other papers that use Scala code (e.g.~\citealt{mania}, ~\citealt[Fig. 2]{cocoa}).

Despite this slowdown, our experiments show state-of-the-art results in convergence per number of iterations.
Furthermore, the speed-up patterns that we observe for our implementation of Hogwild and Kromagnon are similar to the ones given in~\citet{mania,hogwild,smola} (in various languages).

The code we used to run all the experiments is available at \url{http://www.di.ens.fr/sierra/research/asaga/}.

\begin{figure}
\center \includegraphics[width=0.5 \linewidth]{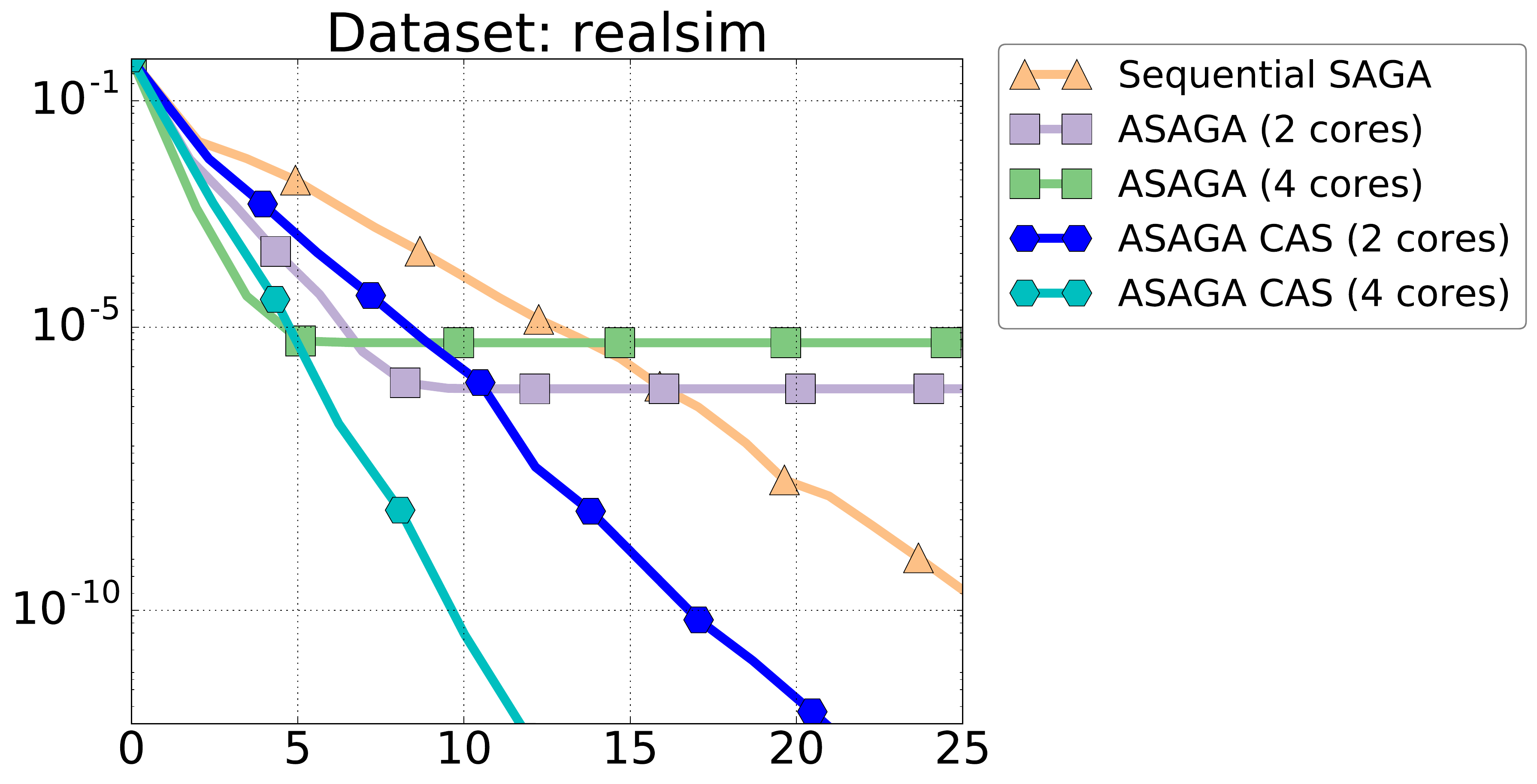}
\caption{{\bf Compare and swap in the implementation of \ASAGA}.
Suboptimality as a function of time for \ASAGA, both using compare-and-swap (CAS) operations and using standard operations.
The graph reveals that CAS is indeed needed in a practical implementation to ensure convergence to a high precision.} \label{fig:cas_comparison}
\end{figure}

\paragraph{\textit{Necessity of compare-and-swap operations.}}
Interestingly, we have found necessary to use compare-and-swap instructions in the implementation of \ASAGA.
In Figure~\ref{fig:cas_comparison}, we display suboptimality plots using non-thread safe operations and compare-and-swap (CAS) operations. The non-thread safe version starts faster but then fails to converge beyond a specific level of suboptimality, while the compare-and-swap version does converges linearly up to machine precision.

For \textit{compare-and-swap} instructions we used the \texttt{AtomicDoubleArray} class from the Google library \texttt{Guava}. This class uses an \texttt{AtomicLongArray} under the hood (from package \texttt{java.util.concurrent.atomic} in the standard Java library), which does indeed benefit from lower-level CPU-optimized instructions.

\paragraph{\textit{Efficient storage of the $\alpha_i$.}}
Storing $n$ gradient may seem like an expensive proposition, but for linear predictor models, one can actually store a single scalar per gradient~\citep[as proposed in][]{laggedsaga}, which is what we do in our implementation of \ASAGA.

\subsection{Biased Update in the Implementation}\label{apx:Bias}
In the implementation detailed in Algorithm~\ref{alg:sagasync}, $\bar \alpha$ is maintained in memory instead of being recomputed for every iteration.
This saves both the cost of reading every data point for each iteration and of computing $\bar \alpha$ for each iteration.

However, this removes the unbiasedness guarantee.
The problem here is the definition of the expectation of $\hat \alpha_i$.
Since we are sampling uniformly at random, the average of the $\hat \alpha_i$ is taken at the precise moment when we read the $\alpha_i^t$ components.
Without synchronization, between two reads to a single coordinate in $\alpha_i$ and in $\bar \alpha$, new updates might arrive in $\bar \alpha$ that are not yet taken into account in $\alpha_i$.
Conversely, writes to a component of $\alpha_i$ might precede the corresponding write in $\bar \alpha$ and induce another source of bias.

In order to alleviate this issue, we can use coordinate-level locks on $\alpha_i$ and $\bar \alpha$ to make sure they are always synchronized.
Such low-level locks are quite inexpensive when $d$ is large, especially when compared to vector-wide locks.

However, as previously noted, experimental results indicate that this fix is not necessary.

\bibliography{asaga}

\end{document}